\documentclass[a4paper,oneside,12pt]{article}
\usepackage{natbib}

\usepackage[USenglish]{babel}
\usepackage[normalem]{ulem}
\usepackage{amsthm}
\usepackage{caption}
\usepackage{subcaption}
\usepackage{booktabs}
\usepackage[flushleft]{threeparttable}
\usepackage{threeparttablex}
\usepackage[toc,page]{appendix}
\usepackage{amsmath}
\usepackage{amssymb}
\usepackage{graphicx}
\usepackage{mathtools}
\usepackage{rotating}
\usepackage{float}
\usepackage[section]{placeins}
\usepackage[colorlinks=true, allcolors=blue]{hyperref}
\usepackage{titlesec}
\usepackage[justification=raggedright,singlelinecheck=false]{caption}
\usepackage[dvipsnames]{xcolor}

\usepackage{appendix}

\titleclass{\subsubsubsection}{straight}[\subsubsection]
\newcounter{subsubsubsection}[subsubsection]
\renewcommand\thesubsubsubsection%
  {\thesubsubsection.\arabic{subsubsubsection}}

\titleformat{\subsubsubsection}
  {\normalfont\normalsize\bfseries}
  {\thesubsubsubsection}{1em}{}

\titlespacing*{\subsubsubsection}
  {0pt}{3.25ex plus 1ex minus .2ex}{1.5ex plus .2ex}

\title{\vspace*{-2cm} Estimating ``Realized'' Skewness using Convolutional Neural Network }

\author{Haozhe Jiang,$^a$ Ostap Okhrin,$^{b}$ Michael Rockinger,$^{c,*}$}

\renewcommand{\baselinestretch}{1.5} 

\theoremstyle{plain}

\newtheorem{definition}{Definition}

\newtheorem{proposition}{Proposition}
\newtheorem{theorem}{Theorem}

\theoremstyle{remark}
\newtheorem{remark}{Remark}

\newcommand{\E}{\mathbb{E}}
\newcommand{\R}{\mathbb{R}}
\newcommand{\Var}{\operatorname{Var}}
\newcommand{\Cov}{\operatorname{Cov}}

\renewcommand{\P}{\mathbb{P}}

\date{October 2026}

\begin{document}

\maketitle

\begin{abstract}
We propose a new estimator of low-frequency skewness (monthly) that exploits high-frequency data (5 min) through a direct functional mapping consisting of layers of Convolutional Neural Networks followed by Multi Layer Perceptrons. Monte Carlo experiments involving stochastic volatility models, with and without jumps, show that the estimator reduces finite-sample bias relative to existing realized-skewness estimators and remains accurate under model misspecification. Empirically, our estimator exhibits temporal stability and delivers cross-sectional performance in actual skewness-sorted portfolios. We characterize the finite-sample bias of conventional plug-in skewness estimators and establish risk consistency properties of the proposed NN estimator under standard regularity conditions.

\vspace*{0.07cm}

\vspace*{0.03cm}
\end{abstract}

\noindent$\mbox{}^{a}\,${\footnotesize 
Technische Universität Dresden. ``Friedrich List'' Faculty of Transportation. Chair of Applied Statistics. E-mail: \url{haozhe.jiang@tu-dresden.de}.
}

\noindent$\mbox{}^{b}\,${\footnotesize 
Technische Universität Dresden. ``Friedrich List'' Faculty of Transportation. Chair of Applied Statistics.  \\ Center for Scalable Data Analytics and Artificial Intelligence (ScaDS.AI) Dresden/Leipzig, Germany. E-mail: \url{ostap.okhrin@tu-dresden.de}.
}

\noindent$\mbox{}^{c}\,${\footnotesize 
University of Lausanne. Faculty of
Business and Economics, CH 1015 Lausanne, Switzerland. E-mail: \url{michael.rockinger@unil.ch}.}

\vspace*{0.5cm}
{\footnotesize  \noindent * Corresponding author. The views expressed in this paper are those of the authors and are not necessarily reflective of views from their employers. We report no conflict of interest. We are most grateful to the Deutsche Forschungsgemeinschaft (DFG, German Research Foundation) who partially funded the project under number 497262150, and to the Collaborative Research Center SFB/TRR 408 (AgiMo: Data-driven agile planning for responsible mobility), who partially funded the project under number F-016273-551-6A7-1161105.}

%
\section{Introduction} \label{sec: introduction}
\cite{pearson1895} defined the skewness of a random variable $Z$ as its standardized third central moment,
\begin{eqnarray}
    \mathrm{Skew}(Z)=
        \frac{\mathbb{E}\left[(Z-\mathbb{E}[Z])^3\right]}
    {\mathbb{E}[(Z-\mathbb{E}[Z])^2]^{\frac{3}{2}}}.\label{TheSkewness}
\end{eqnarray}
The estimation of \eqref{TheSkewness} in the context of finance is however problematic. A long history of data is required. Yet over long time intervals, the estimate can be dominated by a few rare events, such as the 1987 stock market crash \citep{KimWhite2004}, and the true skewness itself may not be stable over time.
For the variance, high-frequency (HF) data resolve this dilemma: realized variance replaces a long history of low-frequency returns by many intraday returns within the same period and converges to the integrated variance
\citep{andersen2001distribution,andersen2003,barndorff2002econometric,Barndorff2002QV}. It is
tempting to use HF returns in the same way to estimate the low-frequency skewness from a short
history. This does not work as well. No path-wise infill limit exists for the third moment analogous to realized variance because of the interaction between returns and squared returns. This cross term has an infill limit but the limit itself does not identify the third moment. Even for a Brownian motion path, a realized third moment that incorporates the leverage effect converges to a random quantity rather than capturing the third moment, which is well-known as zero.\footnote{The example of the realized third moment of a Brownian motion is illustrated in the Supplementary Material Section~S4.1.}
Consequently, observing a single path ever more finely does not reveal its low-frequency skewness. 

A new definition of skewness and a new estimator has been proposed by \cite{Neuberger2012} and extended by \cite{NP2021}.
\cite{Neuberger2012} introduced aggregator functions for the moments of a return and estimates the third and higher moments by using these aggregator functions with HF data as input. For the third moment, the resulting low-frequency measure is the sum of cubed HF returns plus a leverage term capturing the interaction between price increments and variance increments. The result holds when prices follow a martingale and applies to log-price differences or simple returns. Within this information set, the aggregator is unique. \cite{NP2021} build on \cite{Neuberger2012} by using only the stock price data and by dividing the realized third moment by realized volatility raised to the power $3/2$. They call the resulting statistic `NP skewness'. 

\cite{NP2021}'s aggregator combines a small, fixed set of functions of the HF path. The question is whether projecting the skewness onto a larger class of functions of the same path yields a more accurate estimate?

We answer this question by learning the estimator rather than deriving it.  
Following the logic of the realized third moment of \citet{Neuberger2012} and \citet{NP2021}, whose defining property holds in expectation rather than path-wise, we therefore define realized skewness as the output of a measurable aggregator functional $h(\cdot)$
which maps a path of HF returns to an estimate of the global, low-frequency skewness $S^\star$.
Supervised learning is used in this paper to approximate the aggregator functional $h(\cdot)$.
We specify a data generating process (DGP) with a parameter vector $\vartheta$ and with a available closed form for the true skewness. We draw $\vartheta$ from a prior distribution, simulate one year of five-minute returns $r$ from the model, and train a network that maps $\Phi(r)=(\tilde r,\tilde r^{2},\tilde r^{3})$, the returns and their second and third powers arranged by intraday, daily and monthly position, into a prediction of $S^\star$. Under squared loss the best such predictor is the conditional expectation $\mathbb{E}[S^\star\mid\Phi]$, that is, the posterior mean of the model-implied monthly skewness given the observed path. We call the neural approximation of this object High Frequency Network (HFN). The fundamental architecture underlying HFN is the Convolutional Neural Network (CNN), which offers a flexible way to combine many lagged HF returns and their powers, with aggregation weights that are learned rather than fixed. Initially, CNNs were invented for character recognition, \cite{lecun1989backpropagation}. In the context of financial econometrics, they represent a very natural tool to combine layers of powers of returns each layer containing lagged returns. Clearly, the way that we constructed the HFN is based on some prior information namely that the aggregated skewness should be a function of returns, their squares and their cubes. This helps to render the network more efficient.

HFN is therefore a simulation-trained estimator: it is not model-free like the NP estimator, because the specialized DGP prior enters through the training law. In this research, we use the \citet{heston1993closed} model as one DGP 
since its time series generates skewed returns due to the interaction between returns and volatility (the leverage effect) and because it can be simulated efficiently \citep{andersen2008simple} and its moments are available analytically \citep{ORS2022}. Nevertheless, the Heston model is a vehicle for producing labels and not a structural assumption of the method. The core contribution of this paper is the HFN architecture: a convolutional network whose receptive fields follow the intraday-to-daily and daily-to-monthly aggregation of HF returns, and which extracts global distributional characteristics from a HF path. Nothing in the architecture or in our main Theorem~\ref{thm:consistency} is specific to the Heston model. The theorem holds for any training law that satisfies the regularity conditions provided in the Supplementary Material. Training paths can therefore be drawn from a richer mixture of semimartingale models, for example with jumps, rough or multi-factor volatility, or microstructure noise, as long as their skewness can be computed exactly or simulated with high precision. Enlarging the training law in this way moves the learned target $\E[S^\star\mid\Phi]$ towards the one relevant for a broad class of semimartingales. 

The HFN has two distinctive features. First, instead of replacing the numerator and the denominator of \ref{TheSkewness} to get a skewness estimator it directly maps the input to the target skewness. Second, it learns the weights relating the inputs instead of keeping them fixed. 
Both features matter for the comparison with NP. The first is what allows the network to reduce the finite-sample bias of a ratio of separately estimated moments \citep{NP2021}.  Doing so, has the advantage of avoiding regions where the denominator of this expression could become very small (think about very calm stock markets) and, thus, avoids introducing large outliers in the statistics.
It therefore avoids the issue that a ratio of unbiased moments is not necessarily unbiased. Actually, it takes about 20 years of data before the bias in NP statistic disappears.  
The second feature is what makes the resulting estimates smoother over time. By construction, our method operates as a filter, smoothing out such periods. 
%
To anticipate, Figure (\ref{fig:allmethods-skewness}) shows that, fixed at one-year horizon, relative to existing competing skewness estimates, our method exhibits less bias.

We have trained the HFN by using Heston's model but the HFN continues to generate better skewness under model misspecification.
The HFN is trained on Heston simulations with parameters satisfying the Feller condition. 
Interestingly, when tested on trajectories where the Feller condition is violated, thus where volatility may come very close to 0 and stay there for a while, we find that the HFN continues to outperform competing methods. 
Also \cite{KimWhite2004} show that skewness can be influenced by rare events, such as the 1987 stock market crash. See also \cite{conrad2013ex} and  \cite{ghysels2016invest} who look for alternatives to traditional skewness. As we will show, by simulating from the \cite{bates1996jumps} model which contains jumps, the HFN estimator has very good filtering properties with respect to crashes, i.e. outliers. 
Those findings suggest that the HFN has learned a functional which remains largely model-free.

Our research should therefore be useful in various strands of research where skewness appears prominently. In finance, \cite{rubinstein1973fundamental} and \cite{kraus1976skewness} demonstrate the usefulness of skewness in portfolio allocation. In \cite{harvey1999autoregressive} co-skewness, the covariance between returns of one asset and the squared return of another one  plays an essential role. Our model should help the econometrician to test theoretical models where skewness appears such as in \cite{mcqueen2004whence, barberis2008stocks,Albuquerque2019}. 
We initiate such applications by bringing the HFN estimator to actual stock market data. It produces smoother and more persistent skewness estimates than competing methods. In the cross-section, following \cite{Amaya2016}, Fama--MacBeth regressions show that HFN skewness is positively related to next-month returns at the 5\% level, whereas the coefficients on the NP, ACJV and CL estimates are small and insignificant. A skewness-sorted long/short portfolio strategy is also illustrated as one potential application. In some other research we reconsider the question if green companies have a different skewness than brown companies.

There exists a further literature directly related to the work of \cite{Neuberger2012}.  \cite{choelee} obtain an expression involving essentially the leverage effect of HF data. 
In contrast, \cite{Amaya2016} focuses on an expression containing cubed terms aiming to capture the skewness generated by jumps.
Later \cite{fukasawa} discovers that their findings can be cast in the general setting of Bell polynomials and extend their work to any higher order cumulants. 

Our theoretical convergence and approximation results, mostly relegated to the Supplementary Material, rely on a proof strategy that is standard in modern learning-based estimation and closely related to \cite{SchmidtHieber2020Annals}. By considering the structure of our neural network we can prove convergence of the empirical-risk minimizer using the Rademacher-complexity framework of \cite{BartlettMendelson2002}. Our overall pipeline, combining continuity, neural-network approximation, and risk consistency, relates most directly to \cite{SchmidtHieber2020Annals} and recent approximation results such as \cite{Yarotsky2017} and \cite{Jentzen_2021}. The robustness result further parallels Wasserstein-based generalization bounds established in \cite{EsfahaniPeyman_2018}.


The HFN comes with some limitations. In its current form, it is trained for a specific input frequency (5-minute returns) and a fixed output horizon (monthly skewness). Adapting it to other frequencies or horizons, possibly via transfer learning, is left for future research, and we regard the present paper as a proof of concept.

The structure of the paper is as follows. In Section \ref{sec:method}, we describe our CNN model. Section \ref{sec:empirical result} reports simulation results, with robustness checks, and three empirical applications to U.S. stock data. Finally, Section \ref{sec: conclusion} concludes the paper. The Supplementary Material contains the details of the HFN architecture and its training, the detailed misspecification results, a parameter-sensitivity analysis, ablation studies, the proof of the finite-sample ratio bias, and the statistical learning theory of the HFN.
\section{Methodology} \label{sec:method}
%
\subsection{Setup and notation}
Let $(\Omega,\mathcal{F},\{\mathcal{F}_t\}_{t\in[0,T]},\mathbb{P})$ be a filtered probability space supporting 
a semimartingale log-price process $X=(X_t)_{t\in[0,T]}$ which can be decomposed in the form
\begin{equation}\label{eq:semimart}
X_t 
\;=\; X_0 
\;+\; \int_0^t \mu_s\,ds 
\;+\; \int_0^t \sigma_s\,dW_s 
\;+\; \int_0^t\!\int_{\R} z\,\big(\mu^X-\nu^X\big)(ds,dz),
\end{equation}
where $W$ is a Brownian motion, $\sigma$ is c\`adl\`ag and locally bounded, 
$\mu^X$ is the jump measure with compensator $\nu^X(ds,dz)=\lambda_s(dz)\,ds$, 
and $\int_{\R} (1\wedge z^2)\,\lambda_s(dz)<\infty$.  In this decomposition the first term is the finite-variation drift, the second term is the
continuous local martingale component, and the last term collects the jumps. This decomposition has been given by \cite{protter1990}. For ease of exposition, we assume that $X_0=0$. We are interested in the higher-order moment properties of a return over long horizons $[0,T]$.\footnote{Later, we take $T$ equal to one-month horizon unless specified otherwise.} This return is defined  by $R_T:=X_T-X_0$, and its population (central) moments and skewness are defined as
\begin{eqnarray}
m_2(T):=\E\!\big[(R_T-\E[R_T])^2\big], \quad
m_3(T):=\E\!\big[(R_T-\E [R_T])^3\big],\label{eq:moments}\\
S^\star 
:= \operatorname{Skew}(R_T)
= \frac{m_3(T)}{m_2(T)^{3/2}} = \frac{\E[(R_T-\E[R_T])^3]}{\big(\E[(R_T-\E[R_T])^2]\big)^{3/2}}.\label{eq:skewness}
\end{eqnarray}
Those statistics are assumed to be finite and well defined. Moreover, we assume that the integrated variance over $[0,T]$ is not degenerate: $m_2(T)\ge v_{\min}>0$.
To keep the notation simple, in the following we assume either that $\E[R_T]=0$ (martingale assumption), or that all moments are central moments.\footnote{Empirically, this could be obtained by centering returns corresponding to a price path.}

HF data is obtained by introducing a partition of $[0,T]$ denoted
$\Pi_M = \{0 = t_0 < t_1 < \cdots < t_{M} = T\}$, with $||\Pi_M||=\max_i(t_i-t_{i-1})$ being the mesh of $\Pi_M$ and $M$ the number of HF observations within the fixed time horizon $[0,T]$. We define HF (i.e., intraday) log-returns on the partition by $r_i:=X_{t_i}-X_{t_{i-1}}$, where $i=1,\cdots,M$ and consider the situation where $\|\Pi_M\|\downarrow 0$. Summing HF returns yields the low-frequency return over the long horizon $T$, so that $R_T \;=\; \sum_{i=1}^{M} r_i$.
\subsection{Link with the Literature} \label{sec: literature}
To nest the various models from the literature, the easiest is to start from the decomposition of standard semimartingales in equation (\ref{eq:semimart}). 
Suppose $X$ and $Y$ are two semimartingale processes on $[0,T]$.
For any $t \in [0,T]$ and a given partition $\Pi_M$, let $k = k(t,M)$ be the largest index such that
$t_k \le t < t_{k+1}$. Define
\[
Q(X,Y;\Pi_M)(t)
:= \sum_{i=1}^{k}
    \bigl(X_{t_i}-X_{t_{i-1}}\bigr)\bigl(Y_{t_i}-Y_{t_{i-1}}\bigr).
\]
If the limit
\[
[X,Y]_t
:= \lim_{\|\Pi_M\|\to 0} Q(X,Y;\Pi_M)(t)
\]
exists in probability for each $t \in [0,T]$, then $[X,Y]_t$ is called the
\emph{quadratic covariation} process of $X$ and $Y$. The process $[X]_t := [X,X]_t$ is the \emph{quadratic variation} process of $X$.

The quadratic variation of $X_t$ has two components and can be written as:
\begin{equation}
[X]_t 
= \int_0^t \sigma_s^2\,\mathrm{d}s 
+ \sum_{s\le t} (\Delta X_s)^2,
\label{eq:QV}
\end{equation}
where $X_{s-} = \lim_{u\uparrow s}X_u$ denotes the left‐limit of $X$ at $s$ and $\Delta X_s := X_s - X_{s-}$ is the jump size at time $s$.
This decomposition is the foundation of the early contemporaneous research on realized volatility by \cite{andersen2001distribution} and \cite{barndorff2001non, barndorff2002econometric}.

Equation \eqref{eq:QV} provides the exact continuous--jump decomposition of quadratic variation, which serves as the path-wise conditional second moment of returns, given the volatility path. By It\^o calculus, one can easily show \citep{Barndorff2002QV}:
\begin{equation}
X_t^2 
= X_0^2 
+ [X]_t 
+ 2\int_0^t X_{s-}\,\mathrm{d}X_s.
\label{eq:DecomposeX2}
\end{equation}
Under the martingale assumption and $X_0=0$, realized variance (RV) defined by $\sum_{i=1}^{M}\bigl(X_{t_i}-X_{t_{i-1}}\bigr)^2$ serves as a consistent estimator of \emph{quadratic
variation} and an asymptotically unbiased estimator of the conditional variance of returns \citep{Barndorff2002QV}.

An analogous decomposition for higher moments is less straightfroward.
For $p>2$, the $p$-variation of a continuous martingale vanishes faster, so a
na\"ive cubic variation $\sum_{i=1}^{M}\bigl(X_{t_i}-X_{t_{i-1}}\bigr)^3$ captures \textit{only} jump components missing the continuous part of skewness entirely. \citet{Amaya2016} (ACJV) defined realized skewness as the standardized sum of cubed HF returns:
\begin{equation}\label{eq:rskewACJV}
\mathrm{RSkew}_{T}^{ACJV}
= 
\frac{\sum_{i=1}^{M}r_{i}^3}
     {\bigl(\sum_{i=1}^{M}r_{i}^2\bigr)^{3/2}}.
\end{equation}
As mentioned, there is no continuous cubic variation analogous to quadratic variation. 
Using the result of \cite{BKS2010}, the work by \cite{Amaya2016} shows that the limiting behavior of the numerator is given by:
\begin{equation}\label{eq:RM3conv}
  \mathrm{RTM}_T^{ACJV}
  \;=\;\sum_{i=1}^{M}r_{i}^3
  \;\xrightarrow{\mathbb{P}}\;
  \sum_{s\le T}(\Delta X_s)^3.
\end{equation}
It indicates that the cubic variation is trivial for diffusions: If the process $X_t$ has no jumps, then the realized third moment converges to 0.

Therefore, estimating the low frequency third moment requires a different technique. A useful trick to decompose the third moment is proposed by \cite{choelee} who introduce the quadratic covariation of $X$ and $X^2$. For any semimartingale $X_t$, using Itô's lemma
\begin{equation}
X_t^3 
= X_0^3 
+ [X, X^2]_t 
+ \int_0^t X_{s-}\,\mathrm{d}X_s^2 
+ \int_0^t X_{s-}^2\,\mathrm{d}X_s,
\label{eq:DecomposeX3}
\end{equation}
which, after rearranging, yields the canonical semimartingale decomposition
\begin{equation}
X_t^3 
= X_0^3 
+ 3\int_0^t X_{s-}\,\mathrm{d}[X]_s
+ 3\int_0^t X_{s-}^2\,\mathrm{d}X_s
+ \sum_{s\le t} (\Delta X_s)^3.
\label{eq:cubic_semimart}
\end{equation}
Equation \eqref{eq:cubic_semimart} demonstrates that the third moment of a
semimartingale path consists of two components: a continuous part arising from
the interaction between the price level $X_{s}$ and the local quadratic
variation $\mathrm{d}[X]_s$, and a jump part given by the signed cubic jump
sizes. The third term on the right hand side of Equation \eqref{eq:cubic_semimart} is identically zero under the martingale assumption.

Consequently, \cite{choelee} prove that, specialized to the stochastic volatility model (i.e., Heston) and with the martingale assumption, the linearly scaled quadratic covariation between the log-return and the squared log-return processes is an unbiased estimator of the third central moment (also known as the third cumulant) of the long-horizon log-returns:
\begin{equation}\label{eq:3rdmomentCL}
\mathbb{E}[R_T^3] = 
\frac{3}{2}\mathbb{E}[[X,X^2]_T].
\end{equation}
This correlation between return and variance innovations is the so-called leverage effect.\footnote{Because this term misses in \cite{Amaya2016} their estimator converges to 0 when there are no jump components.}
This results in the \cite{choelee} (CL) skewness estimator:
\begin{equation}\label{eq:rskewCL}
\mathrm{RSkew}_{T}^{CL}
= 
\frac{\frac{3}{2}\sum_{i=1}^{M}\bigl(X_{t_i}-X_{t_{i-1}}\bigr)\bigl(X_{t_i}^2-X_{t_{i-1}}^2\bigr)}
     {\bigl(\sum_{i=1}^{M}r_{i}^2\bigr)^{3/2}},
\end{equation}
whose numerator is an asymptotically unbiased estimator of the true third moment when the sampling frequency goes to infinity. However, a bias still exists in finite samples, and this estimator is model-based.

\cite{Neuberger2012} did not specialize the model framework and proposed a discretization of the continuous-time identity, Eq. \eqref{eq:cubic_semimart}, as a measure of the third moment. Replacing $\mathrm{d}[X]_s$ by the change in the conditional variance of $X_T$ between times $t_{i-1}$ and $t_i$, 
i.e., $\bigl(V_{t_i}-V_{t_{i-1}}\bigr)$ where $V_{t_i}:=\Var[X_T|\mathcal{F}_{t_i}]$, and approximating the integral by a Riemann sum yields his definition of the realized third moment (RTM):
\begin{equation}
\mathrm{RTM}_T^{Neuberger} 
= 3\sum_{i=1}^{M}\bigl(X_{t_i}-X_{t_{i-1}}\bigr)\bigl(V_{t_i}-V_{t_{i-1}}\bigr)
\;+\; 
\sum_{i=1}^{M}\bigl(X_{t_i}-X_{t_{i-1}}\bigr)^3,
\label{eq:RTM_discrete}
\end{equation}
so that $\E[\mathrm{RTM}_T^{Neuberger} ]=\mathbb{E}[R_T^3]$.

Based on Equation \eqref{eq:RTM_discrete}, approximating the conditional variances by the aggregation of squared intraday returns,  \citet{NP2021} (NP) introduced a new method
\begin{equation}\label{eq:3rdmomentNP}
\begin{split}
\mathbb{E}[R_T^3] \;=\; M\,\mathbb{E}[r_i^3] + 3M\,\mathbb{E}\bigl[y_i^{(1*)}\, r_i^2\bigr],
\qquad
y_i^{(1*)}\coloneqq \frac{1}{M}\sum_{u=1}^{M} (X_{t_{i-1}}-X_{t_{i-u}}),
\end{split}
\end{equation}
leading to a new definition of the realized skewness:
\begin{equation}\label{eq:rskewNP}
\mathrm{RSkew}_{T}^{NP}
\;=\;
\frac{\sum_{i=1}^{M} r_i^3 + 3\sum_{i=1}^{M} y_i^{(1*)}\, r_i^2}
     {\bigl(\sum_{i=1}^{M} r_i^2\bigr)^{3/2}}. 
\end{equation}
\citet{NP2021} show that if  $X_t$ is a martingale and generates stationary returns, then the numerator is an unbiased estimator of the true third central moment even when the sampling frequency is finite. Intuitively, the first term on the numerator captures the jumps, whereas the second term captures the impact of volatility on returns.


One can complete this literature review of models describing realized skewness by mentioning the work by \cite{fukasawa}. Those authors realize that one can relate 
\cite{Neuberger2012} to Bell's polynomials. By building on the theory behind those polynomials, they are able to further generalize the work by Neuberger to higher order cumulants. 
\subsubsection{Limitations} \label{sec: limitation}
In summary, the numerator of ACJV does not fully capture all sources of asymmetry, particularly the leverage effect, resulting in a systematic bias. The CL and NP numerators incorporate the leverage effect and are (asymptotically) unbiased for the third central moment. However, unbiased numerators and denominators do not guarantee an unbiased ratio \citep{NP2021}. Hence, according to their methods, the plug-in realized skewness exhibits a bias in estimating the true low frequency skewness as demonstrated below.

\begin{proposition}[Finite-sample ratio bias of plug-in skewness]
\label{prop:ratio-bias}

Fix a target horizon \(T\), and let
\[
    S^\star=\frac{m_3}{m_2^{3/2}},
    \qquad m_2>0,
\]
where \(m_2\) and \(m_3\) are the second and third central moments of the
low-frequency return. Let
\[
    \widehat S
    =
    \frac{\widehat m_3}{\widehat m_2^{3/2}}
\]
be a plug-in skewness estimator based on estimators \(\widehat m_2>0\) and
\(\widehat m_3\).

Even if \(\widehat m_2\) and \(\widehat m_3\) are unbiased for \(m_2\) and \(m_3\),
respectively, \(\widehat S\) is generally biased in finite samples. The leading
(second-order) bias is
\[
\frac{15}{8}m_3m_2^{-7/2}\Var(\widehat m_2)
-
\frac{3}{2}m_2^{-5/2}\Cov(\widehat m_2,\widehat m_3).
\]
Let \(\widehat m_q=n^{-1}\sum_{j=1}^n Y_{q,j}\), \(q=2,3\), where
\(Y_{q,j}\) is computed on the \(j\)-th of \(n\) non-overlapping low-frequency
estimation windows. Here \(Y_{q,j}\) denotes a realized moment of window \(j\): \(Y_{2,j}\) is the realized variance \(\sum_{i=1}^M r_i^2\), computed from the HF returns of window \(j\), appearing in the denominator of
Eqs.~\eqref{eq:rskewACJV}, \eqref{eq:rskewCL} or \eqref{eq:rskewNP}, and \(Y_{3,j}\) is the
corresponding numerator. Suppose that
\begin{itemize}
    \item[(a)] \(\{(Y_{2,j},Y_{3,j})\}_j\) is (weakly) stationary with
    \(\E[Y_{q,j}^2]<\infty\), \(q=2,3\);
    \item[(b)] the auto- and cross-covariances are summable,
    \(\sum_{k\in\mathbb Z}\big|\Cov(Y_{q,0},Y_{r,k})\big|<\infty\) for \(q,r\in\{2,3\}\).
\end{itemize}
Then the leading ratio bias is of order \(O(n^{-1})\).
\end{proposition}
 A complete proof is provided in Section~S4.2 of the Supplementary Material.

This proposition highlights that the relevant sample size for the ratio bias is the
number of non-overlapping low-frequency observation windows, not the number of intraday observations. For example, one year of data contains only about twelve
non-overlapping monthly returns. Rolling the monthly window at a five-minute
frequency creates many overlapping monthly returns, but these windows are highly
dependent and do not increase the effective sample size proportionally. Therefore,
HF sampling alone does not eliminate the finite-sample bias of plug-in
skewness ratios.

It should be noted that due to this finite-sample ratio bias the \citet{NP2021} estimator requires roughly 20 years of HF data for a single realized-skewness estimate, implicitly assuming stationarity over long periods. This assumption is very unlikely to hold in practice. With shorter samples, instability can be substantial. How unstable is an empirical question.

\begin{figure}[!htbp]
  \centering
  \begin{subfigure}[b]{0.48\textwidth}
    \includegraphics[width=\textwidth]{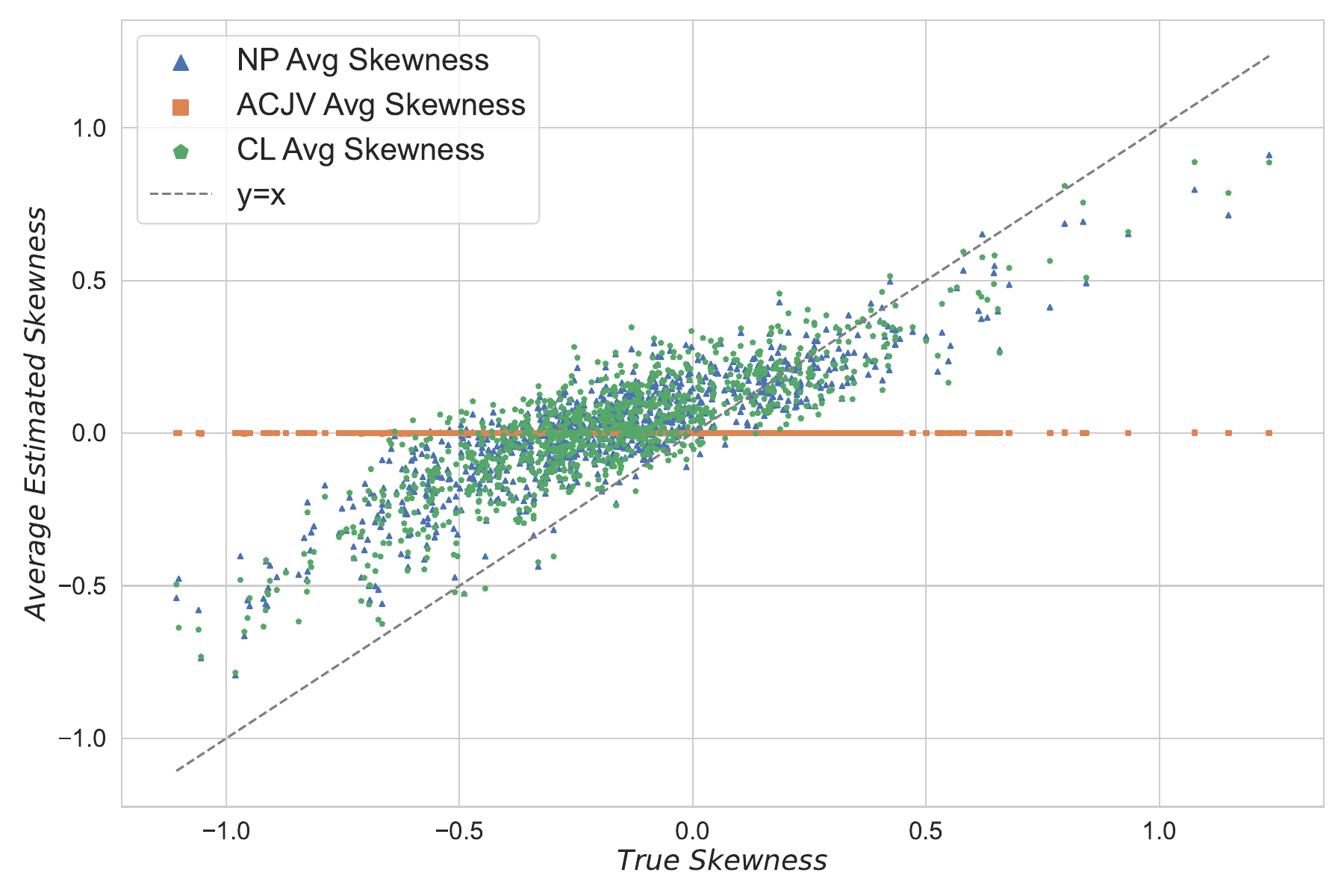}
    \caption{Average estimated skewness}
    \label{fig:othermethods-skewness}
  \end{subfigure}\hfill
  \begin{subfigure}[b]{0.48\textwidth}
    \includegraphics[width=\textwidth]{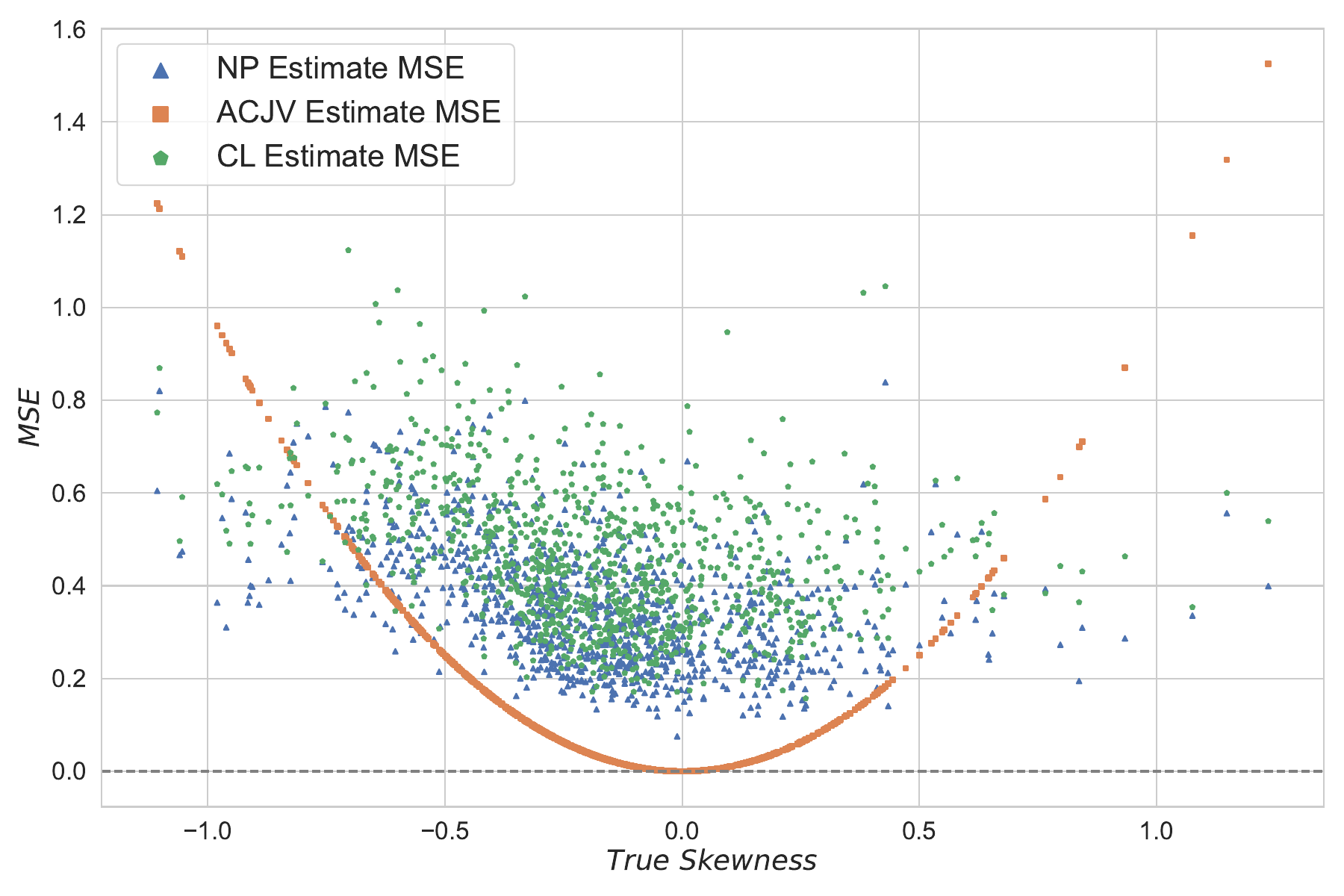}
    \caption{Mean squared error}
    \label{fig:othermethods-mse}
  \end{subfigure}
  \caption{Average estimated skewness and MSE of benchmark methods. \textit{Note:} See details of the test dataset and estimator application in the beginning of Section~\ref{sec:empirical result}.}
  \label{fig:othermethods}
\end{figure}

To investigate this issue, we simulated $50{,}000$ trajectories (for more details see Section \ref{sec:empirical result}). Figure~\ref{fig:othermethods} illustrates the performance of the benchmarks when only one year of HF data is available to estimate low-frequency skewness.
The left part of Figure \ref{fig:othermethods} displays the value of the realized skewness averaged over 50 paths as a function of the true skewness $S^\star$ determined by the formula provided in \cite{ORS2022}. The right part of this figure displays the Mean Squared Error (MSE) between the true value and the estimated value of skewness.

Because we simulate from Heston's model there are no jumps. As a consequence, as expected, the ACJV estimator converges to 0 as theory commands. On the other hand, both the CL and NP estimators are similar but exhibit an upward bias (Figure~\ref{fig:othermethods-skewness}) and are not stable as demonstrated by the MSE (Figure~\ref{fig:othermethods-mse}).

The analytical bias mechanism and empirical instability motivate learning the skewness functional directly
from HF features rather than separately estimating \(m_2\) and \(m_3\)
and forming their nonlinear ratio.
\subsubsection{Theoretical motivation: primitive-feature representation of the Bayes prediction target}
\label{subsec:theoretical-motivation}
As demonstrated above, the CL and NP benchmark estimators may exhibit non-negligible finite-sample bias.
Nevertheless, they provide important structural insight: CL links the third moment to leverage-type quadratic covariations
between $X$ and $X^2$, while NP derives discrete-time identities expressing $\E[R_T^3]$ through filtered covariances of
$r_i$ and $r_i^2$.

Building on these ideas, our goal here is not to claim an exact skewness estimator in full generality given a finite sample of size $M$ or to establish that finite-dimensional HF features are sufficient for exact recovery of the latent skewness,
but rather to formalize the prediction problem induced by economically motivated
HF representations. In particular, we observe that the primitive array $\mathcal R_M=\{(r_i,r_i^2,r_i^3)\}_{i=1}^M$ contains the information used by classical realized skewness constructions such as CL and NP, and therefore provides a natural feature basis for learning conditional skewness functionals from HF data.

\setcounter{proposition}{1}

\begin{proposition}[Optimal prediction of low-frequency skewness using high-frequency features]
\label{prop:motivation-measurable}
    Let $\Phi_M=\phi_M((\{r_i,r_i^2,r_i^3\})_{i=1}^M)$ denote the high-frequency features used as network input. Under the regularity conditions stated in Section S4.2 of the Supplementary Material, the predictor $\tilde h_M$ of $S^\star$ that minimizes expected squared error under the training law $\P_{\mathrm{tr}}$, among all predictors based on $\Phi_M$, exists as
    \begin{align*}
        \tilde h_M(\Phi_M)
        = \E[S^\star\mid\Phi_M] \qquad \text{a.s.}
    \end{align*}
\end{proposition}

The HFN is trained to approximate this conditional mean.
The formal statement and proof are provided in Section S4.2
of the Supplementary Material.
%
In contrast to ratio-based realized skewness constructions (CL/NP), the HFN directly learns
a predictive mapping from transformed HF features to low-frequency
skewness within a supervised learning framework. Under squared loss, the optimal predictor
associated with the primitive-feature representation is the conditional expectation
\(
\E[S^\star\mid \Phi_M]
\),
which depends on the observed HF data through the transformed feature tensor
\(
\Phi_M
\).
The result establishes
that economically transformed HF representations induce a well-defined
prediction target that can be approximated by a neural network architecture.

\subsection{Asymptotic Regimes} \label{sec:contribution}
Before introducing the HFN, we clarify the asymptotic regime in which we study it.
According to the aforementioned definition, realized skewness can be studied under three distinct limits:
\begin{enumerate}
    \item[(a)] \emph{Infill asymptotics}, $M\to\infty$: the number of intraday observations within the fixed horizon $[0,T]$ grows, i.e.\ $\|\Pi_M\|\to 0$;
    \item[(b)] \emph{Long-span asymptotics}, $n\to\infty$: the number of non-overlapping low-frequency windows grows, which is the regime in which the ratio bias of plug-in estimators vanishes (Proposition~\ref{prop:ratio-bias});
    \item[(c)] \emph{Training asymptotics}, $N\to\infty$: the number of simulated training pairs $\big(\Phi_M^{(j)},S^{\star(j)}\big)$ grows.
\end{enumerate}
The first paragraph in Section~\ref{sec: introduction} shows that regime (a) alone cannot deliver skewness. Regime (b) is the one on which the existing estimators implicitly rely. It requires long samples, roughly two decades of data for the NP estimator, and hence stationarity over long periods. The limitation has been discussed in Section~\ref{sec: limitation}. Our question is therefore how to obtain an unbiased and stable skewness estimate for given $M$ and $n$. Accordingly, we hold the sampling design fixed and study the training regime (c). Theorem~\ref{thm:consistency} shows that, as $N\to\infty$, the HFN attains the Bayes risk of the fixed feature representation and is unbiased on average, $\E[\widehat S_{\mathrm{HFN}}]\to\E[S^\star]$. Section~\ref{sec:empirical result} documents its stability. A joint treatment of all three regimes is left for future research.

\subsection{High Frequency Network (HFN) Approach} \label{sec:HFN approach}
According to \cite{Neuberger2012}, a realized characteristic is generated by a function $g(\cdot)$ that maps from HF returns to a low-frequency characteristic while respecting the \emph{aggregation property}, i.e., $\E[g(X_T-X_0)]=\E[ \sum_{i=1}^{M}g(X_{t_i}-X_{t_{i-1}})]$. In \eqref{eq:rskewACJV}, \eqref{eq:rskewCL}, and \eqref{eq:rskewNP}, HF observations enter as quadratic or cubic terms with constant weights in the sum. Inspired by this, a convolutional form should work on the HF data, and a non-linear functional form should be added to capture the latent non-linear feature induced by the fraction operator. A na\"ive fixed-window convolution (product–sum) over rolling blocks is essentially the ``monthly method" of \citet{NP2021}, which has a good average outcome but can have high variance due to dependence across overlapping adjacent windows. 

Convolutional neural networks (CNNs) provide therefore an interesting  alternative: they scan local fields with shared weights and, through depth, pooling, and strides, can reduce sensitivity to local shifts while capturing serial dependence \citep{lecunCNN}. We therefore design a bespoke CNN that maps HF data to a low-frequency skewness estimate, allowing, however, the aggregation weights to be learned from a large-scale simulation. Nonlinear activations in both convolutional and fully connected layers accommodate the nonlinear transformation from third moments to skewness. We call this approach the \emph{High Frequency Network} (HFN).
\subsubsection{Data Generating Process} \label{DGP}
The network weights of the HFN are learned using trajectories simulated from the underlying stochastic models. In this work, we simulate HF trajectories specifically under the jump–diffusion framework of \citet{bates1996jumps}. This way of proceeding is convening since powerful simulation methods exist, see \cite{andersen2008simple}, and there exist exact closed form formulas for the required moments, see \cite{ORS2022}. 
The following DGP specializes our general semi-martingale (\ref{eq:semimart}):
\begin{align}
	\label{SDGP}
	\mathrm{d}S_t &= \left(
	\mu - \bar k\right)S_t\, \mathrm{d}t + \sqrt{v_t} S_t\,\mathrm{d}W^x_t + S_tJ_t\mathrm{d}N_t , \\
		\label{VDGP}
	\mathrm{d}v_t &= \kappa(\theta - v_t)\mathrm{d}t + \sigma \sqrt{v_t} \mathrm{d}W^v_t, \\
	\label{corr}
	&\mathbb{E}[\mathrm{d}W^x_t \, \mathrm{d}W^v_t] = \rho \mathrm{d}t.
\end{align}
Equation~\eqref{SDGP} specifies the equity price $S_t$ with a continuous diffusion component of drift $\mu$ and stochastic variance $v_t$, augmented by multiplicative jumps of relative size $J_t$ arriving according to a homogeneous Poisson process $N_t$ with intensity $\lambda \in \mathbb{R}_+$. We assume $(1+J_t)\sim \mathrm{Lognormal}(1+\mu_J,\sigma_J^2)$ so that jumps are strictly multiplicative. Equivalently, $\log(1+J_t)\sim \mathcal{N}\!\left(\log(1+\mu_J)-\sigma_J^2,\,2\sigma_J^2\right)$. Following \citet{bates1996jumps}, the compensator $\bar{k}=\lambda \mu_J$ offsets the average jump effect so that $\mathbb{E}[\mathrm{d}S_t \mid S_t]=\mu S_t\,\mathrm{d}t$. In~\eqref{VDGP}, the instantaneous variance follows the CIR process \citep{cox1985theory}, where $\kappa$ is the speed of mean reversion, $\theta$ the long-run mean, and $\sigma$ the volatility of variance. The Brownian motions $W_t^x$ and $W_t^v$ are correlated as in~\eqref{corr}, with $\rho \in [-1,1]$.

We treat $(\mu,\kappa,\theta,\sigma)$ as constants. For the CIR dynamics, a strong solution exists for any $(\kappa,\theta,\sigma)\in(0,\infty)^3$ and $v_0>0$ \citep{RevuzYor1999}. If $\gamma \coloneqq 2\kappa\theta/\sigma^2 \ge 1$ (the Feller condition), $v_t$ remains strictly positive; when $0\le\gamma<1$, the origin is attainable and instantaneously reflecting \citep{Feller1951,PitmanYor1982}. We do not impose the Feller condition and assume only $\kappa,\theta,\sigma>0$.

For affine tractability \citep{duffie2000transform}, define $X_t=\log S_t$. By Itô’s lemma \citep{ORS2022},
\begin{align}
\label{XDGP}
\mathrm{d}X_t \;=\; \Bigl(\mu - \tfrac{1}{2} v_t - \tilde{k}\Bigr)\,\mathrm{d}t + \sqrt{v_t}\,\mathrm{d}W_t^x + \log(1+J_t)\,\mathrm{d}N_t,
\end{align}
with $\tilde{k} \coloneqq \lambda\!\left[\log(1+\mu_J)-\sigma_J^2\right]$ compensating the average jump contribution in log space. Equations \eqref{VDGP}, \eqref{corr}, and \eqref{XDGP} now jointly define the system of stochastic differential equations (SDEs) governing our data-generating process. 
Setting $\lambda=0$ yields the \citet{heston1993closed} model. We refer to specifications with jumps as the \emph{Bates model}, and to those without jumps as the \emph{Heston model}.

When we simulate $X_t$ at high frequency, we use the QE scheme of \citet{andersen2008simple} since this simulation technique has been found to have good statistical properties and is conveniently fast, see \cite{okhrin2022simulating}. 

For the training of a base HFN, we simulate trajectories of a Heston process with random parameters generated according to the distributions in the 2nd column in Table~\ref{Table: sampling distribution}. To ensure that during the training the Feller condition is verified, we generate draws for $\gamma \coloneqq 2\kappa\theta/\sigma^2$ from a shifted exponential. We also sample $(\theta,\sigma)$ from a uniform distribution with supports as indicated in the table.  Eventually, we obtain $\kappa = \gamma \sigma^2/(2\theta)$ deterministically. We fix $\mu=\theta/2$ to preserve the martingale property in logs.
Data generation and processing are described in Sections~\ref{data generation} and \ref{data processing}; the network architecture is detailed in Section~\ref{HFN}.
\subsubsection{Data Simulation: Training on the fly} \label{data generation}
Conventionally, the standard training of an NN consists of loading the entire generated dataset, splitting it into a fixed number of mini-batches, then back-propagating the gradients with respect to each batch in a certain epoch, and then iterating this gradient-descent type optimization through a fixed number of epochs. This differs from our setting, as our designed HFN targets learning from unlimited training data. To avoid memory bottlenecks, training data are generated \emph{on the fly} during training \citep{Sahoo2018}. We synthesize a fixed number (i.e., batch size) of new data blocks every time a new epoch begins and then trigger backpropagation. In other words, at each epoch, a rejuvenated batch of data is generated and submitted to the optimizer for the gradient descent. However, employing a fully refreshed data batch at every epoch introduces significant gradient noise, thereby impairing training stability and inducing catastrophic forgetting in the network.

To attenuate gradient noises, we (i) in the first epoch, initialize one dataset and repeat this set $e$ times (creating $e$ identical mini-batches). Then, we (ii) from epoch 2 onward, retain a fraction of $z\%$ of the data from the previous epoch and generate new data for $(1-z)\%$. For instance, setting $e=1$ and $z=0$ produces a single fresh batch per epoch, which empirically exhibits a ``zig-zag'' behavior in the loss curve. 

\begin{figure}[!htbp]
    \centering
    \includegraphics[width=0.75\linewidth, trim=0 180pt 0 0, clip]{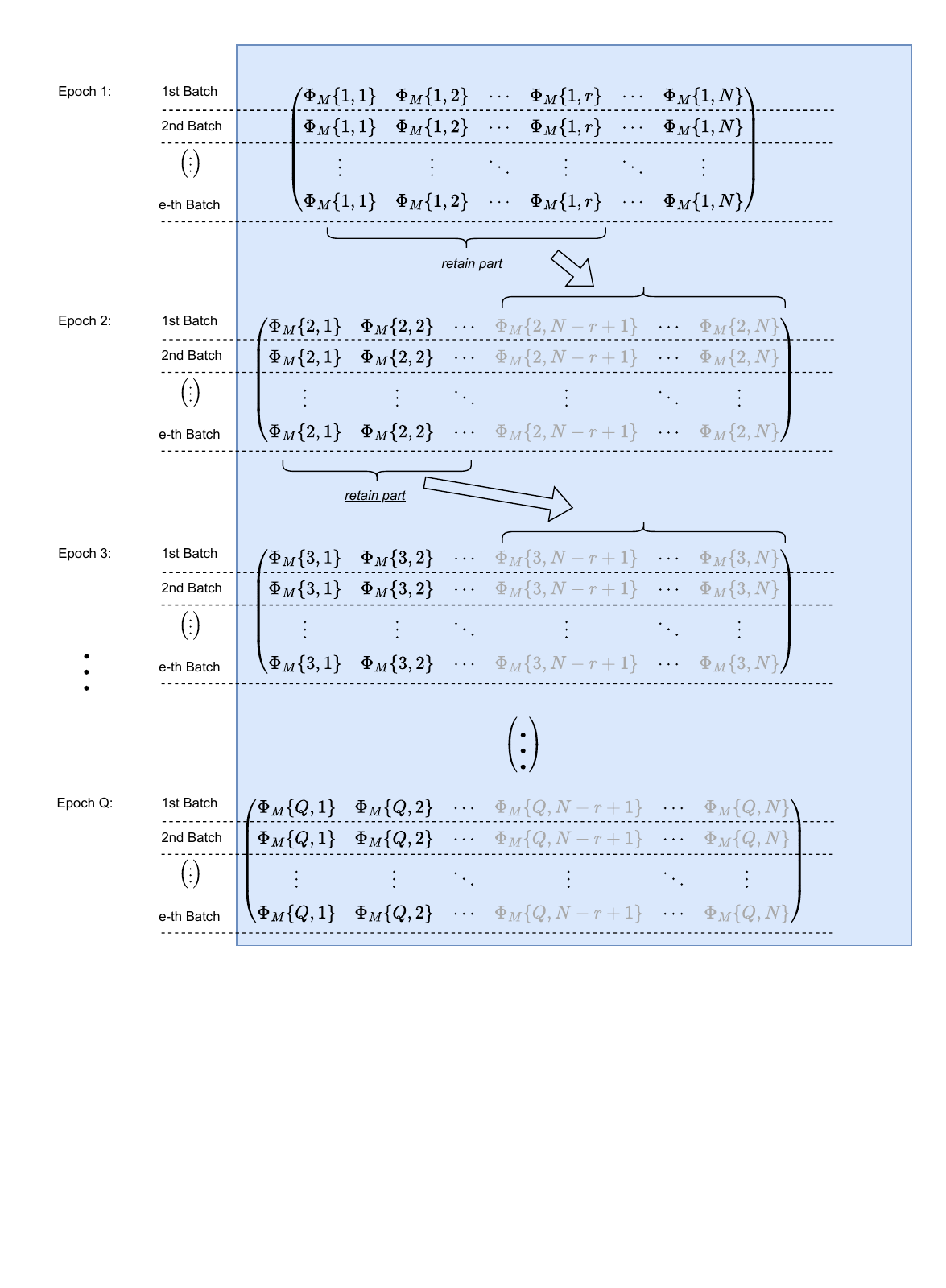}
    \caption{On-the-fly generation of training data.}
    \label{fig:DataGeneration}
\end{figure}

Figure~\ref{fig:DataGeneration} summarizes the flow of data generation, with the ``DataBlock'' $\Phi_M$ being defined in Section~\ref{data processing}.
\subsubsection{Data Processing} \label{data processing}
\cite{Liu2015} have demonstrated that 5-minute realized variance is an empirically dominant volatility estimator. For this reason, we also 
consider prices sampled at a five-minute interval running from 9{:}30 to 16{:}00 EST, for a total of $78$ intraday returns per day. We repeat this for $22$ days per month, and $12$ months per year. This yields series with $20{,}592$ observations. The objective is to estimate the monthly skewness. Other frequencies or lengths are not repeatedly discussed in this paper, as they can be handled analogously. Our focus on a monthly time horizon is justified by the fact that most allocation decisions in finance are conducted at this horizon. 

\paragraph{Reshaping time series}\label{sec:data reshape}
Having a given HF log-price time series $X_t$, we expand the single time index $t$ into a pair $(a,l_b)$ so that the log-price process can be represented as $X_{a,l_b}$ where $a=1,2,\cdots,A$ denotes the month index, $l=1,2,\cdots,L$ denotes the day index, and $b=1,2,\cdots,B$ denotes the intraday index. Following this notation, the short-horizon (i.e., HF intraday) log-return process can be defined by $r_{a,l_b} = X_{a,l_b}-X_{a,l_{b-1}}$ while the low-frequency (e.g., monthly) return is $R_{a,l_b} = X_{a,l_b}-X_{a-1,l_b}$.

\begin{figure}[!htbp]
    \centering
    \includegraphics[width=\linewidth, trim=0 0pt 0 0, clip]{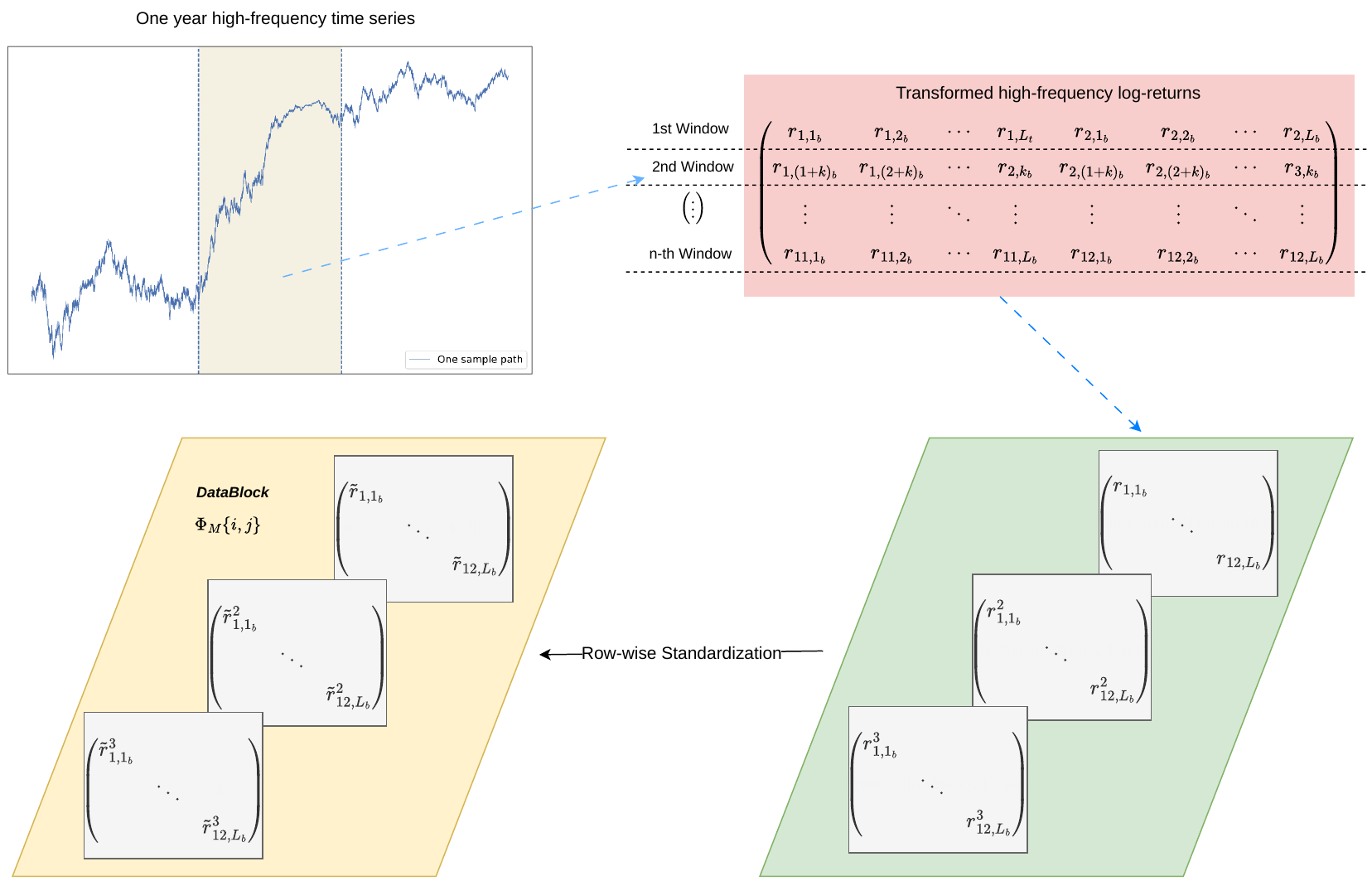}
    \caption{Transforming a high-frequency series into the HFN input.}
    \label{fig:DataTransform}
\end{figure}

As shown in Figure \ref{fig:DataTransform}, we start from the left figure with a given HF time series. The objective is to reshape this time series in a convenient way into a matrix to which we can apply the CNN filters by preserving the nature of the data. Our idea follows a rolling-window approach: we insert the first set of two-month HF log-returns in the first row. Then, we shift the window forward by $k$ days and insert the following two-month data-window into the second row. We repeat this process until the whole time series has been stored in the matrix. Here, monthly shift ($k=22$), biweekly shift ($k=10$), and weekly shift ($k=5$) can be applied with a lag in the kernel convolution, facilitating the CNN’s ability to learn lagged autocovariances at multiple granularities. This feature will be explained in section \ref{HFN}. The arrow from left to right in Figure \ref{fig:DataTransform} visualizes this data transform, reshaping a one-dimensional time series into a matrix form. Within the red box, each element $r_{a,l_b}$ contains 78 intraday returns indexed with $b=\{1,2,\cdots,78\}$. As discussed earlier, we know that the estimation of skewness also requires the second and third powers of returns. So we implement the same transform for both squared returns and cubed returns. Following the top-down arrow in Figure \ref{fig:DataTransform}, we construct a three-channel array by stacking the matrices of returns, squared returns, and cubed returns, i.e., $(r_i, r_i^2, r_i^3)$. The combined object is a tensor. The final object is similar to an RGB-type image and the tools used in machine learning for image processing can now be applied to our data.
\paragraph{Scaling}
In order to obtain good convergence of the optimization of a NN we scale the data, e.g., \citep{Ioffe2015}. As we show in Figure~\ref{fig:DataTransform} (the arrow from right to left), we perform for each row of the elements in the 3-dimensional tensor a transform dividing it by its empirical standard deviation:
\begin{equation}\label{eq:standardization}
    \tilde r^{(q)}=[\tilde{r}^{(q)}_{a,(1+k)_b},\tilde{r}^{(q)}_{a,(2+k)_b},\cdots,\tilde{r}^{(q)}_{a+2,k_b}] = \frac{[r^{(q)}_{a,(1+k)_b},r^{(q)}_{a,(2+k)_b},\cdots,r^{(q)}_{a+2,k_b}]} {\sigma([r^{(q)}_{a,(1+k)_b},r^{(q)}_{a,(2+k)_b},\cdots,r^{(q)}_{a+2,k_b}])},
\end{equation}
where $q={1,2,3}$ represents the used power of each return. At this stage, the yellow box of Figure \ref{fig:DataTransform} contains one sample of what we called in Section \ref{data generation} a ``DataBlock''. 

Thus, we can define the feature extraction function $\phi_M$, that reshapes the one-dimensional sequence of intraday returns $r_i$ into a fixed array with three channels,
\[
    \Phi_M \in \R^{H\times W\times 3},\qquad
    \Phi_M \;=\; \big[\, \tilde r,\; \tilde r^2,\; \tilde r^3\,\big],
\]
where $r$ denotes the vector of intraday returns arranged on the HF $M$ grid within the fixed time horizon $[0,T]$.  Here $H$ and $W$ depend only on the fixed low-frequency (e.g., monthly) horizon and the reshape rule  (e.g., $H=11, W=2\times M$ in the baseline). The three channels correspond to first-, second-, and third-power features of the returns. Row-wise standardization and (if needed) clipping ensure that all entries are bounded, so that the image of $\phi_M$ lies in a compact subset of $\mathcal X\subset\R^{H\times W\times 3}$, a property required in the subsequent learning-theoretic results. 

\subsubsection{HFN Architecture} \label{HFN}
We design a HFN with $c$ convolutional blocks (depth) and $F$ filters per block (width). 
In this architecture, depth governs the complexity of the approximating functional by enabling nested compositions of local operators, whereas width determines the dimensionality of the learned feature basis at each layer. This specification should be viewed as a flexible architectural family: practitioners may increase $c$ or $F$ to enrich the representational capacity and allow the network to approximate more complex nonlinear functionals. For theoretical purposes, we formalize this flexibility through the following HFN family.

\begin{definition}[HFN family]
For each \(p\ge1\), let \(\mathcal F_p\) denote a family of
high-frequency networks (HFNs) mapping the feature tensor
\[
    \Phi_M\in\mathcal X
    \subset
    \mathbb R^{H\times W\times 3}
\]
into an output vector containing a low-frequency skewness prediction.

Each \(f\in\mathcal F_p\) consists of:
(i) a finite sequence of convolutional layers operating across progressively coarser temporal aggregation scales,
(ii) pointwise non-polynomial activation functions,
(iii) pooling operators, and
(iv) a finite sequence of fully connected layers.
\end{definition}

The index \(p\) parameterizes the architectural complexity of the network,
including the number of convolutional layers, the number of convolutional
filters, and the number of hidden units in the fully connected component.

\begin{figure}[!htbp]
    \centering
    \includegraphics[width=\linewidth, trim=0 300pt 0 0, clip]{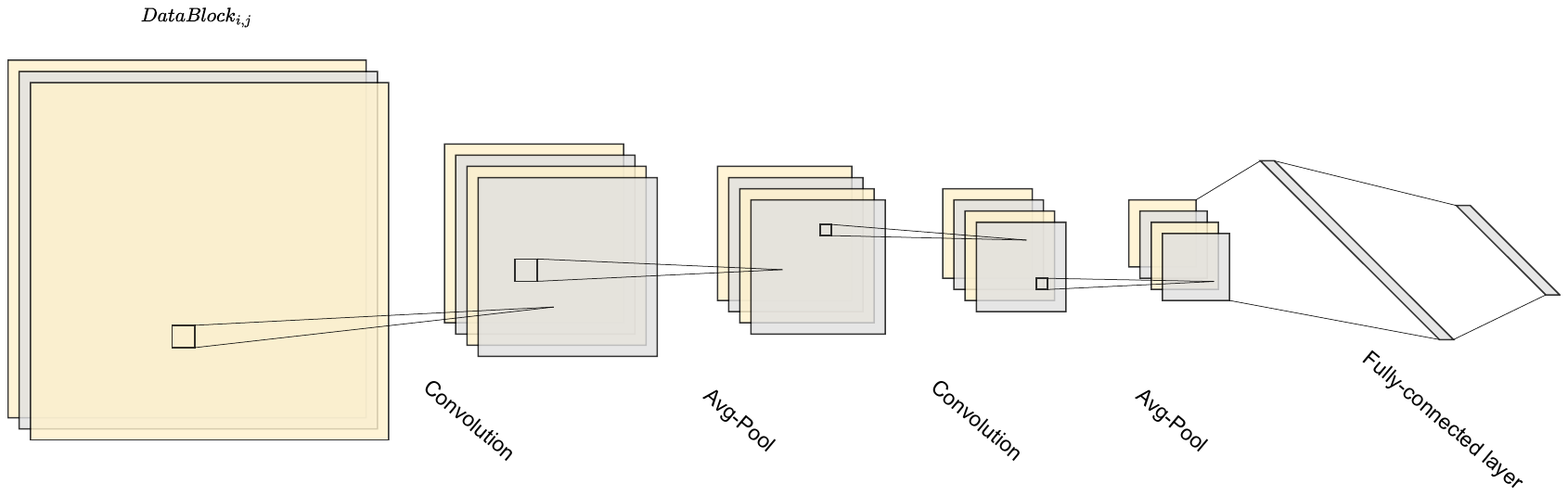}
    \caption{High Frequency Network (HFN) Structure. \textit{Note:} Our baseline HFN consists of two levels of multi-dimensional convolution followed by average pooling.}
    \label{fig:HFN}
\end{figure}

In the empirical implementation studied in this paper, which focuses on estimating monthly skewness, we use a parsimonious member of this family consisting of two convolutional blocks and four filters per block, designed to reflect the frequency-aligned aggregation structure (i.e., intraday-to-daily and daily-to-monthly) of the skewness estimation problem.  
Figure \ref{fig:HFN} depicts the structure of a two-layer baseline HFN ($c=2$).
The detailed formulas of the convolution and the average-pooling operations are provided in Section~S1 of the Supplementary Material. Suffice it to say here that a convolutional layer multiplies local subsets of the entries in a ``DataBlock'' (i.e., feature tensor $\Phi_M$) by small cubes, called ``kernels'', and then sums these products to pass the resulting feature maps to the next layer. After each layer we apply a non-linear activation function ReLU. The convolutional layer substantially reduces the number of parameters through weight sharing and can be interpreted as a learned filtering operation.

The average-pooling layer computes the average value within small patches of the feature maps. This operation further reduces the dimensionality of the data, introduces smoothing by construction, and retains only the most stable local information.

It is important to stress that in our architecture, the first convolution-pooling pair reduces the large raw-data tensor to daily observations, while the second pair aggregates the daily features to the monthly frequency used in the subsequent layers.

The final layers are fully connected layers with the ReLU activation function. The output of the NN is an estimate of skewness and an estimate of volatility. Because realized volatility and realized skewness are expected to share common latent features, we include the true volatility factor in the training target (i.e., labels) to stabilize the loss convergence. In the sequel, this empirical implementation is referred to as the
\textit{baseline HFN} and, for notational simplicity, abbreviated as
HFN in the empirical results section.

Further details on the training, such as learning rate, optimizer, etc, can be found in Section~S2 of the Supplementary Material.

\subsubsection{Theoretical properties of HFN} \label{sec:uniHFN}

The HFN is trained on simulated HF paths generated from the stochastic volatility models in Section~\ref{DGP}. Its inputs are the feature tensors $\Phi_M$ constructed on a partition $\Pi_M$ as described in Section~\ref{data processing}, and training minimizes squared loss on simulated pairs $\big(\Phi_M^{(j)},\,S^{\star(j)}\big)$, where $j=1,\cdots,N$ indexes i.i.d.\ training samples. Please note that the i.i.d.\ assumption applies to the simulated training ensemble, not to the observed market time series themselves.

For a fixed sampling design \(M\), the HFN learning problem is therefore a standard
supervised regression problem on the feature space \(\mathcal X\). As shown in Proposition~\ref{prop:motivation-measurable}, under squared loss,
the corresponding Bayes prediction target is
\[
    \E[S^\star\mid \Phi_M=x].
\]

Given a sampling resolution $\|\Pi_M\|$, the performance of the HFN estimator depends on two ingredients: (i) the number of simulated training samples $N$, and (ii) the expressive capacity of the neural network. The theoretical results below characterize approximation and prediction-risk properties of the HFN relative to this conditional expectation.
\medskip

\setcounter{theorem}{0}
\begin{theorem}[Risk consistency of the HFN estimator]
\label{thm:consistency}
Fix a sampling design \(M\) and the associated feature tensor \(\Phi_M\).
Consider an HFN trained under squared loss on simulated data
$\{(\Phi_M^{(j)},S^{\star(j)})\}_{j=1}^N$ and the resulting estimate denoted as $\widehat S_{\mathrm{HFN}}$.
Suppose the regularity assumptions in
Supplementary Material hold. Then, as
\(N\to\infty\),
\[
\E\!\Big[
\big(
\widehat S_{\mathrm{HFN}}
-
\E[S^\star\mid\Phi_M=x]
\big)^2
\Big]
\longrightarrow 0,
\]

Moreover,
\[
\E[\widehat S_{\mathrm{HFN}}]
\longrightarrow
\E[S^\star].
\]
\end{theorem}
A complete proof of the risk-consistency statement is provided in Section~S4.3 of the Supplementary Material. The second property, which shows unconditional calibration, follows immediately from Cauchy--Schwarz together with the tower
property of conditional expectations.

The theorem identifies the HFN as a risk-consistent approximation of the
Bayes predictor induced by the fixed feature representation \(\Phi_M\).
Within the chosen feature space, the HFN asymptotically recovers the
optimal predictor under squared loss. 

The unconditional calibration
statement further implies that, averaged over the training distribution,
the HFN does not systematically overestimate or underestimate skewness.
This should not be interpreted as conditional unbiasedness for each realized HF path. The pointwise prediction target remains the Bayes predictor
\(\E[S^\star\mid\Phi_M]\). Hence, while individual pathwise predictions may still contain irreducible conditional error, the estimator reproduces the correct first-order skewness level in expectation across the simulated population.
This result establishes an aggregate
calibration property that is particularly relevant for applications
involving cross-sectional ranking, portfolio sorting, or estimation of
average skewness premia.

\subsubsection{Empirical Baseline HFN}
As mentioned in Section \ref{HFN}, for our empirical study of monthly skew estimation throughout this paper, a baseline HFN is selected from the HFN family. It consists of two convolutional
blocks with kernel sizes \((3,78)\) and \((3,22)\), followed by a fully connected head.
These receptive fields are chosen to match the intraday-to-daily and daily-to-monthly
aggregation structure of the problem.
\begin{remark}[Why use a two-layer HFN?]
\label{rem:hfn-vs-generic} The theoretical results above apply to a sufficiently rich HFN family on the feature tensor $\Phi_M$. However, the selected two-layer HFN baseline is better aligned with the economics of the problem and brings a series of advantages:

\begin{enumerate}
    \item By fixing the convolutional receptive fields at the intraday and daily scales $(78,22)$, every filter in the first block ``sees'' an entire trading day, and every filter in the second block aggregates information over an entire month. This matches the structure of realized variance and realized skewness operators, which are themselves built from day-by-day and month-by-month block aggregations.

    \item Because the relevant temporal scales are hard-coded, the architecture family $\mathcal F_p$ does not have to spend capacity on discovering which lags matter. The frequency-aligned receptive fields reduce the burden on the network to learn
    economically relevant temporal aggregation scales from data alone. Empirically,
    this structured design improves training stability and predictive performance
    relative to generic architectures explored in preliminary experiments.

    \item The number of free parameters grows mainly through the numbers of filters and hidden units, while kernel supports remain fixed. This keeps the Rademacher complexity $\mathfrak R_N(\mathcal F_p)$ in Assumption 4 (Supplementary Material) under control and prevents the capacity of the network from exploding as the number of layer increases.
\end{enumerate}

From an optimization perspective, the reduced architectural search space also simplifies training: gradient-based methods operate on a lower-dimensional parameter manifold and do not need to explore configurations with inappropriate temporal receptive fields. While we do not make formal statements about the nonconvex optimization landscape, our empirical results suggest that the baseline HFN trains more stably and reaches good solutions with moderate network sizes and sample sizes.
\end{remark}

\section{Empirical Results} \label{sec:empirical result}
This section evaluates the empirical performance of the proposed HFN method relative to three benchmark approaches: \texorpdfstring{\citet{Amaya2016}'s (ACJV)}{Amaya et al.'s (ACJV)}, \texorpdfstring{\citet{choelee}'s (CL)}{Choe and Lee's (CL)}, and \texorpdfstring{\citet{NP2021}'s (NP)}{Neuberger and Payne's (NP)} methods, using both simulated and real-world data.

The test dataset comprises $1{,}000\times50$ simulated stochastic processes. Each process corresponds to one year of five-minute prices from 9:30 to 16:00 EST, producing $78$ intraday returns per day, $22$ days per month, and $12$ months per year ($20{,}592$ log-returns per path). We draw 1,000 Bates parameter sets from the predefined distributions, and for each parameter set, simulate 50 sample paths. The quantity of interest is the \emph{monthly} realized skewness, i.e., the skewness of 22-day returns.

For each simulated path, realized skewness is derived by reshaping the time series into the prescribed input tensor and applying the trained HFN operator. The benchmark estimators are computed according to equations \eqref{eq:rskewACJV}, \eqref{eq:rskewCL}, and \eqref{eq:rskewNP}, scaled to monthly quantity. The experiments proceed in five stages.

First, a baseline HFN is trained under the Heston framework (without jumps) using the hyperparameters reported in Table~S1 in Section~S2 of the Supplementary Material. Once trained, this model is used for all comparative analyses and interpretability studies. In this section we test the base HFN under the Heston framework (without jumps), namely, we sample the DGP (\ref{XDGP}) model parameters as in the 2nd column in Table \ref{Table: sampling distribution}. The parameter ranges are anchored in economically meaningful priors and selected to be sufficiently broad so as to include the calibration outcomes of \cite{eraker2004}.  Importantly, although the testing environment mirrors the training setup, the random seeds used to generate test data differ, ensuring that results are out-of-sample. The findings are presented in Section~\ref{result: no jump}.

Second, the base HFN is evaluated under two misspecified scenarios. The first includes normally distributed jumps while maintaining the Feller condition; the second removes jumps but allows the Feller condition to be violated. Parameter distributions for both scenarios are summarized in the 3rd and 4th columns in Table~\ref{Table: sampling distribution}. The range of jump parameters in the 3rd column follows \cite{eraker2004}. The results of the skewness prediction are summarized in Section~\ref{result: misspecification} and reported in detail in Section~S3.1 of the Supplementary Material.
\begin{table}[!htbp]
\centering
\caption{Sampling distributions for Bates' model parameters in 3 use cases}
\small 
\begin{tabular}{l|l|l|l}
\textbf{Param.} & \textbf{Base Test} & \textbf{Misspecification Test 1} & \textbf{Misspecification Test 2} \\
\hline
$\gamma=\frac{2\kappa\theta}{\sigma^2}$ & $\mathrm{Exp}(4)+1$ & $\mathrm{Exp}(4)+1$ & $\mathrm{Exp}(0.5)+10^{-6}$ \\
$\theta$ & $\mathcal{U}(0.05,0.5)$ & $\mathcal{U}(0.05,0.5)$ & $\mathcal{U}(0.05,0.5)$ \\
$\sigma$ & $\mathcal{U}(0.1,1)$ & $\mathcal{U}(0.1,1)$ & $\mathcal{U}(0.1,1)$ \\
$\rho$ & $\mathcal{U}(-1,1)$ & $\mathcal{U}(-1,1)$ & $\mathcal{U}(-1,1)$ \\
$\mu$ & $\mu=\theta/2$ & $\mu=\theta/2$ & $\mu=\theta/2$ \\
$\kappa$ & $\gamma\sigma^2/(2\theta)$ & $\gamma\sigma^2/(2\theta)$ & $\gamma\sigma^2/(2\theta)$ \\
$\lambda$ & $0$ & $\mathcal{U}(0.25,1.50)$ & $0$ \\
$\mu_J$ & $0$ & $\mathcal{U}(-0.11,0.01)$ & $0$ \\
$\sigma_J$ & $0$ & $\mathcal{U}(0.04,0.12)$ & $0$ \\
\end{tabular}
\label{Table: sampling distribution}
\end{table}
Third, we explore the interpretability of the HFN by examining its estimated dependence on key Heston parameters ($\rho$, $\kappa$, $\sigma$). For each parameter set, we compute the true skewness and average the estimated skewness across simulated paths to assess whether the HFN learns the correct functional relationships (Section~S3.2 of the Supplementary Material).

Fourth, a series of ablation studies investigates the determinants of HFN performance, focusing on (i) hyperparameter and initialization effects, (ii) input data format, and (iii) network architecture (Section~S3.3 of the Supplementary Material).

Finally, we apply all estimators to HF U.S.\ equity data from 1998–2025 to assess empirical behavior, predictive power, and cross-sectional implications (Section~\ref{result: real-world}).

\subsection{Comparison in the Jump-Free Case} \label{result: no jump}
\begin{figure}[!htbp]
  \centering
  \begin{subfigure}[b]{0.48\textwidth}
    \includegraphics[width=\textwidth]{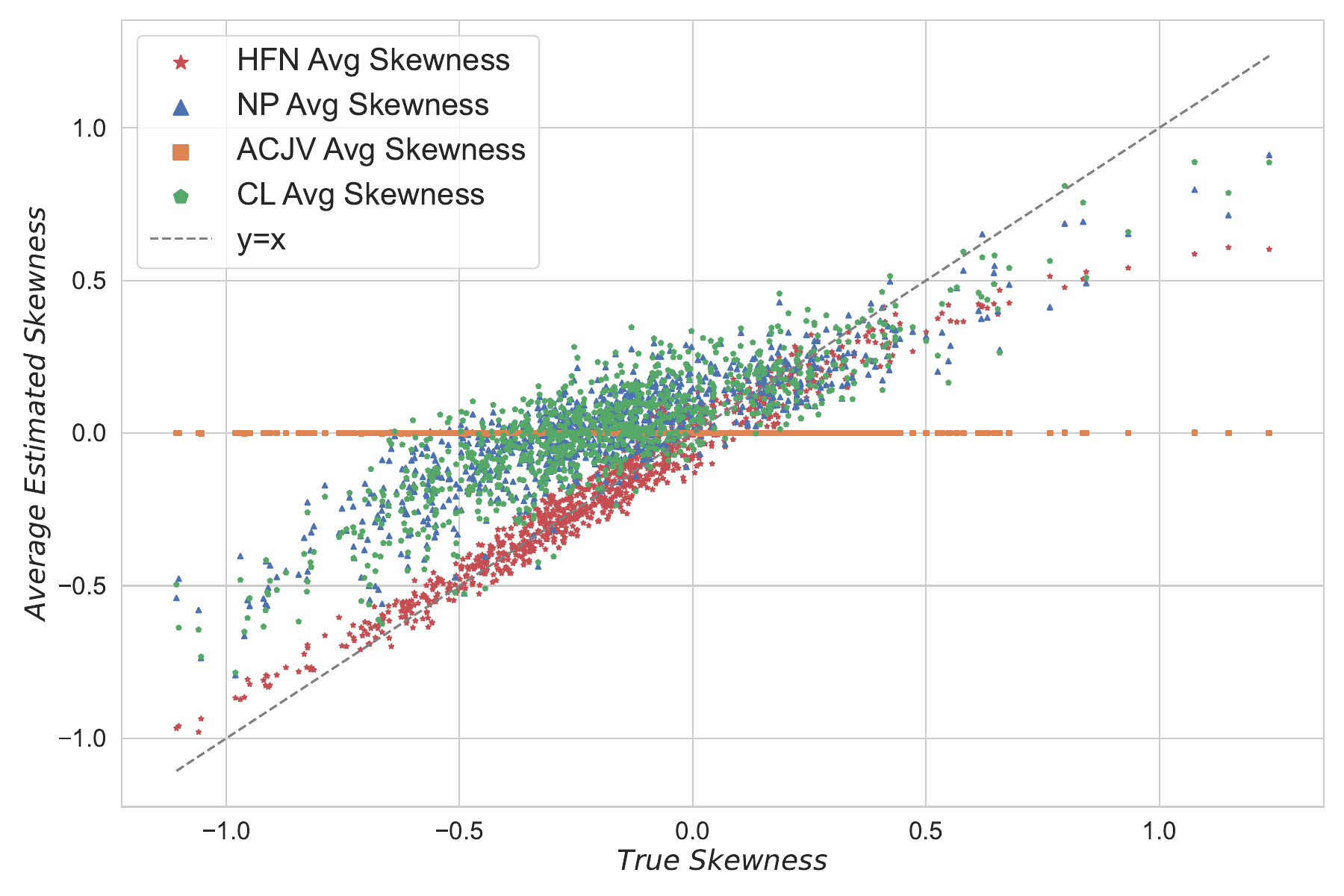}
    \caption{Average estimated skewness}
    \label{fig:allmethods-skewness}
  \end{subfigure}\hfill
  \begin{subfigure}[b]{0.48\textwidth}
    \includegraphics[width=\textwidth]{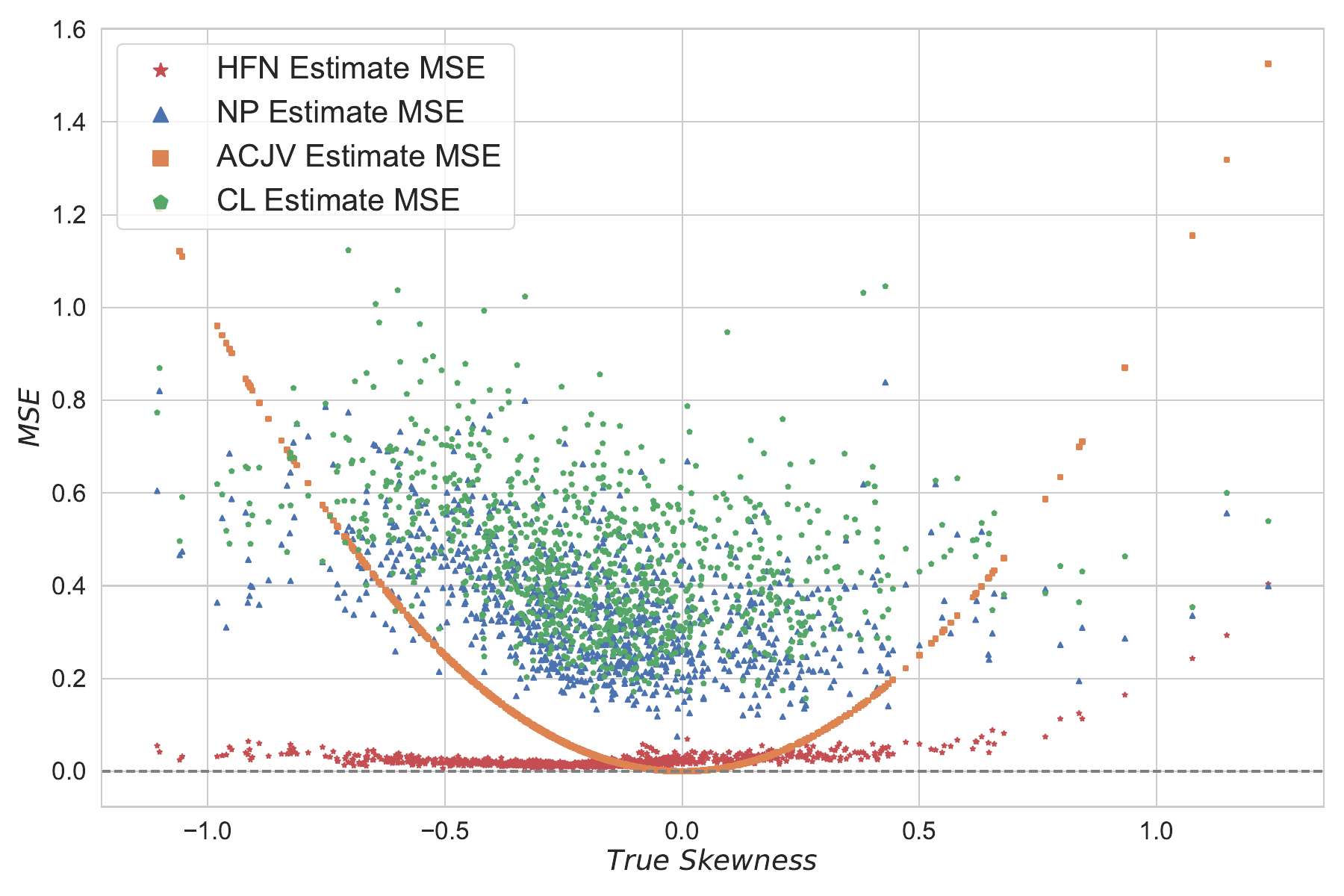}
    \caption{Mean squared error}
    \label{fig:allmethods-mse}
  \end{subfigure}
  \caption{Average estimated skewness and MSE of the HFN and benchmark estimators (no jumps).}
  \label{fig:allmethods}
\end{figure}

In a Heston test environment, as Figure \ref{fig:allmethods} demonstrates, our HFN approach generates the best skewness estimator. The base HFN's realized skewness results are, on average, scattered closely around the true skewness line. Moreover, the mean squared error of the HFN outcomes exhibits no sensitivity to variations in the true skewness level, remaining almost uniformly low across the parameter space. Hence, it not only reduces the upward bias of the CL and NP estimators, but also enhances the estimation stability as visualized in Figure \ref{fig:allmethods-mse}. Especially for a random variable with a skewness close to -1 or even lower than -1, which fits the fact of negatively skewed financial market returns \citep{Neuberger2012}, the base HFN method produces accurate and robust estimated coefficients, with a substantial advantage over all other methods. In contrast, while both the NP and CL methods capture the qualitative pattern of the true skewness, their estimates are systematically shifted upward relative to the true values, resulting in a substantial positive bias. In addition, the high MSEs suggests unstable estimates. The observed pattern in Figure \ref{fig:allmethods-skewness} implies a high likelihood that the NP and CL methods overestimate the skewness of the financial return data, adversely affecting investors who are particularly averse to left-tail events. The ACJV omits the leverage effect, which has been empirically evidenced as a major driver of the skewness over long horizons \citep{NP2021}, and converges to 0 as theory predicts. For high skewness levels, say above 1, our HFN estimator becomes biased, and its robustness also deteriorates. Nonetheless, for stock market data, returns are not in this high skewness range. Thus, the bias should not be of practical concern in this context.

\subsection{Misspecification Tests} \label{result: misspecification}
We also examine the effect of moderate distributional shifts between the simulated training law and the evaluation law. Since the HFN predictors are Lipschitz on a
compact feature domain, their squared prediction risk satisfies a standard
bounded-Lipschitz stability bound under perturbations of the joint law of
\((\Phi_M,S^\star)\).

\setcounter{theorem}{1}
\begin{remark}[Distributional stability of prediction risk]\label{thm:robust} 
    Let $\P_0$ denote the true law of $(\Phi_M,S^\star)$ and $\P_{\mathrm{tr}}$ the training law. Assuming the Lipschitz and boundedness conditions (Assumption 3 in Supplementary Material), with $|S^\star|\le S_{\max}$ almost surely and 
    $\mathcal X$ compact, the squared prediction risk of any predictor \(f\in\mathcal F_p\) satisfies
    \[
        \Big|\E_{\P_0}\big[(\widehat S(\Phi_M)-S^\star)^2\big]-\E_{\P_{\mathrm{tr}}}\big[(\widehat S(\Phi_M)-S^\star)^2\big]\Big|\le C\, W_1(\P_0,\P_{\mathrm{tr}}),
    \]
    where $W_1$ is the 1-Wasserstein distance and the constant \(C\) depends only on the Lipschitz and boundedness constants of the predictor class.
\end{remark}

A formal version of this result, together with a complete proof, is provided in Section~S4.3 of the Supplementary Material. The remark provides a stability interpretation for the following empirical misspecification
experiments. In particular, if the
distributional shift between the training and evaluation laws (e.g., discrepancies in volatility dynamics, jump behavior, or microstructure noise) is moderate in
Wasserstein distance, then the corresponding change in squared prediction risk is
controlled proportionally. The remark itself is a standard bounded-Lipschitz stability bound whose role
here is to connect the empirical misspecification experiments with the
distributional robustness properties implied by the Lipschitz structure
of the predictor class.

In the simulation study summarized below, the training law is generated under the
Heston model, whereas some evaluation experiments are conducted under
the Bates model with additional jump components or under the Heston model with violated Feller condition. The stability bound
therefore suggests that the resulting change in prediction risk should
remain controlled whenever the induced shift in the joint law of
\((\Phi_M,S^\star)\) is moderate in Wasserstein distance. A formal
quantitative analysis of Wasserstein distances between such stochastic-volatility
and jump-diffusion specifications is left for future research.

\paragraph{Summary of the misspecification results.}
The detailed results of both experiments are reported in Section~S3.1 of the Supplementary Material (Figures~S2--S4); we summarize them here. First, when jumps are added (Bates model, 3rd column of Table~\ref{Table: sampling distribution}), the pattern of Figure~\ref{fig:allmethods} is preserved. Although jumps enlarge the admissible parameter space and generate more extreme skewness values, the baseline HFN, which is trained on Heston paths only, remains the most accurate and most stable estimator. Training the HFN directly on Bates-type dynamics is therefore not necessary for jump--diffusion data. Second, when the Feller condition is violated for approximately 40\% of the simulated variance paths (4th column of Table~\ref{Table: sampling distribution}), the variance process simulated with the QE scheme of \cite{andersen2008simple} can remain at zero for extended periods, which produces flat price paths and very negative true skewness. None of the estimators, including the HFN, recovers skewness on these outlying paths. We regard this as desirable, since flat price paths carry no information from which skewness could be recovered. Outside these outliers, i.e.\ for true skewness in $[-2,2]$, the HFN again delivers the smallest bias and MSE.

\subsection{Real-world Application}  \label{result: real-world}
In this section, we apply the proposed low-frequency skewness estimator to U.S. stock market data, sourced from the TAQ database, and benchmark its performance against all competing methods. We select 753 liquid stocks over the period from January 1998 to May 2005 and sample prices at five-minute intervals to construct the corresponding HF return series. As in the simulation study, we exclude the overnight jumps for each time series first. Then, we select a one-year rolling window (i.e., $\text{12 months} \times \text{22 days} \times \text{78 intraday returns}$), following \citet{NP2021}, we demean the intraday returns within the rolling window, and shift this window each step by one month (i.e., $\text{22 days} \times \text{78 intraday returns}$) till we reach the end of a series. We do this for all skewness estimation methods. Overall, for each method, we produce a time series of monthly skewness estimates for each stock. 

We demonstrate the utility of our HFN skewness estimator through four distinct applications. First, we provide a detailed analysis of the results for a single stock in Section \ref{sec:real-world single stock}. Second, we perform cross-sectional Fama-MacBeth regressions in Section \ref{sec: real world famamacbeth}. Finally, following the approaches of \cite{Amaya2016} and \cite{NP2021}, we construct an equally weighted long-short portfolio based on cross-sectional skewness sorting and evaluate its performance in Section \ref{sec:real-world portfolio}. 
%
\subsubsection{Single Stock Time Series} \label{sec:real-world single stock}
\begin{figure}[!htbp]
    \centering
    \includegraphics[width=0.9\linewidth]{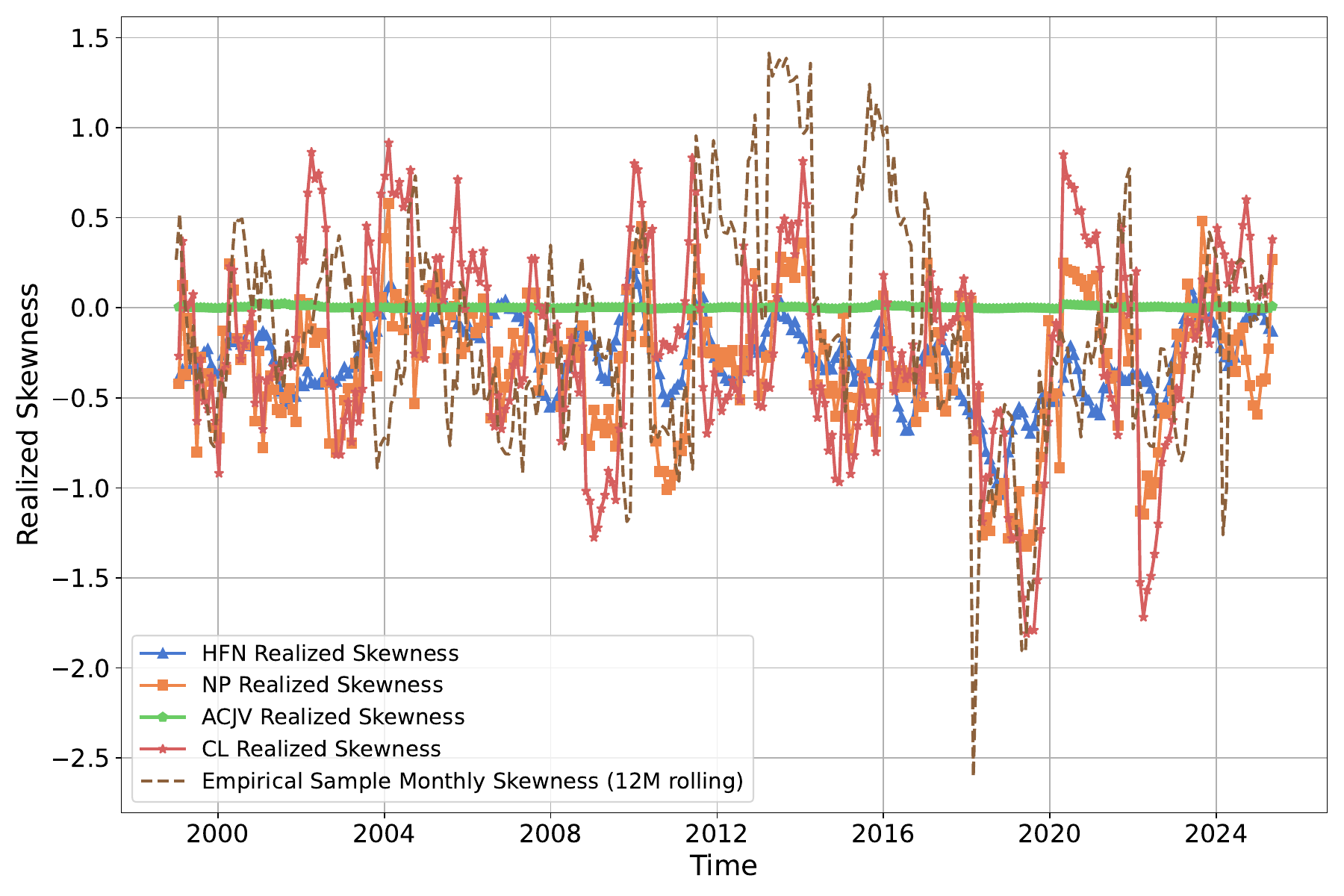}
    \caption{Realized Skewness of MSFT Stock (Rolling 12-Month Window)}
    \label{fig:realworld_single_skew}
\end{figure}

\begin{figure}[!htbp]
    \centering
    \includegraphics[width=0.8\linewidth]{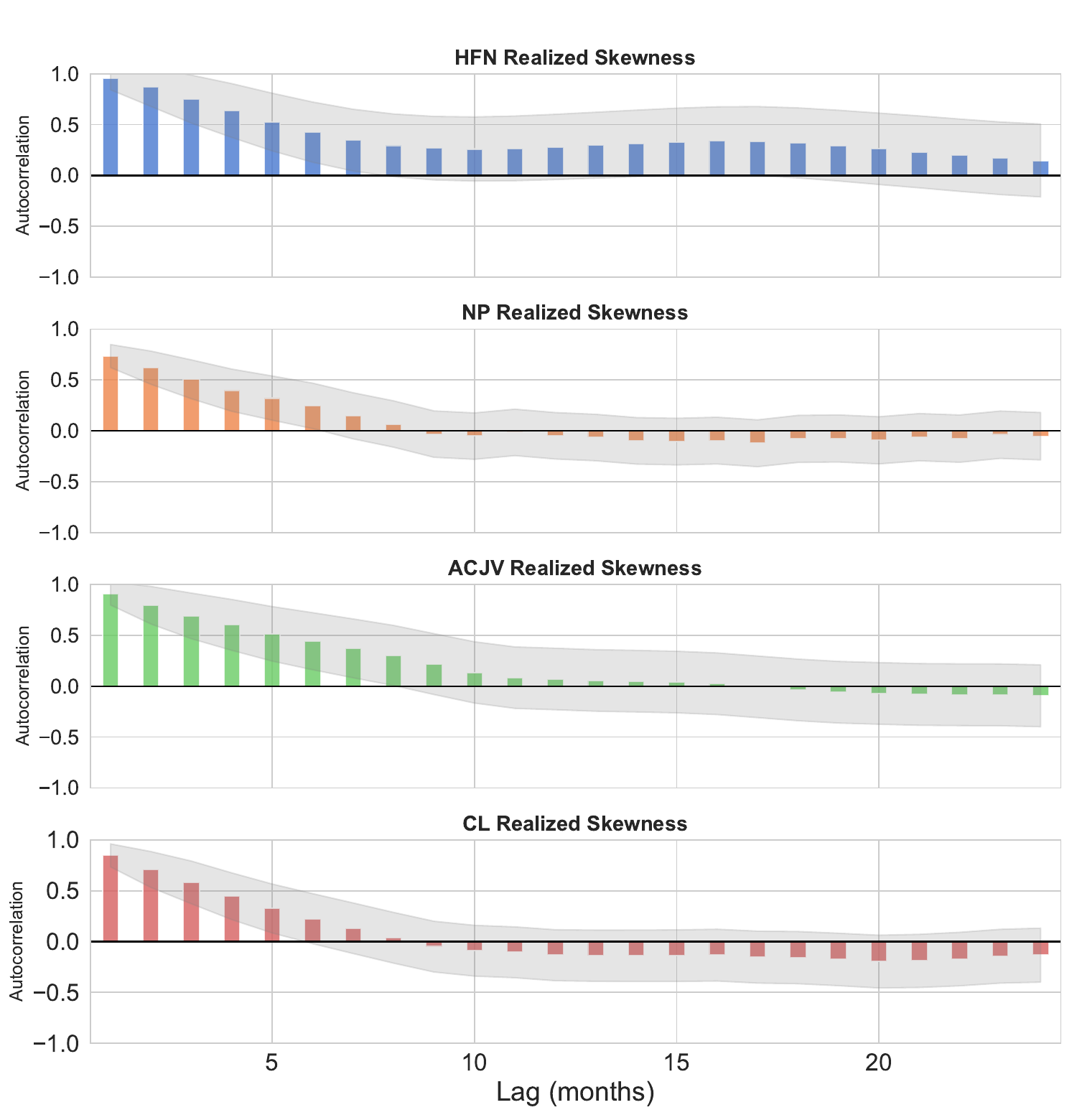}
    \caption{Autocorrelation Function (ACF) of Realized Monthly Skewness (95\% confidence interval marked in grey)}
    \label{fig:realworld_single_acf}
\end{figure}

We illustrate the behavior of the various estimators using HF data for Microsoft (MSFT). For each method, we apply the rolling-window procedure described above to compute realized monthly skewness, and we also calculate the corresponding empirical sample skewness within each window for comparison. The resulting skewness series are shown in Figure~\ref{fig:realworld_single_skew}. Except for the ACJV estimator, all methods display broadly similar dynamics. Among them, our HFN estimator yields the smoothest and most stable series, whereas the alternative estimators exhibit substantial volatility, with the empirical sample skewness showing particularly pronounced fluctuations. Using 327 monthly return observations from 1998 to 2025, the full-sample empirical skewness is approximately $-0.3$, which is consistent with the skewness levels produced by our HFN estimator.
The persistent evolution of HFN skewness, contrasted with the oscillatory patterns of other methods in Figure \ref{fig:realworld_single_skew}, prompts an analysis of their autocorrelations (Figure \ref{fig:realworld_single_acf}). The autocorrelations of the NP, ACJV, and CL estimators decay rapidly and become negative, whereas the HFN's remain positive at all lags, including non-overlapping 12-month horizons. This sustained positive autocorrelation implies that HFN uniquely captures a clustering of skewness. Combined with its proven accuracy, this finding indicates that standard estimators may overlook a critical momentum in tail-risk dynamics, raising concerns of bias in practical applications.
\subsubsection{Cross-sectional Regressions}  \label{sec: real world famamacbeth}
Following \cite{Amaya2016}, we examine the relation between realized skewness and future returns using the \cite{famamacbeth1973} methodology. In each month $t$, we run a cross-sectional regression of individual stock returns in month $t+1$ on realized skewness and a set of control variables, all measured at the end of month $t$:
\begin{equation}\label{eq:famamacbeth}
r_{i,t+1} = \gamma_{0,t} + \gamma_{1,t}\, RSkew_{i,t} + \gamma_{2,t}\, RV_{i,t} + \gamma_{3,t}\, MOM_{i,t} + \gamma_{4,t}\, LagRet_{i,t} + \epsilon_{i,t+1}.
\end{equation}
$RSkew_{i,t}$ is the realized skewness of stock $i$ estimated by HFN, NP, ACJV or CL over the 12-month window ending in month $t$. The control variables are computed from the same HF data. $RV_{i,t}$ is the logarithm of the realized variance over the same 12-month window, i.e.\ the sum of squared demeaned intraday five-minute returns. $MOM_{i,t}$ is the cumulative return from month $t-11$ to month $t-1$, which skips the most recent month. Following \cite{Amaya2016}, $LagRet_{i,t}$ denotes the lagged return, i.e.\ the return of stock $i$ in month $t$, which captures the short-term reversal effect. For each month, all regressors are winsorized at the 1st and 99th percentiles and standardized cross-sectionally to zero mean and unit variance. A slope coefficient therefore measures the change in next-month return associated with a one-standard-deviation increase in the regressor, and the coefficients are directly comparable across the four skewness estimators. The reported coefficients are time-series averages of the monthly estimates $\hat\gamma_{k,t}$. Their $t$-statistics are based on \cite{Newey-West1987} standard errors with 12 lags, which account for the serial correlation induced by overlapping 12-month estimation windows.
 
Our sample includes 165 highly liquid stocks with few zero returns over the five-minute intervals.\footnote{Our dataset differs from those of \cite{NP2021} and \cite{Amaya2016}. From the S\&P 500 constituents, we retain stocks with five-minute return data covering 1998 to 2025, resulting in 165 stocks. Skewness and realized variance are computed from intraday returns, excluding overnight returns and early-closing days. Monthly returns, momentum and the lagged return are computed from close-to-close prices over blocks of 22 trading days.} The first 12-month estimation window ends in April 1999. Unlike studies that consider multiple higher moments as predictors, we focus exclusively on skewness but compare four different estimators. The results are reported in Table~\ref{tab:fmb_skew}.
\begin{table}[!htbp]
\centering
\caption{Fama--MacBeth Cross-Sectional Regressions with Control Variables.}
\small 
\begin{threeparttable}
\begin{tabular}{@{}lcccc@{}}
\toprule
 & HFN & NP & ACJV & CL \\
\midrule
Intercept          & 0.0110*** & 0.0110*** & 0.0110*** & 0.0110*** \\[-1ex]
                   & (4.03)    & (4.03)    & (4.03)    & (4.03)    \\
Realized Skewness  & 0.0011**  & 0.0005    & -0.0001   & -0.0002   \\[-1ex]
                   & (1.99)    & (0.72)    & (-0.21)   & (-0.31)   \\
Realized Variance  & 0.0016    & 0.0018    & 0.0015    & 0.0015    \\[-1ex]
                   & (1.17)    & (1.25)    & (1.10)    & (1.10)    \\
Momentum           & -0.0008   & -0.0009   & -0.0009   & -0.0007   \\[-1ex]
                   & (-0.61)   & (-0.72)   & (-0.71)   & (-0.54)   \\
Lagged Return      & -0.0021***& -0.0021***& -0.0020***& -0.0021***\\[-1ex]
                   & (-2.79)   & (-2.82)   & (-2.68)   & (-2.86)   \\
\bottomrule
\end{tabular}
\begin{tablenotes}[flushleft]\footnotesize
\item \textit{Note:} This table reports cross-sectional regressions of next-month stock returns on realized skewness, estimated by HFN, NP, ACJV and CL, controlling for realized variance (in logs), momentum (cumulative return from month $t-11$ to $t-1$) and the lagged return (month-$t$ return). $t$-statistics based on \cite{Newey-West1987} HAC standard errors with 12 lags are shown in parentheses. $*$: $p<0.10$, $**$: $p<0.05$, $***$: $p<0.01$.
\end{tablenotes}
\end{threeparttable}
\label{tab:fmb_skew}
\end{table}
Table~\ref{tab:fmb_skew} displays very similar patterns for the NP, ACJV and CL estimates. Their skewness coefficients are close to zero and far from statistical significance, so those estimates do not establish the existence of a skewness premium in stocks. This finding aligns with \cite{NP2021}, who argue that investors demand compensation for intraday skewness, whereas longer-term asymmetry does not elicit a similar risk premium. In contrast, the coefficient on the more stable HFN skewness estimate is positive and significant at the 5\% level ($t=1.99$). A one-standard-deviation increase in HFN skewness is associated with a 0.11\% higher return in the following month, or about 1.3\% per year. The positive relation reflects the skewness clustering phenomenon discussed in Section~\ref{sec:real-world single stock}, indicating that stock-level skewness is persistent and may carry upside potential over subsequent periods.
 
The HFN skewness effect is not driven by the control variables. Across all four specifications, the lagged return is the only significant control, with a coefficient of about $-0.21\%$ ($t\approx -2.8$), consistent with the short-term reversal documented by \cite{Amaya2016}. Realized variance and momentum are insignificant once the other regressors are included. Moreover, the average monthly cross-sectional Spearman correlations between HFN skewness and the controls are small: $0.02$ with realized variance, $-0.09$ with momentum and $-0.01$ with the lagged return. Hence, HFN skewness is not a smoothed proxy for past returns or volatility. By contrast, the NP and CL skewness estimates are noticeably correlated with the lagged return ($-0.15$ and $-0.13$), suggesting that these estimators partly reflect the most recent price movements rather than the persistent component of return asymmetry.

\subsubsection{Allocation based on Skewness sorts} \label{sec:real-world portfolio}
The foundation of this investigation is \cite{Amaya2016}, who show that portfolios longing the stocks with high five-minute skew while shorting the stocks with low five-minute skew yield significant negative returns subsequently. Based on this work, \cite{NP2021} have expanded the skew estimation window to a longer horizon (i.e., daily and weekly) and demonstrated that the evidence of a compensation for exposure to negative skew disappears when the estimated skew is targeted at a longer horizon. According to their finding, it implies that investors do not have preference for holding stocks with greater skewness at longer horizons. However, in view of the patterns observed in our simulations and in Section~\ref{sec:real-world single stock}, it may be worth reflecting on whether the long-short portfolio findings in \cite{NP2021} could in part be related to how effectively their estimator captures the persistence of skewness over longer horizons.

Our work uses a longer investigation period (1998–2025) than \cite{NP2021} (2009–2019) and \cite{Amaya2016} (1993–2013) and focuses on monthly skewness estimators for different methods. As in \cite{NP2021}, we construct a long-short portfolio by buying an equally weighted portfolio of the top third of stocks with the highest realized monthly skewness and selling the bottom third. Skewness is estimated each month using a one-year rolling window of HF data, and the portfolio is rebalanced monthly.
\begin{table}[!htbp]
\caption{Skew-sorted Long--Short Portfolio Returns (Monthly): Summary Statistics.}
\label{tab:ls_returns_statistics}
\centering
\small
\begin{threeparttable}
\begin{tabular}{lcccc}
\toprule
& \multicolumn{4}{c}{\textit{Skewness Estimation Method}} \\
\cmidrule(lr){2-5}
\textbf{Period} & \textbf{HFN} & \textbf{NP} & \textbf{ACJV} & \textbf{CL} \\
\specialrule{\heavyrulewidth}{0pt}{0pt}
1998--2025 & 0.24\% & -0.00\% & -0.22\% & -0.14\% \\
           & (1.59) & (-0.02) & (-1.68)  & (-0.84) \\
2020--2025 & 0.61\% & -0.20\% & -0.35\%  & -0.42\% \\
           & (1.36) & (-0.53) & (-0.98)  & (-0.86) \\
\specialrule{\heavyrulewidth}{0pt}{0pt}
\end{tabular}
\begin{tablenotes}[flushleft]\footnotesize
\item \textit{Note:} This table displays the mean and the corresponding $t$-statistics (in parentheses) of the monthly returns for different long/short portfolios constructed by various skewness estimation methods.
\end{tablenotes}
\end{threeparttable}
\end{table}
\begin{figure}[!htbp]
    \centering
    \includegraphics[width=0.9\linewidth]{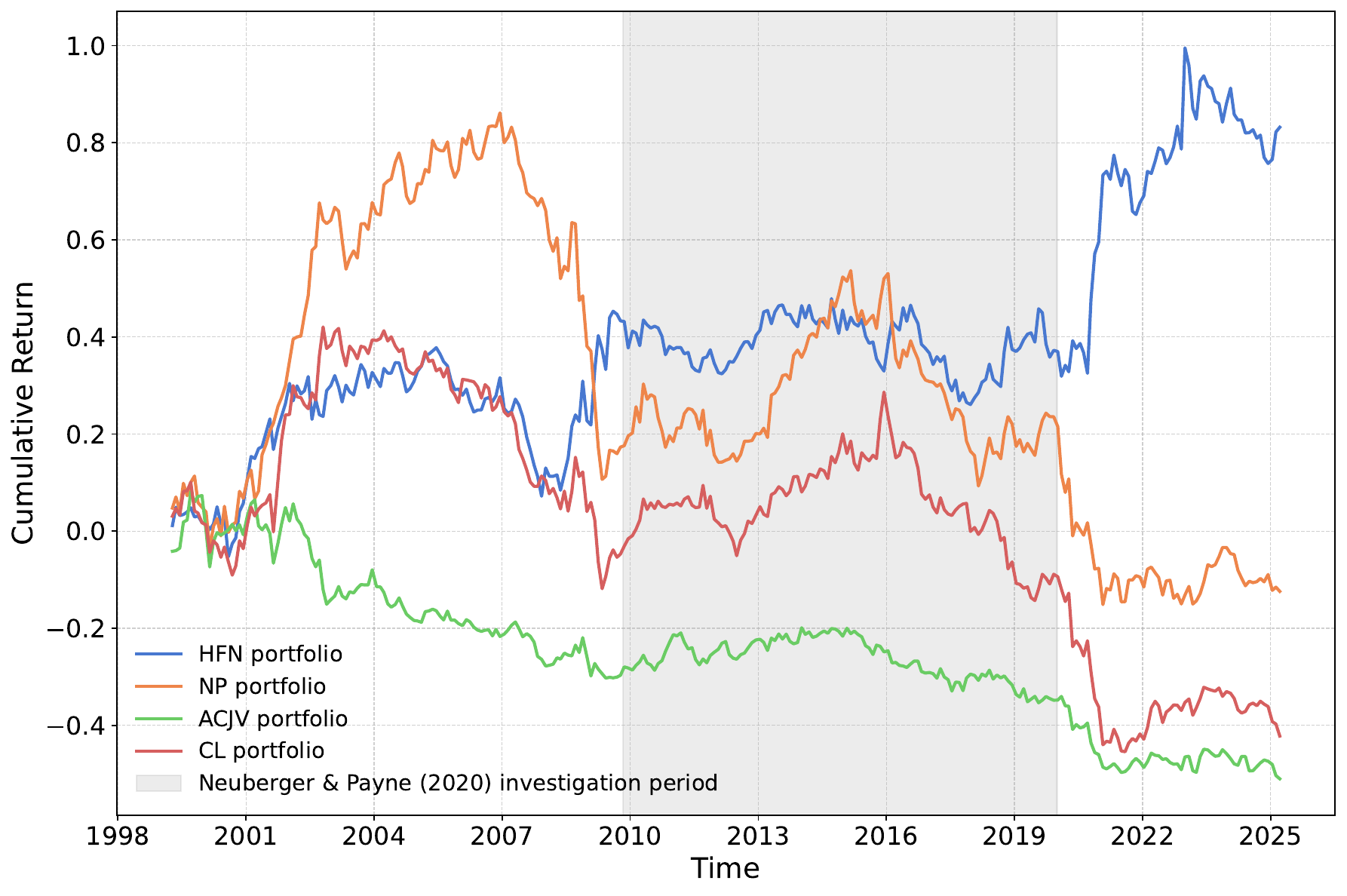}
    \caption{Cumulative long-short returns from sorting stocks on skewness. \textit{Note:} The long-short portfolios result from sorted realized skewness derived respectively with HFN, NP, ACJV, and CL. The realized skewnesses are all based on a monthly time horizon. Portfolio rebalancing is monthly.}
    \label{fig:realworld_portfolio}
\end{figure}
Let us now compare the cumulative returns obtained for the various portfolios depending on the skewness measure. The trajectories are displayed in Figure \ref{fig:realworld_portfolio}. The period over which one can compare our results directly with the ones of \cite{NP2021}, that is, the period 2009-2019, is represented with a gray background. 
In the gray window, the evolution of the cumulative returns obtained by NP in our experiment is similar to what \cite{NP2021} obtained (their Figure 10). The portfolio constructed using the NP method yields an insignificant return. Even if we extend the investigation period to 1998-2025, as illustrated in Table \ref{tab:ls_returns_statistics}, there is no clear evidence of meaningful reward for exposure to negative skewness. 

Unlike earlier studies, our HFN-based long-short portfolio generates positive average returns. This was especially evident during the 2020 post-COVID rebound, where our method showed a higher likelihood of positive performance (Table \ref{tab:ls_returns_statistics}). We reconcile this with existing theory as follows: while investors are compensated for intraday jump risk \citep{Amaya2016, NP2021}, this component is minor at longer horizons. Instead, monthly skewness is primarily driven by the leverage effect. We attribute our positive returns to HFN's capture of skewness clustering, an effect neglected by other estimators.

What can explain this outperformance of the HFN estimator? To investigate this question, we compute out-of-sample return observations with respect to 4 long-short portfolios (HFN, NP, ACJV and CL) and then estimate the out-of-sample skewness of each long-short portfolio using the HFN estimator uniformly. Then, the number of positive and negative skewness estimates is counted. Table \ref{tab:ls_no_skew} shows that the HFN-based portfolio has a 58\% probability of positive subsequent skewness, while all other methods fall below 50\%. Notably, the ACJV method produces portfolios with opposite-sign skewness 84\% of the time. This is intuitive: as ACJV captures only the intraday jump component, it reflects a skewness reversal consistent with compensation for tail risk. This fundamental difference explains the divergence in cumulative returns after 2020 (Figure \ref{fig:realworld_portfolio}). In contrast to the benchmarks, HFN captures skewness persistence, allowing its portfolio to maintain a positively skewed profile and thus generate positive returns during the market rebound. In other words, with the HFN estimate there is a momentum effect the other methods miss.
\begin{table}[!htbp]
\centering
\caption{Out-of-sample Portfolio Skewness.}
\label{tab:ls_no_skew}
\small
\begin{threeparttable}
\begin{tabular}{lcccc}
\toprule
& \multicolumn{4}{c}{\textit{Predicted Skewness}} \\
\cmidrule(lr){2-5}
 & \textbf{HFN} & \textbf{NP} & \textbf{ACJV} & \textbf{CL} \\
\specialrule{\heavyrulewidth}{0pt}{0pt}
No. of positive skew & 166 & 90  & 46  & 84  \\
No. of negative skew & 118 & 194 & 238 & 200 \\
Total Observations   & 284 & 284 & 284 & 284 \\
\specialrule{\heavyrulewidth}{0pt}{0pt}
\end{tabular}
\begin{tablenotes}[flushleft]\footnotesize
\item \textit{Note:} This table displays the counts of positive/negative skewness for different long/short portfolios.
\end{tablenotes}
\end{threeparttable}
\end{table}

\section{Conclusion} \label{sec: conclusion}
This paper proposes a \textit{High Frequency Network} (HFN) framework for learning
low-frequency distributional characteristics from HF financial
data. From a theoretical perspective, we study a HFN family and
show that, under standard universality conditions for convolutional neural
networks, HFN can converge to the optimal predictor
associated with the chosen feature representation under squared loss.
For the empirical application considered in this paper, namely the
estimation of monthly skewness from 5-minute returns, we select a
parsimonious two-layer HFN whose convolutional receptive fields are aligned
with the economically motivated intraday-to-daily and daily-to-monthly
aggregation structure of the problem. This baseline HFN combines
structured temporal aggregation with convolutional feature extraction,
allowing the network to progressively summarize information across multiple
time scales.

Previous research has proposed three benchmark estimators: the ACJV method \citep{Amaya2016}, the CL method \citep{choelee}, and the NP method \citep{NP2021}. However, these estimators exhibit in small samples systematic bias and require long historical samples (up to twenty years of HF observations) to achieve stability, relying on the unrealistic assumption of constant underlying dynamics. To overcome these limitations, we employ the proposed baseline HFN trained on one-year-long HF sequences, enabling the estimation of low-frequency skewnesses over shorter and more realistic sequence lengths.

Although the baseline HFN is trained exclusively on data simulated under a Heston stochastic volatility framework, i.e. without jumps and under a satisfied Feller condition, empirical results demonstrate that it achieves superior accuracy and robustness across a wide range of misspecified environments, including processes with jumps and those violating the Feller constraint. These findings indicate that, despite its model-specific training setup, HFN generalizes effectively to more complex stochastic processes. Particularly for negatively skewed return distributions (e.g., skewness near or below –1), which are typical in financial markets \citep{Neuberger2012}, the HFN substantially outperforms all benchmark methods, producing more stable skewness predictions in finite samples. Furthermore, partial dependence analyses (Supplementary Material, Section~S3.2) confirm that the HFN successfully captures the functional relationships between realized skewness and the underlying Heston parameters ($\rho$, $\kappa$, and $\sigma$), highlighting its potential utility in stochastic model calibration and risk factor analysis. 

Beyond simulation experiments, the HFN estimator is applied to U.S.\ HF equity data from 1998 to 2025. The empirical analysis demonstrates that the HFN produces smoother and more persistent skewness dynamics compared to traditional estimators, capturing a momentum-like clustering behavior in tail risk that benchmark methods fail to identify. In cross-sectional regressions, HFN-estimated skewness exhibits a positive and statistically significant relation with subsequent stock returns, suggesting that skewness persistence may carry predictive information. Moreover, in skew-sorted portfolio analyses, portfolios constructed using HFN-estimated skewness display directionally stronger and more stable returns during recovery periods, such as the post–COVID-19 rebound, reflecting its ability to capture persistent asymmetry. These results underscore the empirical relevance of HFN in both risk assessment and asset pricing applications, bridging the gap between HF data modeling and practical portfolio construction.

Overall, the High Frequency Network provides a robust and scalable framework for establishing a function mapping to realized higher moments from HF data. By integrating stochastic modeling, convolutional architectures, and empirical validation, this study contributes a new methodological perspective to the measurement of return asymmetry and highlights promising avenues for future research in financial econometrics, volatility modeling, and machine learning–based market analysis.

\bibliographystyle{elsarticle-harv} 
{\small
\renewcommand{\baselinestretch}{1} 
\bibliography{literature}
}

\end{document}


\maketitle
\thispagestyle{empty}
\setcounter{remark}{2} 
\noindent
This supplementary material is organized as follows. Section~\ref{sec: appendix DesignHFN} details the HFN architecture and Section~\ref{sec: appendix HFN training} describes its training. Section~\ref{supp:simulation} reports additional simulation results: detailed misspecification tests, a parameter-sensitivity analysis, and ablation studies. Section~\ref{supp:theory} contains the theoretical results: a Brownian example to explain the infill limit of the realized skewness, the proof of the finite-sample ratio bias of plug-in skewness estimators, and the statistical learning theory, approximation results, and robustness properties of the High Frequency Network (HFN) estimator. The numbering of definitions, theorems, lemmas, propositions, corollaries and remarks continues that of the main text: results restated from the main text keep their main-text number, and new results are numbered consecutively after the last result of the same type in the main text. Sections, equations, figures and tables of this supplement carry the prefix ``S''.
\section{HFN Architecture Design}\label{sec: appendix DesignHFN}
We apply convolutions to the \emph{reshaped} feature tensor $\Phi_M\in\R^{H\times W\times 3}$ introduced in Section 2.4.3 of the main text. Concretely, let $R^{(q)}\in\R^{H\times W}$, $q=1,2,3$, denote the three input channels of $\Phi_M$, so that $R^{(1)}_{i,j}=r_{i,j}$, $R^{(2)}_{i,j}=r^2_{i,j}$, and $R^{(3)}_{i,j}=r^3_{i,j}$, where $(i,j)$ index the row/column positions after reshaping (not the original time index). We emphasize that $r_{i,j}$ denotes the entry of the reshaped return matrix (Section 2.4.3 of the main text), not the raw chronological increment $r_i=X_{t_i}-X_{t_{i-1}}$.

The HFN consists of $c$ convolutional blocks (depth) with $F$ filters per block (width). Here $c$ equals the number of frequency-scale transitions required to map the input horizon to the target horizon (in our baseline: intraday $\to$ daily and daily $\to$ monthly, hence $c=2$). The width $F$ is a tunable capacity parameter; in the baseline we use $F=4$ filters per block, which is a small but effective choice in our experiments.

Within each block, there are $F$ numbers of 2D kernels designed to make convolutional operations for each channel. Every kernel has size $(K_h,K_w)$, where $K_w$ equals the number of time points aggregated along the within-row (horizontal) direction to reach the next frequency-scale. The kernel height $K_h$ controls the receptive field across adjacent reshaped rows and depends on the reshape shift $k$ (in trading days) used to build $\Phi_M$:
\begin{equation} \label{eq:kernelsize}
    K_h = 
    \begin{cases}
    3, & \qquad k=22 \quad \text{(monthly shift)},\\
    6, & \qquad k=10 \quad \text{(biweekly shift)},\\
    12, & \qquad k=5 \quad \text{(weekly shift)}.
    \end{cases}
\end{equation}
In our main application (5-minute intraday returns mapped to monthly skewness), we use $c=2$ blocks and $F=4$ filters per block.

\subsection{Block 1 (intraday $\to$ daily)}\label{Appx::Block1}
The first block employs a two-dimensional convolution with three input channels and kernel size $(K_h^{(1)},K_w^{(1)})=(3,78)$ to aggregate intraday patterns over the 5-minute grid. For each output filter $f$, the convolution performs a weighted summation across both the spatial kernel window and the three input channels, yielding $F$ daily-level feature maps. For filter $f\in\{1,\cdots,F\}$ we define
\begin{equation}
    \label{eq:first_conv}
    \begin{aligned}
    DailyFeature_{f;i,j} &= b^{(1)}_{f}
    +\sum_{q=1}^{3}\sum_{u=1}^{K^{(1)}_h}\sum_{v=1}^{K^{(1)}_w}
    W^{(1)}_{f,q;u,v}\; R^{(q)}_{(i-1)s^{(1)}_{h}+u,\,(j-1)s^{(1)}_{w}+v},\\
    &\qquad i=1,\cdots,H_{\rm out}^{(1)},\quad j=1,\cdots,W_{\rm out}^{(1)}.
    \end{aligned}
    \end{equation}
Here $W^{(1)}\in\R^{F\times 3\times K_h^{(1)}\times K_w^{(1)}}$ are convolution weights, $b^{(1)}\in\R^{F}$ are biases, and $(s^{(1)}_{h},s^{(1)}_{w})$ are the vertical and horizontal strides.



After convolution we apply a point-wise activation $\sigma(\cdot)$ (ReLU in our implementation) and average pooling with pooling window $(P_h^{(1)},P_w^{(1)})$:
\begin{equation}
    \label{eq:avgpool1}
    \begin{aligned}
    AvgDailyFeature_{f,i,j}
    &= \frac{1}{P^{(1)}_{h}P^{(1)}_{w}} \sum_{u'=1}^{P^{(1)}_{h}} \sum_{v'=1}^{P^{(1)}_{w}}
    \sigma\!\Bigl(DailyFeature_{f,\,i+u'-1,\,j+v'-1}\Bigr),
    \\
    &
    i=1,\cdots,H_{\rm out}^{(1)}-P_h^{(1)}+1,\;
    j=1,\cdots,W_{\rm out}^{(1)}-P_w^{(1)}+1.
    \end{aligned}
\end{equation}
In the baseline HFN, the intraday $\to$ daily aggregation is primarily achieved by the convolutional receptive field $(K_h^{(1)},K_w^{(1)})=(3,78)$ and stride $(s_h^{(1)},s_w^{(1)})$. The subsequent average pooling serves as a stabilizing local smoother. In our main application, we use a pooling window $(P_h^{(1)},P_w^{(1)})=(K_h^{(1)},K_w^{(1)})=(3,78)$.

\subsection{Block 2 (daily $\to$ monthly)}\label{Appx::Block2}
The second block uses kernel size $(K_h^{(2)},K_w^{(2)})=(3,22)$ to aggregate daily patterns into monthly-level features. For output filter $g\in\{1,\cdots,F\}$ we define
\begin{equation}
    \label{eq:second_conv}
    \begin{aligned}
    MonthlyFeature_{g;p,q} &= b^{(2)}_{g}
    +\sum_{f=1}^{F}\sum_{u=1}^{K^{(2)}_{h}}\sum_{v=1}^{K^{(2)}_{w}} W^{(2)}_{g,f;u,v}\; AvgDailyFeature_{f,\,(p-1)s^{(2)}_{h}+u,\,(q-1)s^{(2)}_{w}+v}, \\
    &\qquad p=1,\cdots,H_{\rm out}^{(2)},\quad q=1,\cdots,W_{\rm out}^{(2)}.
    \end{aligned}
\end{equation}
Here $W^{(2)}\in\R^{F\times F\times K_h^{(2)}\times K_w^{(2)}}$ and $b^{(2)}\in\R^{F}$, and $(s^{(2)}_{h},s^{(2)}_{w})$ are the strides in the second block.

After the second convolution we apply the same activation $\sigma(\cdot)$ and average pooling with pooling window $(P_h^{(2)},P_w^{(2)})$:
\begin{equation}
    \label{eq:avgpool2}
    \begin{aligned}
    AvgMonthlyFeature_{g,i,j} &= \frac{1}{P^{(2)}_{h}P^{(2)}_{w}} \sum_{u''=1}^{P^{(2)}_{h}} \sum_{v''=1}^{P^{(2)}_{w}}
    \sigma\!\Bigl(MonthlyFeature_{g,\,i+u''-1,\,j+v''-1}\Bigr).
    \end{aligned}
\end{equation}
In the baseline HFN, the daily $\to$ monthly transition is primarily achieved by the convolutional receptive field $(K_h^{(2)},K_w^{(2)})=(3,22)$ and stride $(s_h^{(2)},s_w^{(2)})$. Average pooling is again used as a local smoother. In our main application, we use a pooling window $(P_h^{(2)},P_w^{(2)})=(K_h^{(2)},K_w^{(2)})=(3,22)$.

\subsection{Fully connected layers}
After the second activation and pooling step, we obtain the pooled
feature tensor
\[
AvgMonthlyFeature
\in
\R^{F\times H_{\rm pool}^{(2)}\times W_{\rm pool}^{(2)}},
\] where $H_{\rm pool}^{(2)}$ and $W_{\rm pool}^{(2)}$ are the height/width after pooling (determined by the pooling window, stride, and padding convention). We flatten it into a feature vector
\[
    \mathbf{x}=\mathrm{vec}\bigl(AvgMonthlyFeature\bigr)\in\R^{d},\qquad d = F\,H_{\rm pool}^{(2)}\,W_{\rm pool}^{(2)}.
\]
We then map $\mathbf{x}$ through two hidden fully-connected layers to obtain a two-dimensional output $\hat{\mathbf{y}}\in\R^2$:
\begin{equation}
    \label{eq:fc_map}
    \hat{\mathbf{y}} = F_{\mathrm{fc}}(\mathbf{x}) = W^{(3)}\,\sigma\!\bigl(W^{(2)}\,\sigma\!\bigl(W^{(1)}\mathbf{x} + \mathbf{b}^{(1)}\bigr) + \mathbf{b}^{(2)}\bigr) + \mathbf{b}^{(3)}.
\end{equation}
Here $\sigma(\cdot)$ is the same activation as above (ReLU in the baseline), $W^{(1)}\in\R^{L_1\times d}$, $W^{(2)}\in\R^{L_2\times L_1}$, $W^{(3)}\in\R^{2\times L_2}$, and $\mathbf{b}^{(1)}\in\R^{L_1}$, $\mathbf{b}^{(2)}\in\R^{L_2}$, $\mathbf{b}^{(3)}\in\R^{2}$. The two outputs correspond to $[\widehat{Vol},\widehat{Skew}]^\top$; the volatility head is used as an auxiliary output to stabilize training, while our theoretical results focus on the skewness output.

\section{HFN Training} \label{sec: appendix HFN training}
The hyperparameters chosen to train the baseline HFN are illustrated in Table \ref{Table: hyperparameters base HFN}.
\begin{table}[htbp]
\centering
\caption{Hyperparameters used for training the base HFN.}
\begin{tabular}{l|l}
\textbf{Hyperparameter} & \textbf{Value}  \\
\hline
Batch size ($N$) & 256 \\
No. of repeats ($e$) & 8 \\
Retained fraction ($z\%$) & 50\% \\
Reshape shift & Monthly ($k=22$) \\
Labels & True skew, True std \\
Activation function & ReLU \\
Loss function & MSE \\
Weight decay ($L_2$) & 0.0001 \\
Gradient clipping & 1.0 \\
Optimizer & Adam \\
Momentum parameters & $\beta_1=0.9$, $\beta_2=0.999$  \\
Learning rate & 0.0005 \\
Compute resources & 4 × Intel Xeon Platinum 8470 CPU (2.00 GHz) \\
Training duration & 168 hours \\
\end{tabular}
\label{Table: hyperparameters base HFN}
\end{table}

\subsection{Empirical Risk}
Because our training data are generated on the fly by simulation, the relevant optimization target is the population risk under the simulated training law rather than the standard empirical risk associated with a single finite training sample. The on-the-fly procedure therefore produces stochastic approximations to this population objective. To assess whether this procedure leads to convergence, we implement a training design experiment for a comparison with training on a fixed pre-generated dataset and with alternative retain fractions ($z\% = \{50\% (Base), 80\%, 0\%\}$).

\begin{figure}[H]
  \centering
  \begin{subfigure}[b]{0.48\textwidth}
    \includegraphics[width=\textwidth]{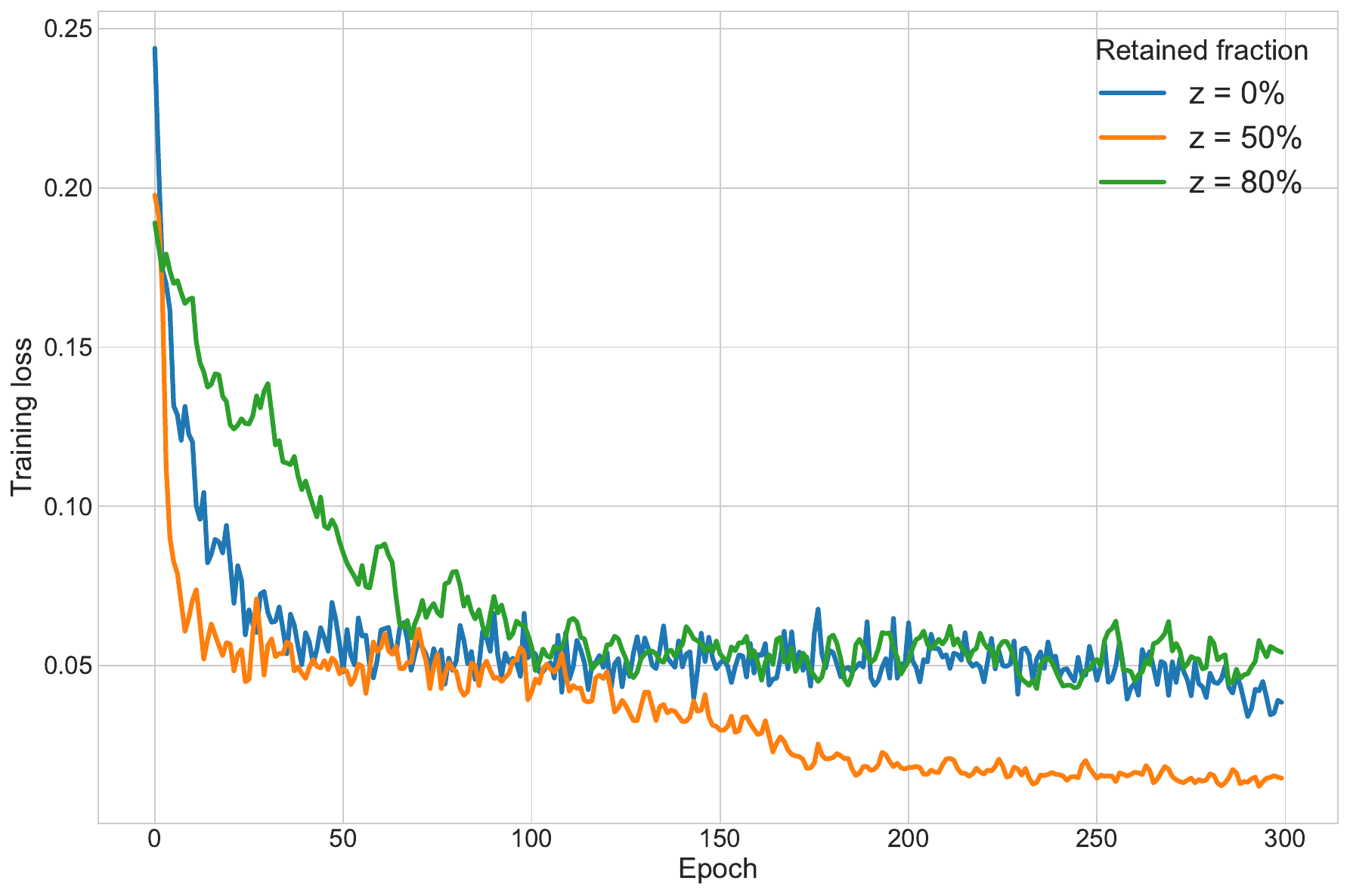}
    \caption{Train loss w.r.t. different retain fractions}
    \label{fig:retain}
  \end{subfigure}\hfill
  \begin{subfigure}[b]{0.48\textwidth}
    \includegraphics[width=\textwidth]{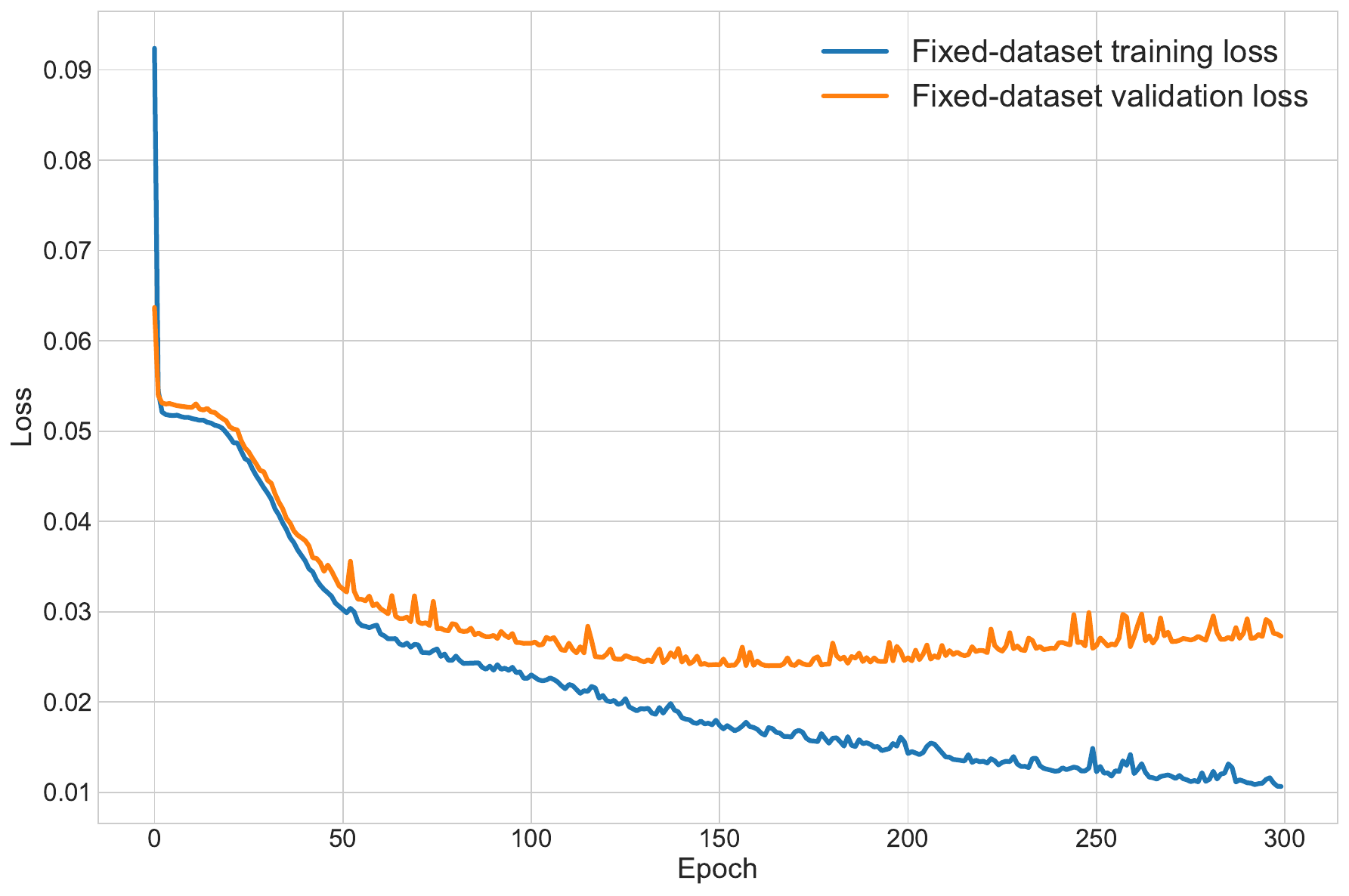}
    \caption{Loss progress of the fixed-dataset training}
    \label{fig:fixeddataset}
  \end{subfigure}
  \caption{Empirical loss observation}
  \label{fig:traindeprogress}
\end{figure}

Figure \ref{fig:traindeprogress} supports the use of the hyperparameters shown in  Table \ref{Table: hyperparameters base HFN}. On the right, the result suggests the presence of overfitting under a static training sample. While the training loss continues to decrease throughout optimization, the validation loss stabilizes at a noticeably higher level and exhibits weaker improvement over later epochs. Using the on-the-fly freshly generated training data avoids this overfitting problem because the network never sees the exactly same sample twice. However, as the plot illustrates on the left, the fully refreshed design ($z=0\%$) exhibits relatively noisy and less efficient convergence, while retaining too large a fraction of the previous batch ($z=80\%$) slows adaptation to newly simulated samples. The intermediate choice $z=50\%$ achieves the lowest training loss, suggesting that it provides a better balance between sample diversity and gradient stability.

\section{Additional simulation results} \label{supp:simulation}
This section reports the simulation experiments summarized in Section 3 of the main text. Unless stated otherwise, the design follows the main text: the test set comprises $1{,}000\times 50$ simulated paths of one year of five-minute log-returns, model parameters are drawn from the distributions in Table 1 of the main text, and the baseline HFN is the network trained under the Heston model with the hyperparameters in Table~\ref{Table: hyperparameters base HFN} (Section~\ref{sec: appendix HFN training}).

\subsection{Misspecification tests: detailed results} \label{supp:misspec}
The theoretical motivation for these experiments, namely the Wasserstein stability bound of Remark 2, is given in Section 3.2 of the main text; a formal version is stated in Section~\ref{sec:distDrift}.

\subsubsection{Presence of Jumps}\label{supp:jumps}
\begin{figure}[!htbp]
  \centering
  \begin{subfigure}[b]{0.48\textwidth}
    \includegraphics[width=\textwidth]{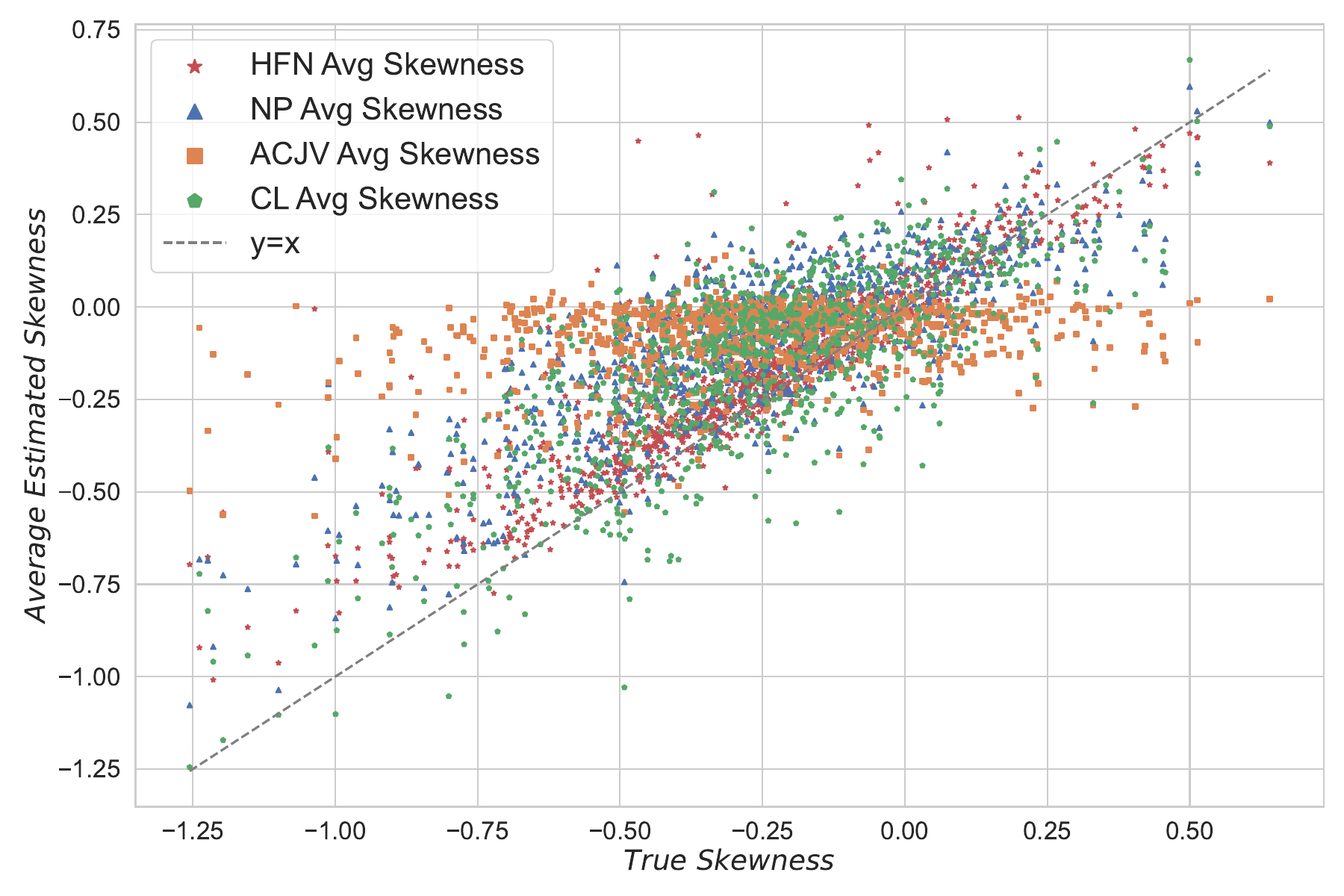}
    \caption{Average estimated skewness}
  \end{subfigure}\hfill
  \begin{subfigure}[b]{0.48\textwidth}
    \includegraphics[width=\textwidth]{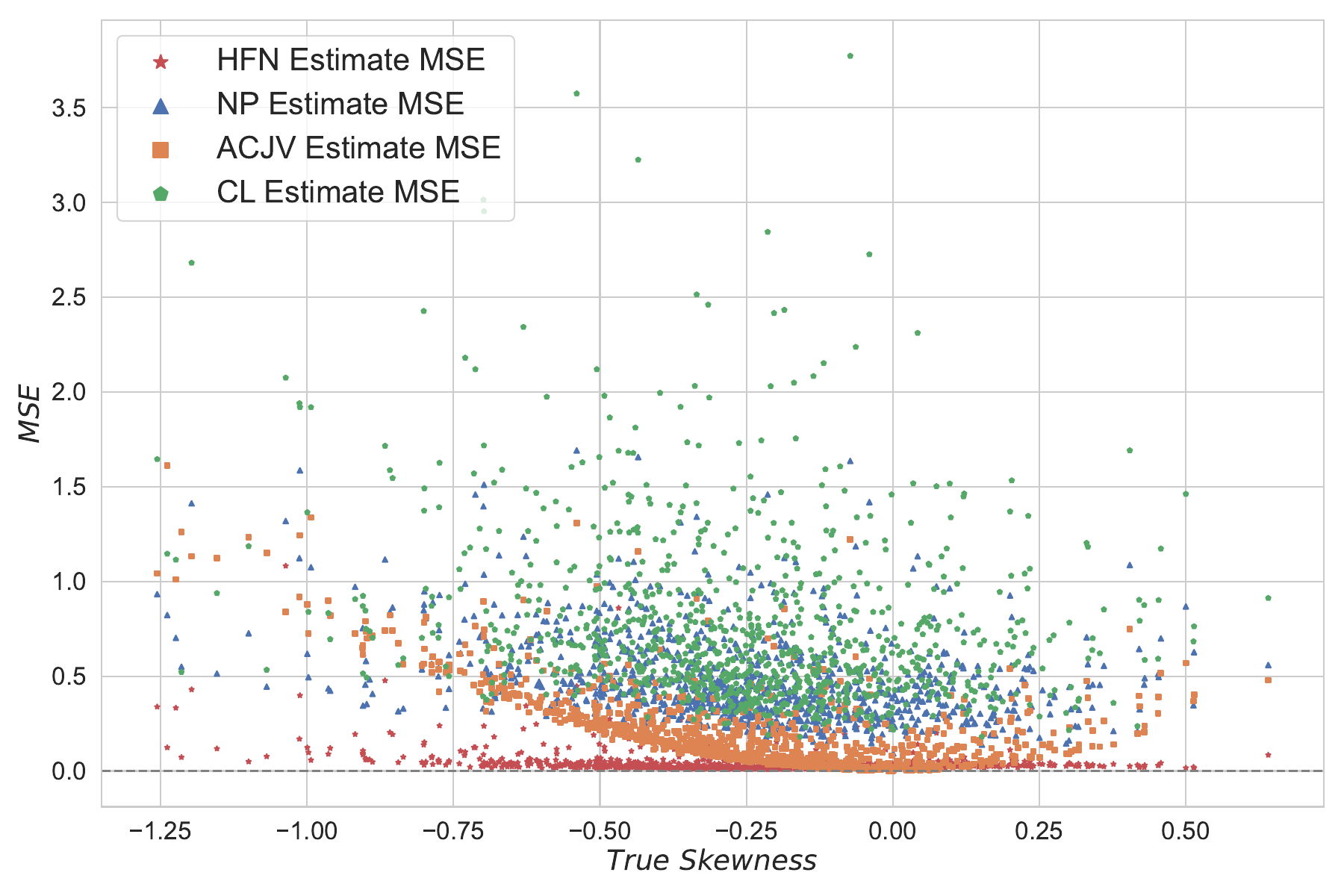}
    \caption{Mean squared error}
  \end{subfigure}
  \caption{Performance comparison under jump–diffusion misspecification.}
  \label{fig:misspecification1}
\end{figure}

By introducing jumps we move to the Bates' process. The performance of each method is depicted in Figure \ref{fig:misspecification1}. In general, compared to Figure 6 of the main text, the observed pattern remains the same in the presence of jumps. Although incorporating jumps enlarges the admissible parameter space and can generate more extreme skewness values, the baseline HFN still delivers the highest accuracy. This result highlights the strong generalization ability of the HFN: a network trained under Heston dynamics retains empirical relevance even when evaluated in a Bates environment. As discussed in Section 2.4 of the main text, extending the estimator to a High Frequency Network trained directly on Bates-type SDEs would be conceptually straightforward. However, the performance of the misspecified HFN test suggests that such an extension is not necessary. The robustness of the estimator indicates that training under the simpler Heston dynamics is sufficient for practical applications involving jump–diffusion data.
\subsubsection{Non-respect of the Feller Condition}\label{supp:feller}
\begin{figure}[!htbp]
    \centering
    \includegraphics[width=\linewidth]{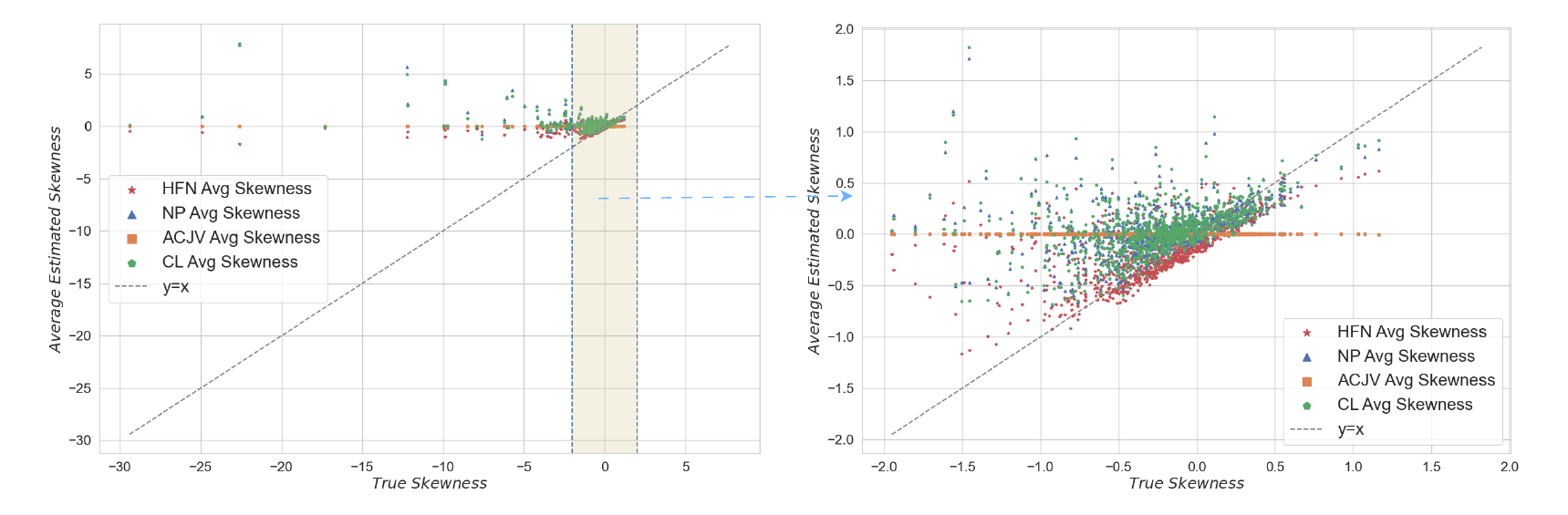}
    \caption{Average estimated skewness (unrestricted Feller conditions). In the figure on the RHS we represent a zoomed version of the colored region in the LHS figure.}
    \label{fig:misspecification2skew}
\end{figure}

\begin{figure}[!htbp]
    \centering
    \includegraphics[width=\linewidth]{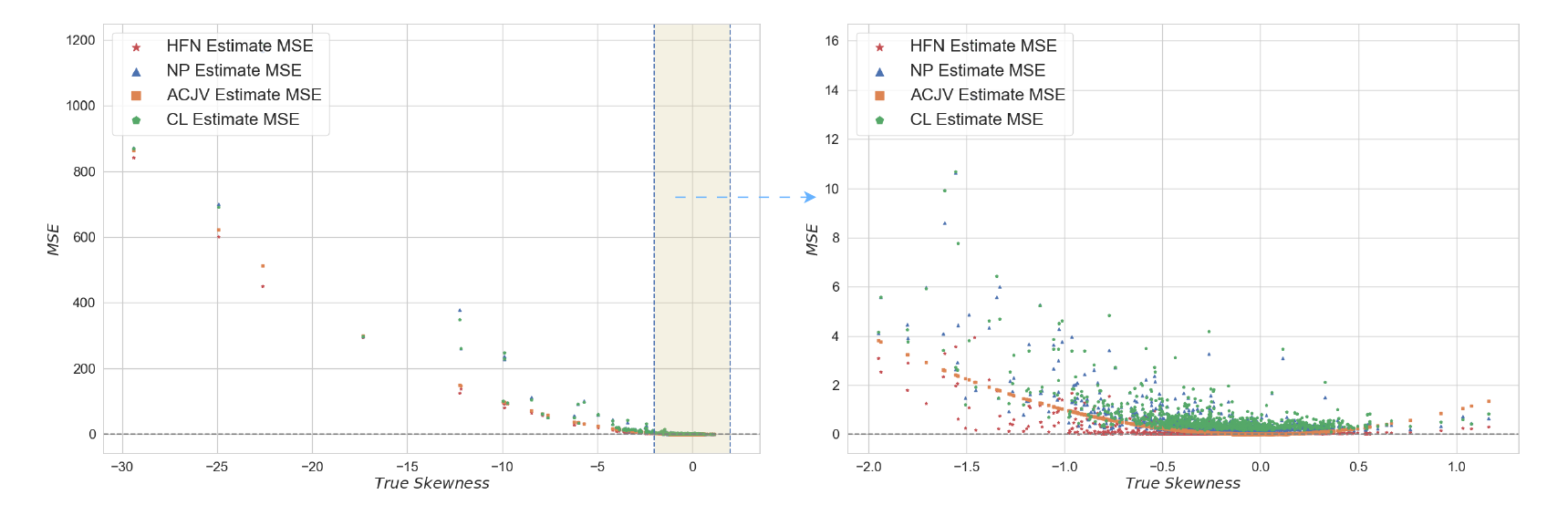}
    \caption{Mean squared error (unrestricted Feller conditions). In the figure on the RHS we zoom into the colored region of the LHS.}
    \label{fig:misspecification2mse}
\end{figure}

In this experiment, we relax the Feller condition by allowing approximately 40\% of the simulated stochastic variance processes to violate it. We simulate, as indicated in Table 1 of the main text (4th column), the Feller ratio $\gamma$ following a shifted exponential. In the worst case, when the draw of the exponential would be 0, the Feller condition would take a value of $10^{-6}$.
Technically, this enlarges the admissible range of the true skewness. As shown in the left panels of Figure \ref{fig:misspecification2skew} and Figure \ref{fig:misspecification2mse}, some realizations exhibit extremely negative true skewness. These observations correspond to outliers for which none of the estimation methods, including our HFN approach, adequately captures their nontrivial behavior.

For the simulations, we employ the QE scheme of \cite{andersen2008simple}. As already noted by Andersen, when the Feller ratio ($\gamma \coloneqq \frac{2\kappa\theta}{\sigma^2}$) is less than 1, simulating trajectories accurately becomes challenging. Inspection of individual sample paths shows that, when the Feller condition is violated, the variance process can hit zero with positive probability. For low Feller ratios, the approximation in the QE scheme produces a degenerate distribution at zero.\footnote{Because of numerical discretization, the variance process may hit zero in simulations even though this cannot occur in continuous time.} As a consequence, zero volatility can persist for extended periods. This produces flat price paths, which in turn generate long stretches of zero returns. Under such conditions, the baseline HFN is unable to outperform the competing methods.

This behavior is in fact desirable. When the input consists of flat price paths, the data contain no information from which meaningful skewness can be recovered. A reliable estimator should reflect this lack of signal rather than fabricate structure. Mapping near-zero returns to substantial skewness would imply a non-invertible transformation and raise concerns about overfitting.

Except for the outliers, as displayed in the right part of Figure \ref{fig:misspecification2skew} and Figure \ref{fig:misspecification2mse} (a zoomed-in view of the true skewness $\in [-2,\ 2]$), the HFN estimator remains best. It continues to yield the strongest performance with respect to both unbiasedness and stability.

\subsection{Parameter sensitivity and estimation methods} \label{result: explainability}
\begin{figure}[!htbp]
    \centering
    \includegraphics[width=0.9\linewidth]{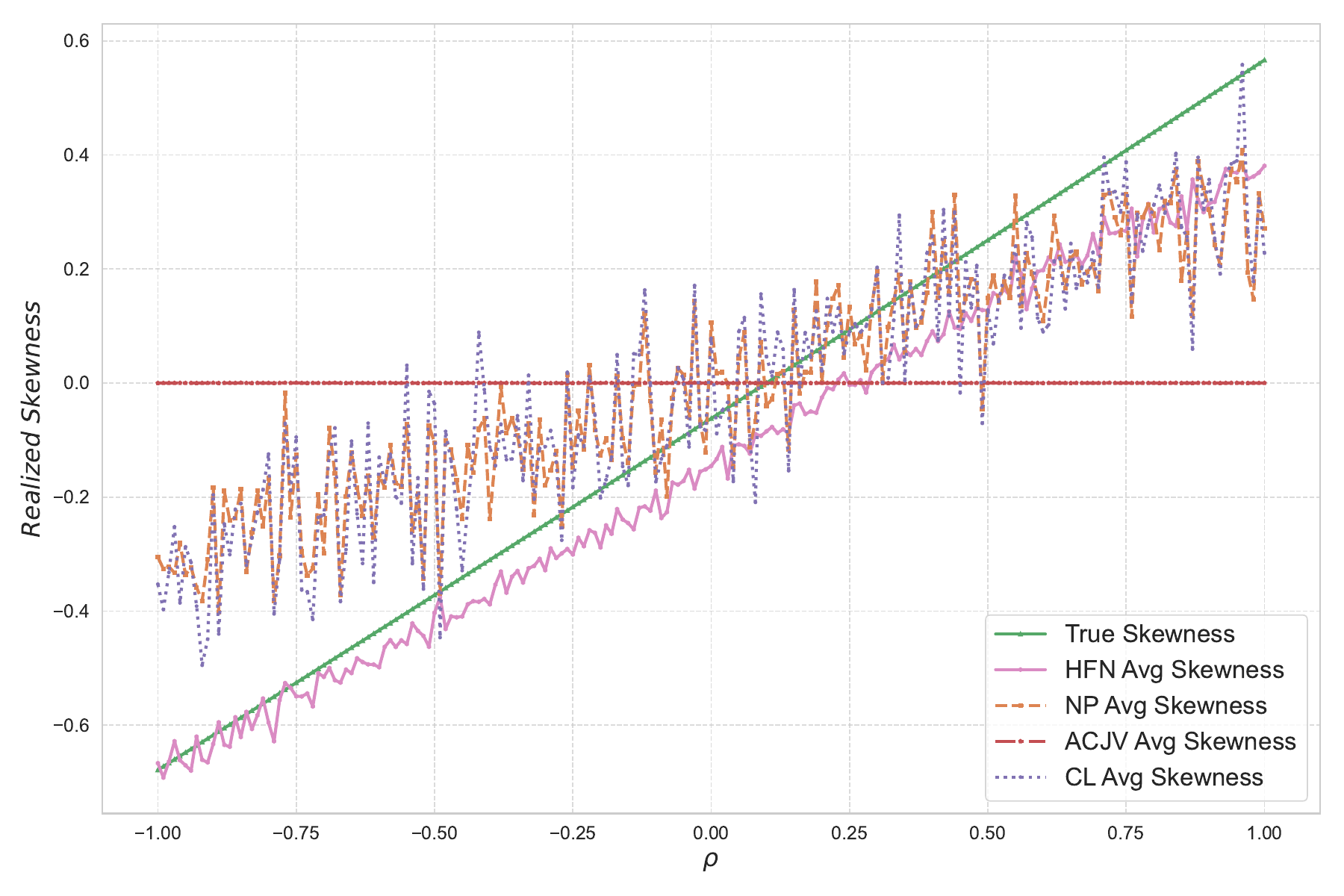}
    \caption{Estimated skewness as a function of $\rho$.}
    \label{fig:explainability1}
\end{figure}

Since the correlation $\rho$ between the two Brownian innovations controls the leverage effect, $\rho$ is considered the primary driver of return asymmetry. Hence, we first investigate the partial dependence of the skewness with respect to $\rho$. As illustrated in Figure \ref{fig:explainability1}, there is a clearly linear relationship between the true skewness and $\rho$. The ACJV's estimates remain anchored at the zero level, demonstrating that it is incapable of capturing diffusive skewness. Although both the NP and CL estimators reproduce a qualitative trend with respect to the parameter $\rho$, the slope of this trend deviates from that of the true skewness, thereby revealing their estimation bias. In addition, the pronounced oscillations in their average estimates highlight their instability. In contrast, our base HFN closely matches the true slope. It produces robust results within the parameter range $[-1, -0.75]$, but the divergence emerges from the point $\rho=-0.7$ upward. Compared to the NP and CL methods, which tend to overestimate the skewness when the true skewness is negative, the HFN slightly underestimates the true level. This conservative bias may be advantageous for risk-averse investors in the event of left-tail occurrences. Regarding positively skewed variables, the HFN method is weakened by increasing bias and decreasing stability. Nevertheless, this shortcoming is of limited practical relevance, due to the typically negative skewness observed in the financial return data \citep{Neuberger2012}.

\begin{figure}[!htbp]
    \centering
    \includegraphics[width=\linewidth]{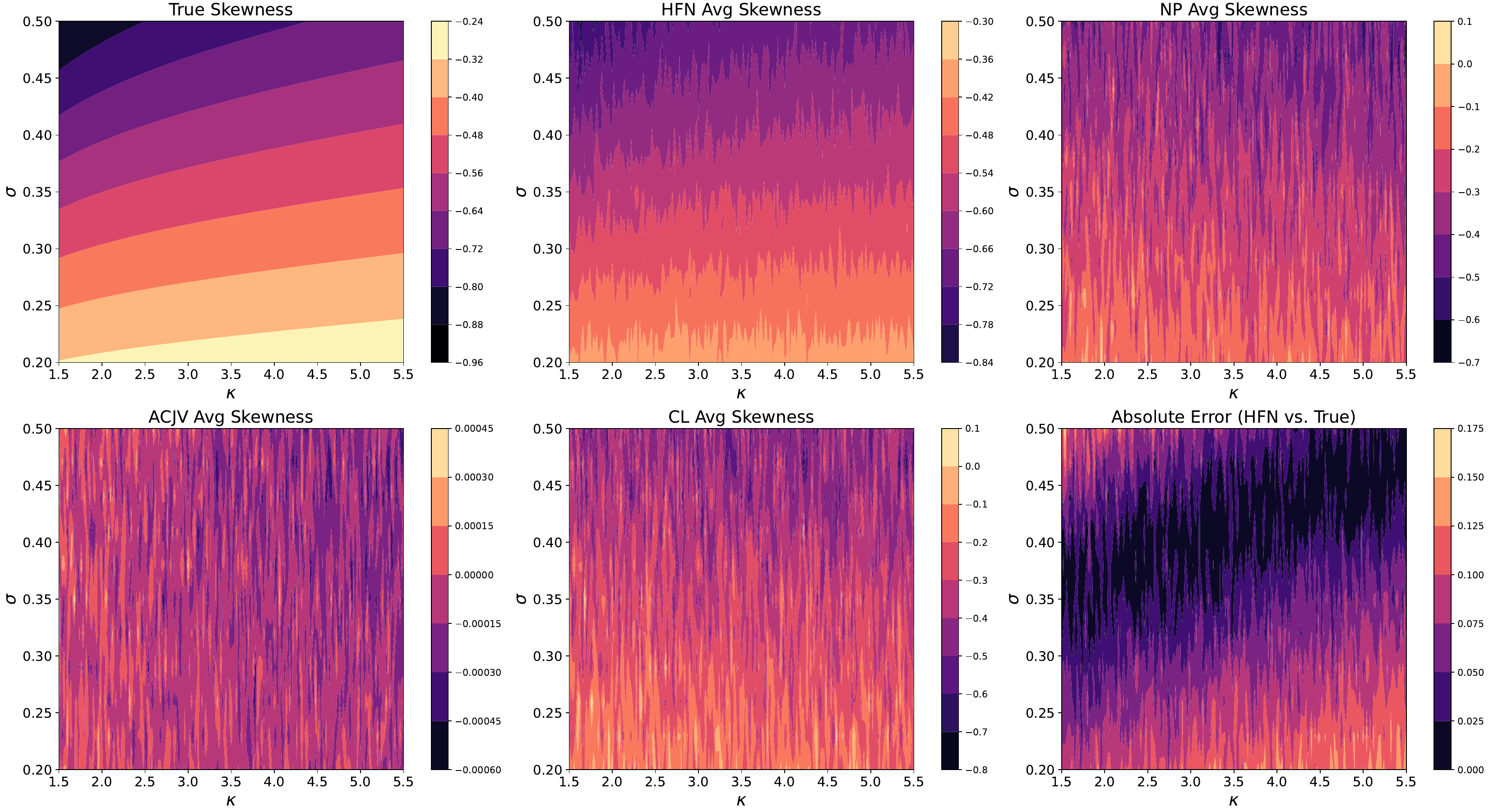}
    \caption{Estimated skewness as a function of $\kappa$ and $\sigma$.}
    \label{fig:explainability2}
\end{figure}

We next explore partial dependence with respect to $\kappa$ and $\sigma$, fixing $\rho=-0.7$ and $\theta=0.04$. The whole comparison is visualized in Figure \ref{fig:explainability2}. The upper left graph displays the true functional sensitivity of the skewness to $\kappa$ and $\sigma$. Clearly, a higher $\sigma$ amplifies the fluctuation of the volatility, and thus magnifies the skewness level. On the contrary, a direct relationship between the skewness and $\kappa$ is not obviously observed. If the volatility of the variance remains at a low level, the return skewness is insensitive to the change of $\kappa$. This behavior indicates a pattern that fast or slow mean reversion speed has minor effect on the skewness when the volatility randomness is low because the low stochasticity in the variance process will hardly generate heavy tails regardless of the direction. However, if the $\sigma$ is high, the mean reversion speed $\kappa$ becomes influential as a fast mean reversion can lead to short-lived volatility spikes and thus reduces the long-run skewness. The upper middle graph in Figure \ref{fig:explainability2} provides evidence of the explanatory power of our base HFN. It successfully captures the sensitivity pattern of true skewness with respect to $\sigma$ and $\kappa$: while $\sigma$ differentiates between distinct skewness levels, $\kappa$ exerts influence primarily when $\sigma$ is large. This interaction is not reflected by any of the other benchmark methods. Recovering such structural dependence, the HFN operator indicates its potential to provide correct values of skewness even in situations where other methods would have difficulties doing so.

\subsection{Ablation Studies} \label{result: ablation}
We conduct ablation analyses to disentangle the drivers of HFN performance across three dimensions: (i) hyperparameters and random initialization, (ii) data format, and (iii) network architecture.
\subsubsection{Hyperparameters and Initialization}
\begin{figure}[!htbp]
  \centering
  \begin{subfigure}[b]{0.48\textwidth}
    \includegraphics[width=\textwidth]{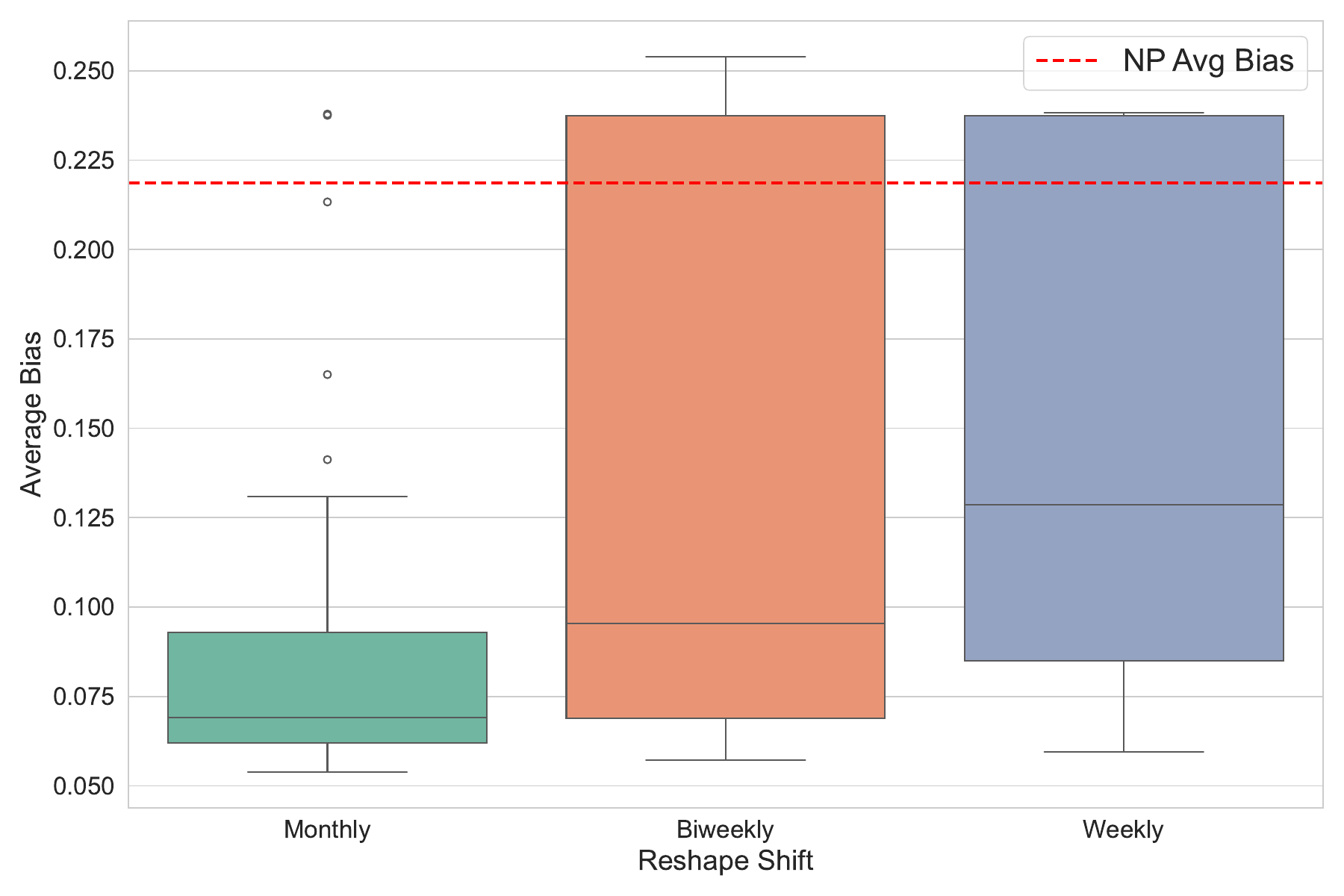}
    \caption{Average bias}
    \label{fig:ablation1bias}
  \end{subfigure}
  \hfill
  \begin{subfigure}[b]{0.48\textwidth}
    \includegraphics[width=\textwidth]{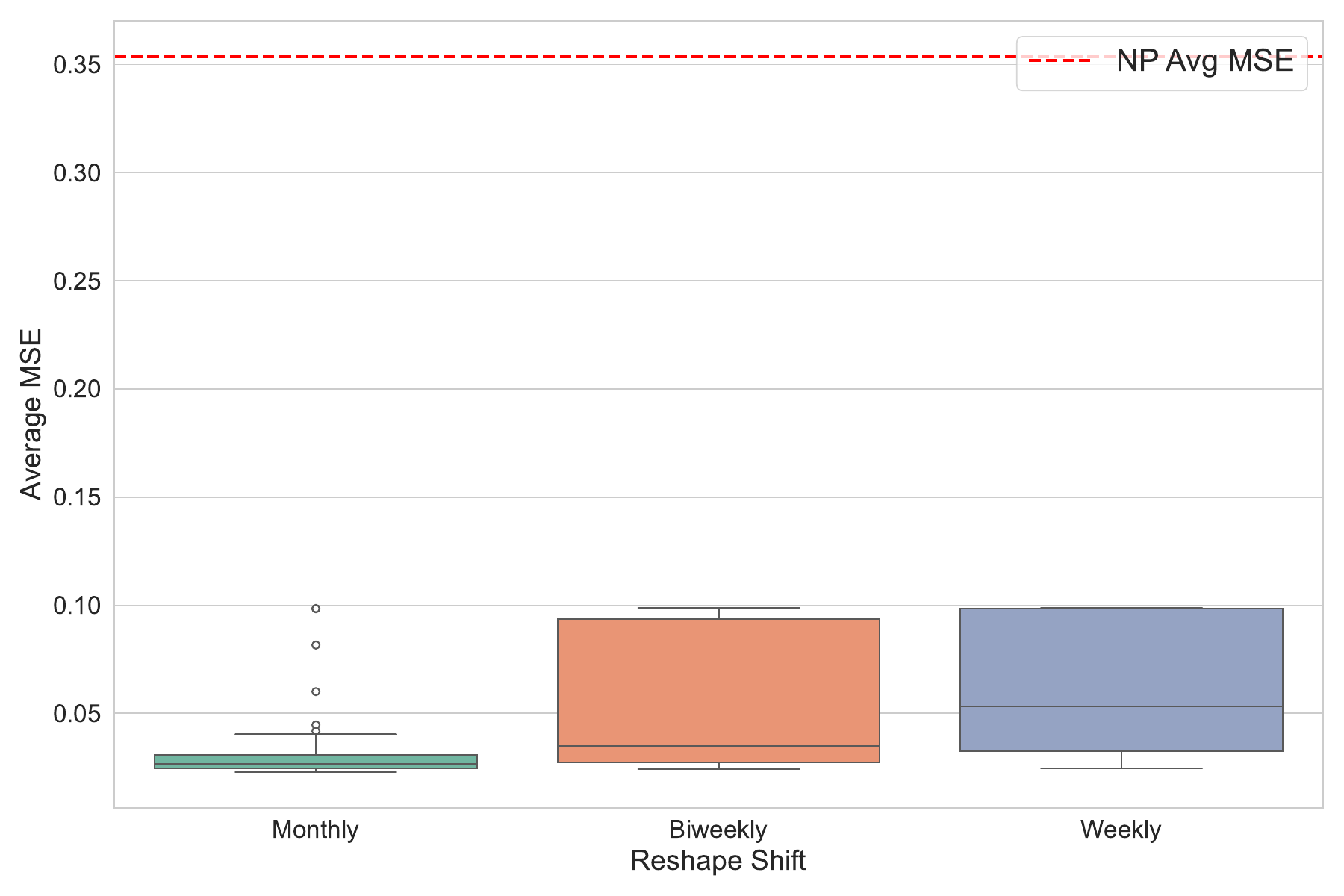}
    \caption{Average MSE}
    \label{fig:ablation1mse}
  \end{subfigure}
  \caption{HFN performance under alternative reshape shifts and hyperparameter settings benchmarked against the NP estimate.}
  \label{fig:ablation1}
\end{figure}

As described in Section 2.4.3 of the main text, the HFN operator is designed to hit input matrices after reshaping the high-frequency time series by shifting a fixed two-month window forward by $k$ days to the next row iteratively until covering the whole series. The reshape-shift index $k$, which controls how the two-month window advances through time, is a key hyperparameter to capture the autocovariance of high-frequency log-returns. Hence, to assess the impact of this factor together with other hyperparameters, we split the first ablation study into three parent cases with reshape shift index $k=22$ (monthly shift), $k=10$ (biweekly shift), and $k=5$ (weekly shift). Following each parent case, we also experiment with varying batch sizes $N = \{256 (Base), 512\}$, number of repeats $e = \{8 (Base), 4\}$, retain percentages $z\% = \{50\% (Base), 80\%\}$, and activation functions $\{ReLU (Base), Tanh\}$. With respect to each configuration, five independent trainings are run with different seeds under identical conditions. All the results are evaluated in the Heston environment. Figure~\ref{fig:ablation1} summarizes bias and MSE outcomes.
The monthly-shift HFN ($k=22$) achieves the lowest bias and MSE. Figure \ref{fig:ablation1bias} demonstrates that using a monthly shift, the HFN method can almost surely beat the benchmark NP method by reducing its estimation bias. In this parent case, other hyperparameters, e.g., batch sizes \& number of repeats, exert relatively small effects on the test outcomes as shown by the green box-plot. 

Applying biweekly or weekly shifts deteriorates the network's performance, and in those two cases, other hyperparameters' influence becomes more significant according to the observation of larger variances of the biweekly and weekly model outcomes. This can be explained by the fact that with a lower reshape shift ($k=10$ or $5$), the corresponding convolutional kernel processes twice (for the biweekly model) and four times (for the weekly model) as many lagged autocovariances within local receptive fields compared to the monthly model, while aiming to estimate the same target—monthly skewness. This redundancy in input length, without a corresponding gain in performance, highlights the relative inefficiency of the biweekly and weekly models. In addition, the use of denser short-term lagged inputs in local pooling amplifies the model's sensitivity to hyperparameter tuning, which in turn contributes to higher output variability in terms of the average bias. Figure \ref{fig:ablation1mse} indicates that regardless of the selected reshaping scheme, the trained model outperforms the benchmark NP method in both accuracy and stability, underscoring the robustness of the CNN-based design. In particular, when the monthly shift is applied, the model yields consistently low MSE across varying hyperparameter settings and random initializations.

%
\subsubsection{Data Format}
\begin{figure}[!htbp]
    \centering
    \includegraphics[width=1\linewidth]{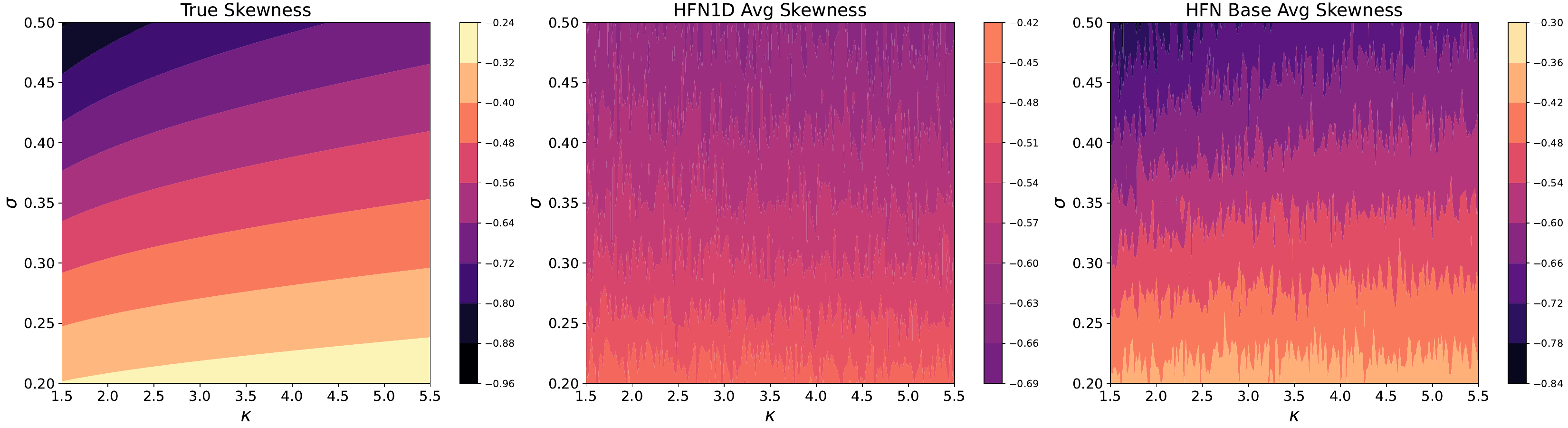}
    \caption{Average estimated skewness of HFN base and HFN1D with respect to $\kappa$ and $\sigma$}
    \label{fig:explainability_1d}
\end{figure}

As described in Section 2.4.3 of the main text, we reshape each input time series into a three-channel, RGB-type matrix using a rolling monthly window. To assess the impact of this data representation, we conduct an experiment in which the input is instead formatted as a one-dimensional time series with three features: returns, squared returns, and cubed returns. We then construct an alternative HFN, called HFN1D, with the same architecture, except that the convolutional kernels are one-dimensional to match the input format. This experiment is conducted within a Heston-model framework. Although the HFN1D method reduces bias to a similar extent and achieves comparably low mean squared errors as our baseline HFN approach, it lacks interpretability. Figure \ref{fig:explainability_1d} displays the findings. 

Comparison between the central and rightmost heatmap reveals that the HFN1D model has difficulties to capture  changes in $\sigma$. It is also unable to capture the effect of $\kappa$ for high $\sigma$. This finding indicates that the use of a one-dimensional time series input induces overfitting, thereby reducing the model’s ability to generalize.

This illustration demonstrates the relevance of having cast the original data into a three-dimensional tensor.

\subsubsection{Network Architecture}
Our HFN estimator is  based on a convolutional neural network. To evaluate the network architecture as one driving factor of the performance, we also used the Long Short-Term Memory (LSTM) specification. The LSTM model was proposed by \cite{hochreiterlstm} and is a specialized form of recurrent neural network (RNN) which is designed to model time series. Moreover, its architecture introduces memory cells equipped with gating mechanisms that regulate the information flow over time. This particular design enables LSTM networks to overcome the vanishing and exploding gradient problems commonly encountered in standard RNNs. Hence, we expect this network architecture to capture global features from high-frequency time series when modeling long-term dependencies and serve as a benchmark against our base HFN framework.

\begin{figure}[!htbp]
  \centering
  \begin{subfigure}[b]{0.48\textwidth}
    \includegraphics[width=\textwidth]{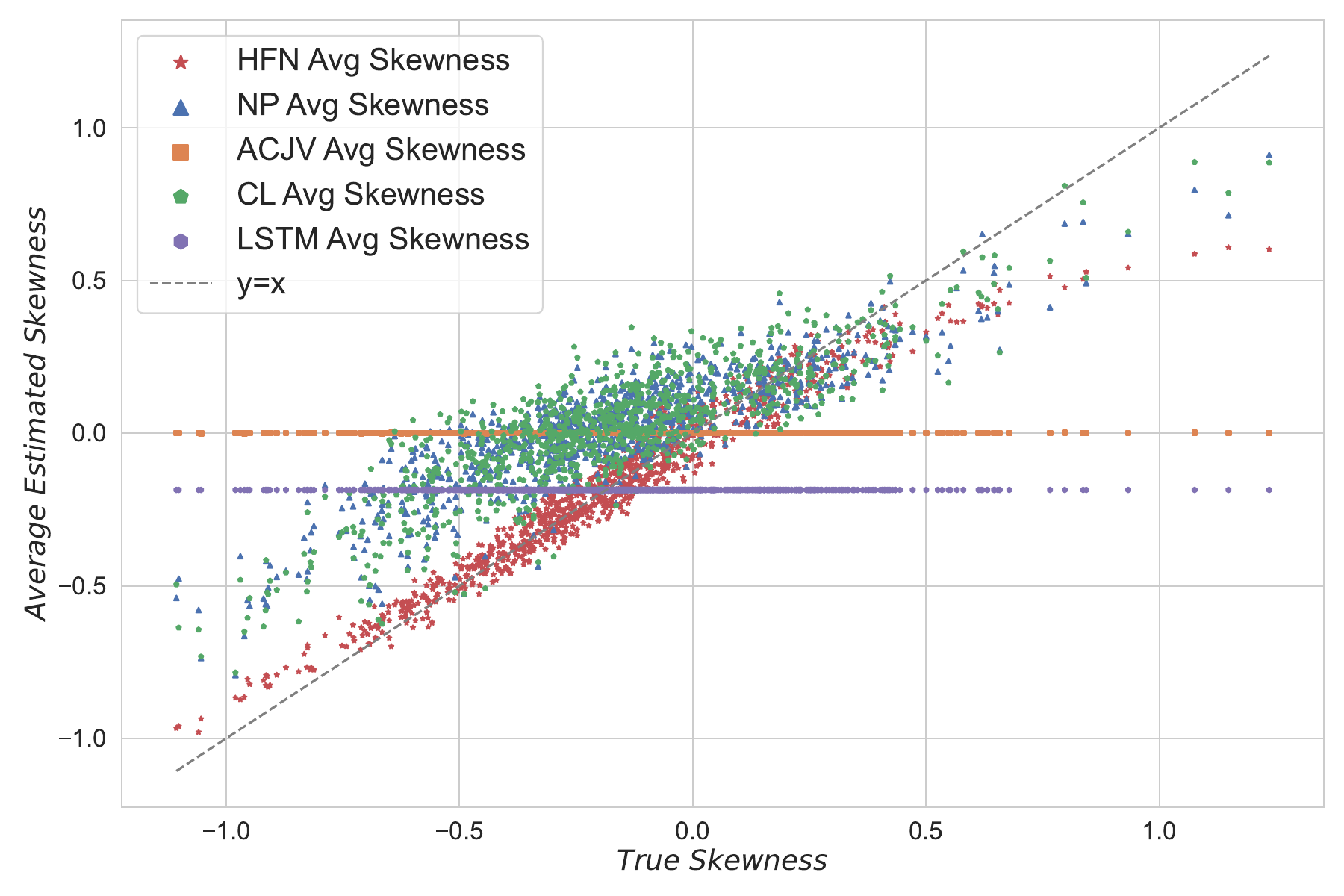}
    \caption{Average estimated skewness}
    \label{fig:lstm-skewness}
  \end{subfigure}
  \hfill
  \begin{subfigure}[b]{0.48\textwidth}
    \includegraphics[width=\textwidth]{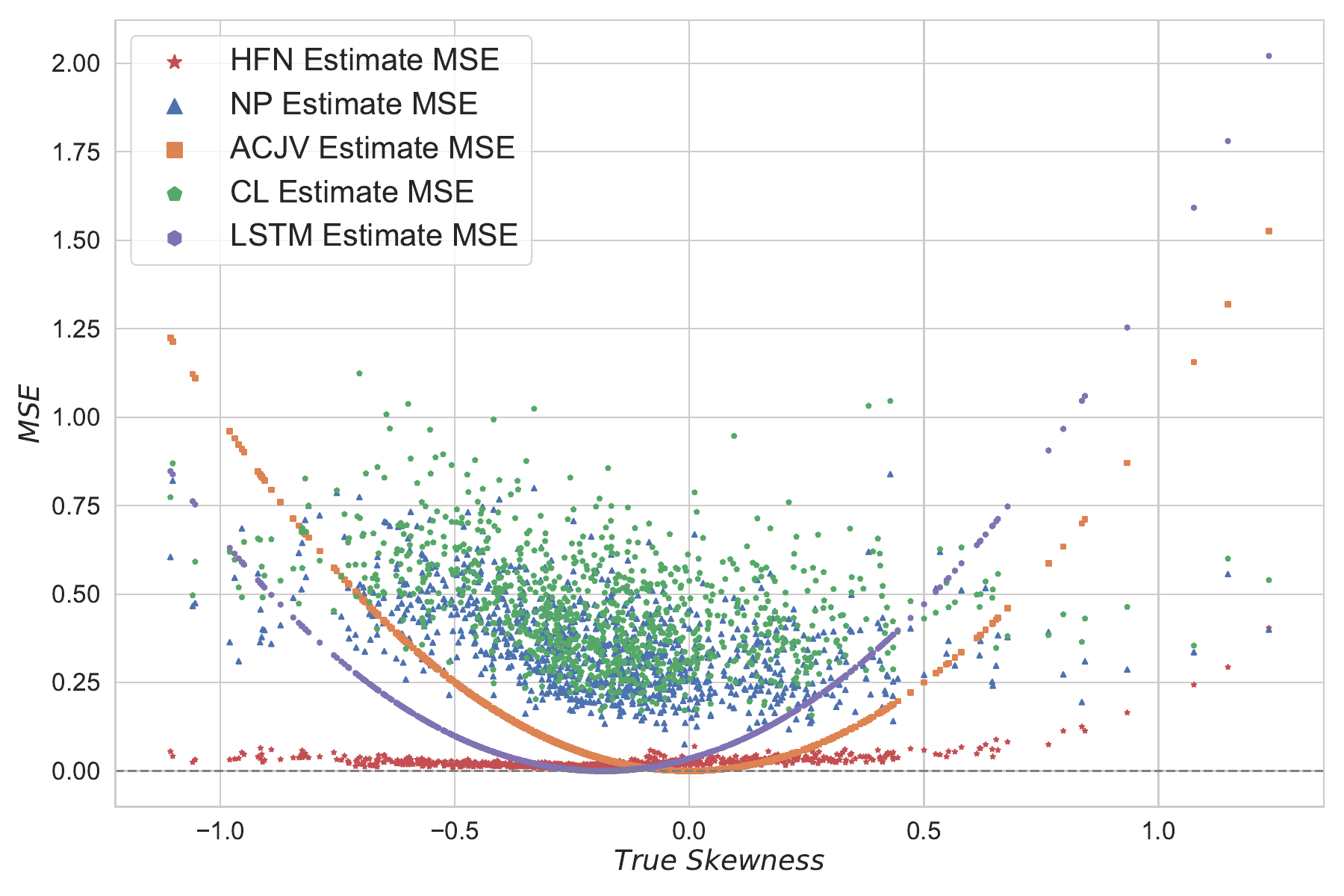}
    \caption{Mean squared error}
    \label{fig:lstm-mse}
  \end{subfigure}
  \caption{Performance of LSTM versus HFN and benchmarks.}
  \label{fig:lstm}
\end{figure}

To ensure the setup is the same as for our HFN, we also endow the LSTM structure with 4 hidden neurons and 2 layers to process 3-channel sequence data. The comparison is depicted in Figure \ref{fig:lstm}. 
The LSTM model yields outcomes similar to those of the ACJV method, for which trivial estimates are produced across the entire parameter space. However, relative to ACJV, the LSTM’s estimates are systematically shifted downward, suggesting that the LSTM architecture tends to exaggerate negative jump sizes while failing to capture the leverage effect—despite being trained under the same conditions as the HFN. We also inspected the training loss of the LSTM. This revealed that the LSTM architecture has difficulties to converge. 

These findings highlight the importance of our custom-designed network architecture as a key factor in achieving superior performance.

%

\section{Theoretical properties} \label{supp:theory}
This section collects the theoretical results. Section~\ref{supp:skewlimit} shows, by an example, why realized third moment and skewness estimators that account for the leverage term have no path-wise infill limit analogous to realized variance that identifies the low-frequency skewness. Section~\ref{supp:ratio-bias} proves the finite-sample ratio bias of plug-in skewness estimators stated in Proposition 1 of the main text. Section~\ref{supp:theory-hfn} contains the statistical learning theory, approximation results, and robustness properties of the HFN estimator.

\subsection{Infill limits of realized third moments: a Brownian example} \label{supp:skewlimit}
Consider a Brownian motion without drift, $X_t=\sigma W_t$, whose skewness at any horizon is
zero, and fix the horizon $T$. By It\^o's formula,
\begin{equation}
  X_T^3 \;=\; 3\sigma^2\int_0^T X_s\,ds \;+\; 3\int_0^T X_s^2\,dX_s
  \qquad\text{a.s.}
  \label{eq:ito_bm}
\end{equation}
The first term is the infill limit of the leverage term in the CL and NP realized third moments,
and the second is a martingale with zero mean. Realized variance converges to the constant
$[X]_T=\sigma^2T$ on every path, so sampling more finely removes all noise from the denominator
of the skewness. The numerator behaves differently. Its limit
$3\sigma^2\int_0^T X_s\,ds\sim\mathcal N(0,3\sigma^6T^3)$ has mean zero, as it must, but remains
a non-degenerate random variable. Infill sampling reduces its variance relative to the single
low-frequency cube $X_T^3$, for which $\Var(X_T^3)=15\,\sigma^6T^3$, but does not remove it.
Consequently, the realized skewness converges to
\[
    \frac{3\sigma^2\int_0^T X_s\,ds}{(\sigma^2T)^{3/2}}
    =\frac{3}{T^{3/2}}\int_0^T W_s\,ds\;\sim\;\mathcal N(0,3)
\]
rather than to the true value zero, whatever the horizon and volatility. By contrast, the sum of
cubed increments converges to zero at rate $1/M$, since it is the cubic variation of a continuous
process; in this model it is therefore correct in the limit, but only because there is no
leverage effect.

Hence, as the sampling interval goes to zero with $T$ fixed, a realized third moment that accounts
for the leverage term converges to a non-degenerate random variable, namely the cube of the
period's return minus a zero-mean martingale term, rather than to the population third moment.
It is unbiased across periods, but within one period it remains random: finer sampling reduces its
dispersion relative to $X_T^3$, but not to zero. This error cannot be averaged away within the
period, only across periods. It explains why estimators of this type, including that of
\cite{NP2021}, need long samples to be accurate. An estimator trained on a model, such as the HFN, targets $\E[S^\star\mid\Phi]$ instead: it uses the structure of the training distribution to map the path to the model-implied skewness, trading the path-wise noise for a shrinkage towards the training distribution.

\subsection{Finite-sample ratio bias of plug-in skewness} \label{supp:ratio-bias}
\setcounter{proposition}{0}
\begin{proposition}[Finite-sample ratio bias of plug-in skewness]
\label{prop:ratio-bias-supp}

Fix a target horizon \(T\), and let
\[
    S^\star=\frac{m_3}{m_2^{3/2}},
    \qquad m_2>0,
\]
where \(m_2\) and \(m_3\) are the second and third central moments of the
low-frequency return. Let
\[
    \widehat S
    =
    \frac{\widehat m_3}{\widehat m_2^{3/2}}
\]
be a plug-in skewness estimator based on estimators \(\widehat m_2>0\) and
\(\widehat m_3\).

Even if \(\widehat m_2\) and \(\widehat m_3\) are unbiased for \(m_2\) and \(m_3\),
respectively, \(\widehat S\) is generally biased in finite samples. The leading
(second-order) bias is
\[
\frac{15}{8}m_3m_2^{-7/2}\Var(\widehat m_2)
-
\frac{3}{2}m_2^{-5/2}\Cov(\widehat m_2,\widehat m_3).
\]
Let \(\widehat m_q=n^{-1}\sum_{j=1}^n Y_{q,j}\), \(q=2,3\), where
\(Y_{q,j}\) is computed on the \(j\)-th of \(n\) non-overlapping low-frequency
estimation windows. Here \(Y_{q,j}\) denotes a realized moment of window \(j\): \(Y_{2,j}\) is the realized variance \(\sum_{i=1}^M r_i^2\), computed from the high-frequency returns of window \(j\), appearing in the denominator of
Eqs. (7), (12) or (15) in the main text, and \(Y_{3,j}\) is the
corresponding numerator. Suppose that
\begin{itemize}
    \item[(a)] \(\{(Y_{2,j},Y_{3,j})\}_j\) is (weakly) stationary with
    \(\E[Y_{q,j}^2]<\infty\), \(q=2,3\);
    \item[(b)] the auto- and cross-covariances are summable,
    \(\sum_{k\in\mathbb Z}\big|\Cov(Y_{q,0},Y_{r,k})\big|<\infty\) for \(q,r\in\{2,3\}\).
\end{itemize}
Then the leading ratio bias is of order \(O(n^{-1})\).
\end{proposition}

\begin{proof}
Define
\[
    S^\star=m_3m_2^{-3/2}.
\]
Write
$\widehat m_2=m_2+\Delta_2$, $\widehat m_3=m_3+\Delta_3$. The relevant derivatives are $\frac{\partial S^\star}{\partial m_2} = -\frac{3}{2}m_3m_2^{-5/2}$, $\frac{\partial S^\star}{\partial m_3} = m_2^{-3/2}$,
and $\frac{\partial^2 S^\star}{\partial m_2^2} = \frac{15}{4}m_3m_2^{-7/2}$, $\frac{\partial^2 S^\star}{\partial m_2\partial m_3} = -\frac{3}{2}m_2^{-5/2}$, $\frac{\partial^2 S^\star}{\partial m_3^2}=0$. Since \(m_2>0\), the map \((x,y)\mapsto yx^{-3/2}\) is
infinitely differentiable in a neighbourhood of \((m_2,m_3)\). Therefore, a second-order
Taylor expansion around \((m_2,m_3)\) gives
\[
\begin{aligned}
    \widehat S-S^\star
    =
    &\underbrace{-\frac{3}{2}m_3m_2^{-5/2}\Delta_2
    +
    m_2^{-3/2}\Delta_3}_{=:L}
    +
    \underbrace{\frac{15}{8}m_3m_2^{-7/2}\Delta_2^2
    -
    \frac{3}{2}m_2^{-5/2}\Delta_2\Delta_3}_{=:Q}
    +
    R,
\end{aligned}
\]
where, on the event \(\{|\Delta_2|\le m_2/2,\ |\Delta_3|\le 1\}\), the Lagrange remainder
satisfies \(|R|\le C\,(|\Delta_2|+|\Delta_3|)^3\) for a constant \(C<\infty\) bounding the
third derivatives of \((x,y)\mapsto yx^{-3/2}\) on
\([m_2/2,3m_2/2]\times[m_3-1,m_3+1]\).

The leading bias of \(\widehat S\) is the expectation of the linear and quadratic terms,
\(\E[L+Q]\). If \(\widehat m_2\) and \(\widehat m_3\) are unbiased, then
$\E[\Delta_2]=0$, $\E[\Delta_3]=0$, so \(\E[L]=0\). Thus,
\[
    \E[L+Q]
    =
    \frac{15}{8}m_3m_2^{-7/2}\Var(\widehat m_2)
    -
    \frac{3}{2}m_2^{-5/2}\Cov(\widehat m_2,\widehat m_3).
\]
Hence, unless the second-order variance and covariance terms offset each other, the plug-in
skewness estimator is biased in finite samples.

Now suppose \(\widehat m_2\) and \(\widehat m_3\) are computed from \(n\) non-overlapping
low-frequency windows as in (a)--(b), and write \(\gamma_{qr}(k):=\Cov(Y_{q,0},Y_{r,k})\).
Since \(\widehat m_q=n^{-1}\sum_{j=1}^nY_{q,j}\) and, by stationarity,
\(\Cov(Y_{q,j},Y_{r,l})=\gamma_{qr}(l-j)\), grouping the \(n-|k|\) pairs \((j,l)\) with
lag \(l-j=k\) yields
\[
    \Cov(\widehat m_q,\widehat m_r)
    =\frac{1}{n^2}\sum_{j=1}^n\sum_{l=1}^n\gamma_{qr}(l-j)
    =\frac1n\sum_{|k|<n}\Big(1-\frac{|k|}{n}\Big)\gamma_{qr}(k).
\]
By (b),
\(\big|\Cov(\widehat m_q,\widehat m_r)\big|\le n^{-1}\sum_{k\in\mathbb Z}|\gamma_{qr}(k)|\),
so that
\[
    \Var(\widehat m_2)=O(n^{-1}),
    \qquad
    \Cov(\widehat m_2,\widehat m_3)=O(n^{-1}).
\]
By Chebyshev's inequality, $\Delta_2=O_p(n^{-1/2})$ and $\Delta_3=O_p(n^{-1/2})$. In
particular, the event \(\{|\Delta_2|\le m_2/2,\ |\Delta_3|\le 1\}\) has probability tending
to one, and \(R=O_p(n^{-3/2})\) is of smaller order than \(L+Q\). Consequently,
the leading ratio bias is \(O(n^{-1})\).
Therefore the finite-sample ratio bias vanishes as the number of non-overlapping low-frequency observations increases. If \(\widehat m_3\) is itself biased, as for the ACJV numerator, which omits the leverage contribution, then \(\E[L]\neq0\) and an additional bias of order \(O(1)\) remains.
\end{proof}

\subsection{Theoretical properties of the HFN estimator} \label{supp:theory-hfn}
\setcounter{definition}{1}
\begin{definition}[Primitive high-frequency features]\label{def:feat-primitives-new}
Let $\Pi_M=\{t_i\}_{i=0}^M$ be a (possibly irregular) high-frequency partition of $[0,T]$ and define returns $r_i = X_{t_i}-X_{t_{i-1}}$. The primitive observable feature sequence is
\[
    \mathcal R_M := \{(r_i, r_i^2, r_i^3)\}_{i=1}^M .
\]
These channels are not assumed to be sufficient for exact recovery of the latent
low-frequency skewness. Rather, they form an economically motivated representation:
standard realized variance, leverage, and jump-type quantities used in existing
realized skewness constructions can be written as Borel-measurable functionals of
\(\mathcal R_M\).
In particular, we consider:
\begin{itemize}
    \item[(a)] \textbf{Variance-type functionals:} realized variance, semivariances, and (when used) robust realized-kernel variants;
    \item[(b)] \textbf{Leverage-type functionals:} CL/NP covariation terms (e.g.\ $[X,X^2]$) and equivalent representations based on linear filters of $(r_i,r_i^2)$;    
    \item[(c)] \textbf{Jump-type functionals:} stabilized cubic sums such as $J_M=\sum_{i=1}^M r_i^3\mathbf 1\{|r_i|>u_M\}$ or pre-averaged/thresholded cubic variants targeting the jump third cumulant.
\end{itemize}
\end{definition}

\begin{remark}[From primitive features to network inputs]\label{rem:Phi}
The sequence \(\mathcal R_M\) represents the primitive feature representation
constructed from the observed high-frequency path on the grid \(\Pi_M\).
The network does not operate directly on \(\mathcal R_M\), but on a reshaped and
standardized tensor
\[
    \Phi_M:=\phi_M(\mathcal R_M),
\]
where \(\phi_M:(\mathbb R^3)^M\to\mathbb R^{H\times W\times 3}\) is a fixed
Borel-measurable preprocessing map.
If clipping/truncation is applied in preprocessing, then $\Phi_M$ takes values in a fixed compact set $\mathcal X\subset\R^{H\times W\times 3}$.

In particular, $\phi_M$ is a deterministic measurable function of $\mathcal R_M$ and introduces no information beyond the three primitive channels $(r_i,r_i^2,r_i^3)$.
\end{remark}

\subsubsection{Assumptions}   \label{sec: appendix theory-formal: assumptions}

\begin{assumption}[Semimartingale structure and finite moments]\label{ass:sm}
The log-price process $X$ admits the semimartingale decomposition (2) in the main text. Moreover,
\[
    \int_0^T \E[\sigma_t^4]\,dt<\infty,\qquad  \int_0^T\!\int_{\R} |z|^3\,\lambda_t(dz)\,dt<\infty, \qquad  \int_0^T |\mu_t|\,dt<\infty\ \text{a.s.}
\]
In particular, $R_T=X_T-X_0$ has finite second and third (central) moments.
\end{assumption}
For the investigated processes the third central moment of the monthly return $R_T=X_T-X_0$ exists and is finite, see \citet{ORS2022}.  
\begin{assumption}[High-frequency sampling, reshape, and stabilization]\label{ass:phi}
\begin{enumerate}
    \item The observation grid $\Pi_M=\{t_i\}_{i=0}^M$ is non-anticipative and satisfies $\|\Pi_M\|\to0$
    as $M\to\infty$ for a fixed horizon $T$.

    \item The reshape/standardization map $\Phi_M=\phi_M(\mathcal R_M)$ is Borel-measurable. Specifically, $\phi_M$ implements the blocking/reshaping described in Section 2.4.3 (main text) followed by row-wise centering and $\varepsilon$-stabilized scaling: for each row $\ell$ let $\hat\mu_{M,\ell}$ and $\hat\sigma_{M,\ell}$ be its sample mean and standard deviation (computed from the pre-standardization row entries), and define
    \[
        \hat\sigma_{M,\ell}^{\varepsilon}:=\max\{\hat\sigma_{M,\ell},\varepsilon\} \qquad\text{for some fixed }\varepsilon>0,
    \]
    so that each standardized row entry has the form
    \[
        z \ \mapsto\  \frac{z-\hat\mu_{M,\ell}}{\hat\sigma_{M,\ell}^{\varepsilon}}.
    \]

    If clipping/truncation is applied, it is performed componentwise by a fixed 1-Lipschitz map (e.g.\ $u\mapsto \min\{\max\{u,-c\},c\}$ for some $c<\infty$). Consequently, after (optional) clipping the feature tensors $\Phi_M$ take values in a fixed compact set $\mathcal X\subset\R^{H\times W\times 3}$.
    
    Moreover, the standardized rows are uniformly square-integrable under $\P_{\mathrm{tr}}$, i.e.
    \[
        \sup_{M}\ \E_{\P_{\mathrm{tr}}}\big[\|\Phi_M\|_F^2\big]\ <\ \infty.
    \]

    \item Standard realized-variation type functionals computed from $\mathcal R_M$
    (realized variance, semivariances, and, when used, realized-kernel variants) are well-defined Borel-measurable functionals of $\mathcal R_M$ and satisfy their usual high-frequency consistency statements for the corresponding quadratic-(co)variation targets under the semimartingale model. When microstructure-noise robustness is required, the corresponding robust variants are used.

    \item (Jump functional / stabilization).
    There exists a sequence $u_M\downarrow 0$ such that the thresholded cubic functional $J_M := \sum_{i=1}^{M} r_i^3\,\mathbf 1\{|r_i|>u_M\}$ is well-defined and satisfies
    $J_M \xrightarrow{P} \sum_{0<s\le T} (\Delta X_s)^3$ as $\|\Pi_M\|\to 0$ under $\P_{\mathrm{tr}}$. A sufficient rate condition is $u_M\to0$ and $\sqrt{M}\,u_M\to\infty$. When microstructure-noise robustness is required, a pre-averaging plus threshold variant of $J_M$ is used.    
\end{enumerate}
\end{assumption}

\begin{assumption}[Network class: measurability, boundedness, Lipschitz continuity]\label{ass:cnn}
    For each $p$, the hypothesis class $\mathcal F_p$ consists of measurable predictors $f_{\mathrm{sk}}:\mathcal X\to\R$ (skewness head). Moreover, there exist constants $L_p<\infty$ such that for all $x,x'\in\mathcal X$,
    \[
        |f_{\mathrm{sk}}(x)-f_{\mathrm{sk}}(x')|\le L_p\,\|x-x'\|,
    \]
    where $\|\cdot\|$ denotes the Euclidean norm after vectorization (equivalently, the Frobenius norm on $\R^{H\times W\times 3}$). Since $\mathcal X$ is compact, each $f_{\mathrm{sk}}$ is automatically bounded on $\mathcal X$; we denote $B_p:=\sup_{x\in\mathcal X}|f_{\mathrm{sk}}(x)|<\infty$. In the HFN architecture, Lipschitz continuity is ensured by finite operator norms of convolutional/linear layers together with $1$-Lipschitz activations (e.g.\ ReLU, $\tanh$) and non-expansive pooling.
\end{assumption}

\begin{remark}[Lipschitz property of the squared loss]\label{rem:lipschitz-loss}
Let $\mathcal Z:=\mathcal X\times[-S_{\max},S_{\max}]$ equipped with
$d\big((x,s),(x',s')\big)=\|x-x'\|+|s-s'|$.
For $f\in\mathcal F_p$ define $v_f(x,s):=(f_{\mathrm{sk}}(x)-s)^2$.
Under Assumption~\ref{ass:cnn}, for all $(x,s),(x',s')\in\mathcal Z$,
\begin{align*}
|v_f(x,s)-v_f(x',s')|
&\le 2(B_p+S_{\max})\Big(L_p\|x-x'\|+|s-s'|\Big) \\
&\le 2(B_p+S_{\max})\max\{L_p,1\}\, d\big((x,s),(x',s')\big).
\end{align*}
Hence $\operatorname{Lip}(v_f)\le 2(B_p+S_{\max})\max\{L_p,1\}$ on $(\mathcal Z,d)$.
\end{remark}

We will use the induced Lipschitz property of the squared loss in the robustness bound, see Section \ref{sec:distDrift}.

\begin{assumption}[Capacity control]\label{ass:comp}
    Let $\mathfrak{R}_N(\mathcal F_p)$ denote the empirical Rademacher complexity of $\{x\mapsto f_{\mathrm{sk}}(x):f\in\mathcal F_p\}$ on $\mathcal X$. There exists a schedule $p=p_N\uparrow\infty$ such that
    \[
        \E\big[\mathfrak{R}_N(\mathcal F_{p_N})\big]\to0 \qquad \text{as } N\to\infty,
    \]
    where the expectation is taken under the training law $\P_{\mathrm{tr}}$.
\end{assumption}

\begin{remark}
    In practice, the capacity schedule is implemented by controlling architecture size and by standard regularization (e.g., weight decay), but the theoretical results below only require the complexity condition in Assumption~\ref{ass:comp}.
\end{remark}

This assumption ensures uniform convergence of empirical to population risk under the standard bounded-Lipschitz framework, \citet{BartlettMendelson2002}. 

\begin{assumption}[Regularity of the Bayes predictor]\label{ass:bayes}
    Fix $M$ (hence the input tensor $\Phi_M$ and compact domain $\mathcal X$). The Bayes regression function
    \[
        \tilde h_M(x):=\E[S^\star\mid \Phi_M=x],\qquad x\in\mathcal X,
    \]
    admits a continuous version on $\mathcal X$, i.e.\ $\tilde h_M\in C(\mathcal X)$.
\end{assumption}

This regularity is mild under the simulation designs used here (smooth dependence of $S^\star$ on the model parameters and continuous dependence of $\Phi_M$ on these parameters), and is standard when invoking uniform approximation on compacta.

\begin{assumption}[Non-degeneracy of the skewness label]\label{ass:decomp}
Under \(\P_{\mathrm{tr}}\), the second central moment of the low-frequency return
satisfies
\[
    m_2(T):=\Var(R_T)\ge v_{\min}>0.
\]
Hence the skewness label
\[
    S^\star=\frac{m_3(T)}{m_2(T)^{3/2}}
\]
is well-defined. Moreover, \(S^\star\in L_2(\P_{\mathrm{tr}})\), and, for the bounded-loss
arguments below, we assume \(|S^\star|\le S_{\max}\) almost surely.
\end{assumption}

\begin{assumption}[Sampling of training pairs]\label{ass:iid-train}
    Training samples $\{(\Phi_M^{(j)},S^{\star(j)})\}_{j=1}^N$ are i.i.d.\ from $\P_{\mathrm{tr}}$. For a fixed $M$, each feature tensor $\Phi_M^{(j)}$ is constructed using the same grid size (mesh) $\|\Pi_M\|$ and the same reshape/standardization rule. The target is bounded: $|S^\star|\le S_{\max}$.
\end{assumption}

\subsubsection{Primitive-feature representation and Bayes prediction target}
\begin{remark}[Why plain $\sum_i r_i^3$ is not enough]
    For a semimartingale $X$, the na\"ive cubic power sum satisfies
\[
    \sum_{i=1}^{M} r_i^3 \xrightarrow{P} \sum_{0<s\le T} (\Delta X_s)^3 \quad\text{as }\|\Pi_M\|\to0,
\]
    so it captures only the jump third variation and omits the diffusive leverage contribution (cf.\ Eq. (8) in the main text). Therefore, recovering the third moment requires additional leverage-type functionals (CL/NP), and jump contributions may be isolated via thresholded (stabilized) cubic sums.
\end{remark}

The results of the following Lemma follow directly from \citet{choelee} for the 2nd part and from \citet{NP2021} for the 3rd part.
 
\begin{lemma}[Measurability of classical realized skewness proxies]\label{lem:proxy-meas}
Let \(\mathcal R_M=\{(r_i,r_i^2,r_i^3)\}_{i=1}^M\). The ACJV cubic variation,
the CL realized covariation term, and the NP filtered covariance term can be written
as Borel-measurable functionals of \(\mathcal R_M\), provided the initial value \(X_0\)
is fixed or observed.

In particular,
\[
    \sum_{i=1}^M r_i^3
\]
is directly measurable from \(\mathcal R_M\). Moreover, since
\[
    X_{t_{i-1}}=X_0+\sum_{j=1}^{i-1}r_j,
\]
the CL term
\[
    \sum_{i=1}^M r_i\bigl(X_{t_i}^2-X_{t_{i-1}}^2\bigr)
\]
is also a measurable functional of \(\mathcal R_M\). Finally, the NP-type filtered term
\[
    \sum_{i=1}^M y_i^{(1\ast)}r_i^2
\]
is measurable whenever \(y_i^{(1\ast)}\) is constructed as a deterministic linear filter of past returns.
\end{lemma}
\begin{proof}
The ACJV term is immediate. For the CL term, the price level \(X_{t_{i-1}}\) is
reconstructed from cumulative returns when \(X_0\) is fixed or observed, so
\(X_{t_i}^2-X_{t_{i-1}}^2\) is a deterministic function of the return sequence.
For the NP term, \(y_i^{(1\ast)}\) is a deterministic linear filter of past returns,
and multiplication by \(r_i^2\) preserves measurability. Hence all three objects are
Borel-measurable functionals of \(\mathcal R_M\).
\end{proof}

\begin{remark}[Target as a random label under $\P_{\mathrm{tr}}$]
    Under $\P_{\mathrm{tr}}$, the label $S^\star$ is random (through the parameter draw). For fixed $M$, the Bayes regression function $\tilde h_M(x):=\E[S^\star\mid \Phi_M=x]$ is well-defined on $\mathcal X$. The following results justify approximation of $\tilde h_M$ from the primitive features.
\end{remark}

\setcounter{proposition}{1}
\begin{proposition}[Optimal prediction of low-frequency skewness using high-frequency features]
\label{prop:motivation-measurable}
    Let $\mathcal R_M=\{(r_i,r_i^2,r_i^3)\}_{i=1}^M$ and assume $\Phi_M=\phi_M(\mathcal R_M)$ for a Borel-measurable map $\phi_M$. Assume $S^\star\in L_2(\P_{\mathrm{tr}})$. Then $h_M(\mathcal R_M):=\E[S^\star\mid \mathcal R_M]$ is well-defined and $\sigma(\mathcal R_M)$-measurable. Moreover, since $\Phi_M$ is $\sigma(\mathcal R_M)$-measurable, we have $\sigma(\Phi_M)\subseteq \sigma(\mathcal R_M)$ and there exists a measurable $\tilde h_M:\mathcal X\to\R$ such that
    \[
        \E[S^\star\mid \Phi_M] = \tilde h_M(\Phi_M)\quad\text{a.s.}
    \]
    This conditional
    expectation is the best predictor under squared loss within the information set
    generated by the chosen features.
\end{proposition}

\begin{proof}
Since $S^\star\in L_2(\P_{\mathrm{tr}})$, we have $S^\star\in L_1(\P_{\mathrm{tr}})$, so the conditional expectation $h_M(\mathcal R_M):=\E[S^\star\mid \mathcal R_M]$ exists and is $\sigma(\mathcal R_M)$-measurable by definition. Because $\Phi_M=\phi_M(\mathcal R_M)$ with $\phi_M$ Borel-measurable, $\Phi_M$ is $\sigma(\mathcal R_M)$-measurable, hence $\sigma(\Phi_M)\subseteq \sigma(\mathcal R_M)$. Moreover, $\E[S^\star\mid \Phi_M]$ is $\sigma(\Phi_M)$-measurable, and by the Doob--Dynkin lemma there exists a measurable map $\tilde h_M:\mathcal X\to\R$ such that $\E[S^\star\mid \Phi_M]=\tilde h_M(\Phi_M)$ almost surely.

For any $\sigma(\mathcal R_M)$-measurable predictor $g(\mathcal R_M)$ with $\E[g(\mathcal R_M)^2]<\infty$,
\[
\E\big[(S^\star-g(\mathcal R_M))^2\mid \mathcal R_M\big]
= \Var(S^\star\mid \mathcal R_M) + \big(\E[S^\star\mid \mathcal R_M]-g(\mathcal R_M)\big)^2.
\]
Taking expectations yields
\[
\E[(S^\star-g(\mathcal R_M))^2]
= \E[\Var(S^\star\mid \mathcal R_M)] + \E\big[\big(h_M(\mathcal R_M)-g(\mathcal R_M)\big)^2\big],
\]
which is minimized uniquely by $g=h_M$. The statement for conditioning on $\Phi_M$ follows analogously. 
\end{proof}

Proposition \ref{prop:motivation-measurable} formalizes the supervised learning target associated
with the chosen feature representation. Since \(\Phi_M\) is a measurable
transformation of \(\mathcal R_M\), the optimal squared-loss predictor based on the
network input is the conditional expectation
\[
    \tilde h_M(\Phi_M)=\E[S^\star\mid \Phi_M].
\]
This is a prediction target induced by the feature representation.



\begin{remark}[Implications for the HFN architecture]\label{rem:arch}
    Proposition \ref{prop:motivation-measurable} identifies the learning target under squared loss as the Bayes regression function
    \[
        h_M(\mathcal R_M)=\E[S^\star\mid \mathcal R_M], \qquad\text{equivalently}\qquad \E[S^\star\mid \Phi_M]=\tilde h_M(\Phi_M).
    \]
    Lemma~\ref{lem:proxy-meas} illustrates that classical high-frequency leverage and jump functionals used in the literature are measurable functionals of the primitive sequence $\mathcal R_M$. Therefore, the HFN operates on a feature representation that contains the same
    primitive ingredients used by classical realized skewness proxies. The learning
    problem is to approximate the Bayes predictor
    \[
        \tilde h_M(x)=\E[S^\star\mid \Phi_M=x].
    \]
    Approximation of this conditional expectation is addressed by the risk-consistency
    result below.
\end{remark}

\subsubsection{Approximation and risk consistency (training asymptotics)}
The HFN maps the fixed feature tensor
\(\Phi_M\in\mathcal X\subset\mathbb R^{H\times W\times 3}\) into a skewness prediction.
As introduced in Assumption \ref{ass:cnn}, for \(f\in\mathcal F_p\), let \(f_{\mathrm{sk}}\) denote the skewness output coordinate.
The learning target induced by squared loss is the Bayes regression function
\[
    \tilde h_M(x)=\E[S^\star\mid \Phi_M=x],
    \qquad x\in\mathcal X.
\]
Under Assumption~\ref{ass:bayes}, this function admits a continuous version on the
compact feature domain \(\mathcal X\). The role of the approximation result below is
therefore not to establish that the HFN identifies \(S^\star\) exactly, but to state
that the HFN family can approximate the Bayes
prediction target associated with the chosen feature representation.

\begin{proposition}[Universal approximation of the HFN family]
\label{prop:approx}
Under Assumptions~\ref{ass:phi}\textup{(ii)},
\ref{ass:cnn}, and~\ref{ass:bayes}, let $\mathcal F_p$ denote the HFN family introduced in
Section 2.4.4 of the main text. The induced function class
\[
    \bigcup_{p\ge1}\{x\mapsto f_{\mathrm{sk}}(x):f\in\mathcal F_p\}
\]
is dense in \(C(\mathcal X)\) under the uniform norm. Then, for every \(\varepsilon>0\),
there exist \(p\) and \(f\in\mathcal F_p\) such that
\[
    \sup_{x\in\mathcal X}
    \left|f_{\mathrm{sk}}(x)-\tilde h_M(x)\right|
    \le \varepsilon .
\]
In particular,
\[
    \E_{\P_{\mathrm{tr}}}
    \left[
        \left(f_{\mathrm{sk}}(\Phi_M)-\tilde h_M(\Phi_M)\right)^2
    \right]
    \le \varepsilon^2 .
\]
\end{proposition}

\begin{proof}
    Since $\mathcal X$ is compact (Assumption~\ref{ass:phi}\textup{(ii)}) and $\tilde h_M\in C(\mathcal X)$ (Assumption~\ref{ass:bayes}), the claimed uniform approximation follows from a universal-approximation theorem for the architecture class $\bigcup_p\mathcal F_p$. Classical results \citep{Cybenko1989,Hornik1991,Leshno1993} yield universality for shallow networks with non-polynomial activations, and modern results \citep{Yarotsky2017,Zhou2020} provide analogous approximation properties for convolutional architectures under standard constraints. The $L_2(\P_{\mathrm{tr}})$ bound follows immediately since $\Phi_M\in\mathcal X$ almost surely:
    \[
        |f_{\mathrm{sk}}(\Phi_M)-\tilde h_M(\Phi_M)|\le \sup_{x\in\mathcal X}|f_{\mathrm{sk}}(x)-\tilde h_M(x)|\le \varepsilon.
    \]
\end{proof}

\begin{remark}[Interpretation]
The \(L_2(\P_{\mathrm{tr}})\) bound statement translates the universal approximation property into the statistical learning. Proposition~\ref{prop:approx} should be interpreted as an approximation
statement for the Bayes predictor generated by the fixed feature map \(\phi_M\).
\end{remark}

Building on the universal approximation result above, we now establish risk consistency of the empirical-risk minimizer under bounded-Lipschitz conditions.
\setcounter{theorem}{0}
\begin{theorem}[Risk consistency of ERM for fixed $M$]\label{lem:risk} 
    Fix $M$ (hence $\Phi_M$ and $\mathcal X$). Let $\hat f_{p,N}$ be a (regularized) empirical risk minimizer over $\mathcal F_p$:
    \[
        \widehat{\mathcal R}_N(f):= \frac{1}{N}\sum_{i=1}^{N} \big(f_{\mathrm{sk}}(\Phi_M^{(i)}) - S^{\star(i)}\big)^2 + \lambda_p\|f\|_{\mathrm{wd}}^2,
    \]
    with $\lambda_p\to0$ and $p=p_N$ satisfying Assumption~\ref{ass:comp}. Under Assumptions~\ref{ass:iid-train} and~\ref{ass:cnn}, as $N\to\infty$,
    \[
        \mathcal R_M(\hat f_{p_N,N}) := \E_{\P_{\mathrm{tr}}}\!\big[(\hat f_{p_N,N}(\Phi_M) - S^\star)^2\big] \ \longrightarrow\ \inf_{f\in\mathcal F_{p_N}}\mathcal R_M(f).
    \]
    If, in addition, $\bigcup_p \mathcal F_p$ is dense in $C(\mathcal X)$ as in Proposition~\ref{prop:approx}, then
    \[
        \lim_{N\to\infty}\mathcal R_M(\hat f_{p_N,N}) =\inf_{f\in\cup_p\mathcal F_p}\mathcal R_M(f) =\mathcal R_M(\tilde h_M) =\E_{\P_{\mathrm{tr}}}\!\big[(\tilde h_M(\Phi_M) - S^\star)^2\big].
    \]
\end{theorem}

\begin{proof}
    Define the loss $\ell_f(x,s):=(f_{\mathrm{sk}}(x)-s)^2$. Because $|f_{\mathrm{sk}}(x)|\le B_{p_N}$ and $|S^\star|\le S_{\max}$, the loss class $\mathcal L_{p_N}:=\{\ell_f:f\in\mathcal F_{p_N}\}$ is uniformly bounded:
    \[
        0\le \ell_f(x,s)\le (B_{p_N}+S_{\max})^2.
    \]
    The map $\phi(u)=u^2$ is $L_\phi=2(B_{p_N}+S_{\max})$-Lipschitz on this range. By the Ledoux-Talagrand contraction inequality (\citealp[Thm.~12(4)]{BartlettMendelson2002}),
    \[
        \mathfrak R_N(\mathcal L_{p_N})\ \le\ 4(B_{p_N}+S_{\max})\,\mathfrak R_N(\mathcal F_{p_N}),
    \]
    so the loss-class complexity is controlled by that of the predictor class.  Standard Rademacher bounds (\citealp[Thms.~8 and~12]{BartlettMendelson2002}) then yield, with probability at least $1-\delta$,
    \[
        \sup_{f\in\mathcal F_{p_N}}\big|\widehat{\mathcal R}_N(f)-\mathcal R_M(f)\big|\ \le\ 2\,\mathfrak R_N(\mathcal L_{p_N}) + c\sqrt{\tfrac{\log(1/\delta)}{2N}} = O\!\big(\mathfrak R_N(\mathcal F_{p_N})\big).
    \]
    By Assumption~\ref{ass:comp}, $\mathfrak R_N(\mathcal F_{p_N})\to0$, hence
    \[
        \sup_{f\in\mathcal F_{p_N}}\big|\widehat{\mathcal R}_N(f)-\mathcal R_M(f)\big|\xrightarrow[]{P}0.
    \]
    Let $f^\star_{p_N}\in\arg\min_{f\in\mathcal F_{p_N}}\mathcal R_M(f)$. Approximate empirical optimality of $\hat f_{p_N,N}$ and the above deviation bound imply
    \[
        \mathcal R_M(\hat f_{p_N,N}) - \mathcal R_M(f^\star_{p_N}) \ \le\ 2\sup_{f\in\mathcal F_{p_N}}\big|\widehat{\mathcal R}_N(f)-\mathcal R_M(f)\big| + o(1)\ \xrightarrow[]{} 0,
    \]
    so $\mathcal R_M(\hat f_{p_N,N})\to\inf_{f\in\mathcal F_{p_N}}\mathcal R_M(f)$. 
    
    If $\mathcal F_{p_N}$ is dense in $C(\mathcal X)$ (Proposition~\ref{prop:approx}), then $\cup_p \mathcal F_p$ is dense in $L_2(\P_{\Phi_M})$. Since the Bayes predictor under squared loss for inputs $\Phi_M$ is
    \[
        \tilde h_M(\Phi_M)=\E[S^\star\mid \Phi_M] \quad\text{(Assumption~\ref{ass:bayes})},
    \]
    we obtain
    \[
        \inf_{f\in\cup_p\mathcal F_p}\mathcal R_M(f) = \E_{\P_{\mathrm{tr}}}\!\big[(\tilde h_M(\Phi_M) - S^\star)^2\big].
    \]
    This identifies the limiting risk of $\hat f_{p_N,N}$ with the Bayes risk for the feature input $\Phi_M$.

\end{proof}

\begin{remark}[Fixed sampling design]
All risk-consistency statements are for a fixed \(M\), corresponding to the fixed
sampling frequency and fixed input tensor.
Extending the result to a joint limit \(M,N\to\infty\) would require additional
conditions ensuring stability of the conditional expectations
\(\E[S^\star\mid\Phi_M]\), which we do not impose here.
\end{remark}

\begin{remark}[What is and is not implied by Theorem~\ref{lem:risk}]
    Theorem~\ref{lem:risk} establishes risk consistency for a \emph{fixed} feature input $\Phi_M$: for suitable $p=p_N$, the ERM $\hat f_{p_N,N}$ achieves asymptotically the best-in-class risk $\inf_{f\in\mathcal F_{p_N}}\mathcal R_M(f)$, and under density it attains the Bayes risk $\mathcal R_M(\tilde h_M)$ for the regression function $\tilde h_M(x)=\E[S^\star\mid \Phi_M=x]$. It makes no statement about the behavior of the HFN under high-frequency asymptotics $M\to\infty$.
\end{remark}

\subsubsection{Robustness to distribution shift}\label{sec:distDrift}
Let $\P_0$ denote the true law of $(\Phi_M,S^\star)$ and $\P_{\mathrm{tr}}$ the training law. Equip $\mathcal Z := \mathcal X \times [-S_{\max}, S_{\max}]$ with the product metric
\[
    d\big((x,s),(x',s')\big) := \|x-x'\|_2 + |s-s'|.
\]
Let $\mathcal P(\mathcal Z)$ denote the set of Borel probability measures on $(\mathcal Z,\mathcal B(\mathcal Z))$ with finite first moment, where $\mathcal B(\mathcal Z)$ is the Borel $\sigma$-algebra induced by $d$, and let $W_1$ denote the 1-Wasserstein distance on $\mathcal P(\mathcal Z)$. For a function $g:\mathcal Z \to \mathbb R$, define its Lipschitz seminorm by
\[
    \operatorname{Lip}(g):= \sup_{\substack{z,z'\in\mathcal Z\\ z\neq z'}}\frac{|g(z)-g(z')|}{d(z,z')}.
\]
\setcounter{remark}{1}
\begin{remark}[Distributional stability of the prediction risk]\label{thm:robustrisk}
    Under Assumption~\ref{ass:cnn}, with $|S^\star|\le S_{\max}$ a.s.\ and $\mathcal X$ compact, for any $f\in\mathcal F_p$ we have
    \[
        \Big|\E_{\P_0}\!\big[(f_{\mathrm{sk}}(\Phi_M) - S^\star)^2\big] - \E_{\P_{\mathrm{tr}}}\!\big[(f_{\mathrm{sk}}(\Phi_M) - S^\star)^2\big]\Big|\ \le\ C(L_p,B_p,S_{\max})\, W_1(\P_0,\P_{\mathrm{tr}}),
    \]
    for a constant $C$ depending only on the Lipschitz constant $L_p$ and the uniform bounds $B_p,S_{\max}$. This result is a standard bounded-Lipschitz stability bound. It is not
specific to the HFN architecture beyond the boundedness and Lipschitz restrictions
imposed on the predictor class. Its role is to clarify how prediction risk changes
under distributional perturbations of the joint law of features and labels.
\end{remark}

\begin{proof} 
    Define $v(x,s):=(f_{\mathrm{sk}}(x)-s)^2$. On $\mathcal Z=\mathcal X\times[-S_{\max},S_{\max}]$, $v$ is Lipschitz because for any $(x,s),(x',s')\in\mathcal Z$,
    \begin{align*}
        |v(x,s)-v(x',s')|&\le 2|f_{\mathrm{sk}}(x)-s|\,\big(|f_{\mathrm{sk}}(x)-f_{\mathrm{sk}}(x')|+|s-s'|\big)\\
        &\le 2(B_p+S_{\max})(L_p\|x-x'\|_2+|s-s'|).
    \end{align*}
    Hence 
    \[
        \operatorname{Lip}(v)\ \le\ 2(B_p+S_{\max})\,\max\{L_p,1\}\ :=\ C(L_p,B_p,S_{\max}).
    \]
    By Kantorovich–Rubinstein duality,
    \[
        \big|\E_{\P_0}[v]-\E_{\P_{\mathrm{tr}}}[v]\big|\ \le\ \operatorname{Lip}(v)\,W_1(\P_0,\P_{\mathrm{tr}}),
    \]
    which yields the stated inequality.
    \end{proof}

\begin{corollary}[Misspecification bound for ERM]\label{cor:misspec}
    If $\hat f_{p,N}$ is $\eta_N$-optimal under $\P_{\mathrm{tr}}$ (i.e.\ its risk satisfies $\E_{\P_{\mathrm{tr}}}[(\hat f_{p,N}(\Phi_M)-S^\star)^2]\le\inf_{f\in\mathcal F_p}\E_{\P_{\mathrm{tr}}}[(f(\Phi_M)-S^\star)^2]+\eta_N$), then under the shifted distribution $\P_0$,
    \begin{align*}
    &\E_{\P_0}\!\big[(\hat    f_{p,N}(\Phi_M) - S^\star)^2\big]\ \le\ \inf_{f\in\mathcal F_p}\E_{\P_0}\!\big[(f(\Phi_M) - S^\star)^2\big]\ +\\
    &\qquad\qquad\ \eta_N\ + \ C(L_p,B_p,S_{\max})\,W_1(\P_0,\P_{\mathrm{tr}}).
    \end{align*}
\end{corollary}
\setcounter{remark}{11}
\begin{remark}[Interpretation]
    The bound quantifies first-order stability of the learned predictor under mild covariate–target distribution shift: the excess risk changes at most linearly with the Wasserstein distance between training and testing laws. This bound provides a formal stability interpretation: for predictors with controlled
    Lipschitz constants and bounded outputs, small shifts in the joint feature-label law
    lead to small changes in squared prediction risk. The empirical misspecification
    experiments in Section 3.2 of the main text and Section~\ref{supp:misspec} examine whether such stability is observed in the simulated
    settings considered here.
\end{remark}

\bibliographystyle{elsarticle-harv} 
{\small
\renewcommand{\baselinestretch}{1} 
\bibliography{literature}
}

\renewcommand{\baselinestretch}{1.5} 
\newpage
\appendix
\renewcommand{\thesubsection}{\Alph{section}.\roman{subsection}}